\input amstex

%
   \def\tableline{\hbox to \hsize}
%
%
\newbox\hdbox%
\newcount\hdrows%
\newcount\multispancount%
\newcount\ncase%
\newcount\ncols
\newcount\nrows%
\newcount\nspan%
\newcount\ntemp%
\newdimen\hdsize%
\newdimen\newhdsize%
\newdimen\parasize%
\newdimen\spreadwidth%
\newdimen\thicksize%
\newdimen\thinsize%
\newdimen\tablewidth%
\newif\ifcentertables%
\newif\ifendsize%
\newif\iffirstrow%
\newif\iftableinfo%
\newtoks\dbt%
\newtoks\hdtks%
\newtoks\savetks%
\newtoks\tableLETtokens%
\newtoks\tabletokens%
\newtoks\widthspec%
%
%
%
%
\tableinfotrue%
\catcode`\@=11
%
%
\def\tstrut{\vrule height3.1ex depth1.2ex width0pt}%
\def\and{\char`\&}
\def\tablerule{\noalign{\hrule height\thinsize depth0pt}}%
\thicksize=1.5pt
\thinsize=0.6pt
\def\thickrule{\noalign{\hrule height\thicksize depth0pt}}%
\def\ctr#1{\hfil\ #1\hfil}%
%
%
%
%
\tablewidth=-\maxdimen%
\spreadwidth=-\maxdimen%
\def\tabskipglue{0pt plus 1fil minus 1fil}%
%
%
\centertablestrue%
%
%
%
%
\parasize=4in%
\gdef\ARGS{########}
\gdef\headerARGS{####}
\def\@mpersand{&}
{\catcode`\|=13
\gdef\letbarzero{\let|0}
\gdef\letbartab{\def|{&&}}%
\gdef\letvbbar{\let\vb|}%
}
{\catcode`\&=4
\def\ampskip{&\omit\hfil&}
\catcode`\&=13
\let&0
\xdef\letampskip{\def&{\ampskip}}%
\gdef\letnovbamp{\let\novb&\let\tab&}
}
\def\begintable{
   \begingroup%
   \catcode`\|=13\letbartab\letvbbar%
   \catcode`\&=13\letampskip\letnovbamp%
   \def\multispan##1{
      \omit \mscount##1%
      \multiply\mscount\tw@\advance\mscount\m@ne%
      \loop\ifnum\mscount>\@ne \sp@n\repeat%
   }
   \def\|{%
      &\omit\widevline&%
   }%
   \ruledtable
}
\long\def\ruledtable#1\endtable{%
%
%
%
   \offinterlineskip
   \tabskip 0pt
   \def\widevline{\vrule width\thicksize}
   \def\endrow{\@mpersand\omit\hfil\crnorm\@mpersand}%
   \def\crthick{\@mpersand\crnorm\thickrule\@mpersand}%
   \def\crnorule{\@mpersand\crnorm\@mpersand}%
   \let\nr=\crnorule
   \def\endtable{\@mpersand\crnorm\thickrule}%
   \let\crnorm=\cr
%
%
   \edef\cr{\@mpersand\crnorm\tablerule\@mpersand}%
   \the\tableLETtokens
%
%
   \tabletokens={&#1}
%
%
   \countROWS\tabletokens\into\nrows%
   \countCOLS\tabletokens\into\ncols%
%
%
   \advance\ncols by -1%
   \divide\ncols by 2%
   \advance\nrows by 1%
%
%
   \iftableinfo %
      \immediate\write16{[Nrows=\the\nrows, Ncols=\the\ncols]}%
   \fi%
%
%
   \ifcentertables
      \ifhmode \par\fi
      \tableline{
      \hss
   \else %
      \hbox{%
   \fi
      \vbox{%
         \makePREAMBLE{\the\ncols}
         \edef\next{\preamble}
         \let\preamble=\next
         \makeTABLE{\preamble}{\tabletokens}
      }
      \ifcentertables \hss}\else }\fi
   \endgroup
   \tablewidth=-\maxdimen
   \spreadwidth=-\maxdimen
}
\def\makeTABLE#1#2{
   {
   \let\ifmath0
   \let\header0
   \let\multispan0
%
%
   \ncase=0%
   \ifdim\tablewidth>-\maxdimen \ncase=1\fi%
   \ifdim\spreadwidth>-\maxdimen \ncase=2\fi%
   \relax
%
   \ifcase\ncase %
      \widthspec={}%
   \or %
      \widthspec=\expandafter{\expandafter t\expandafter o%
                 \the\tablewidth}%
   \else %
      \widthspec=\expandafter{\expandafter s\expandafter p\expandafter r%
                 \expandafter e\expandafter a\expandafter d%
                 \the\spreadwidth}%
   \fi %
   \xdef\next{
      \halign\the\widthspec{%
      #1
      \noalign{\hrule height\thicksize depth0pt}
      \the#2\endtable
%
      }
   }
   }
   \next
}
\def\makePREAMBLE#1{
   \ncols=#1
   \begingroup
   \let\ARGS=0
   \edef\xtp{\widevline\ARGS\tabskip\tabskipglue%
   &\ctr{\ARGS}\tstrut}
   \advance\ncols by -1
   \loop
      \ifnum\ncols>0 %
      \advance\ncols by -1%
      \edef\xtp{\xtp&\vrule width\thinsize\ARGS&\ctr{\ARGS}}%
   \repeat
   \xdef\preamble{\xtp&\widevline\ARGS\tabskip0pt%
   \crnorm}
   \endgroup
}
\def\countROWS#1\into#2{
   \let\countREGISTER=#2%
   \countREGISTER=0%
   \expandafter\ROWcount\the#1\endcount%
}%
\def\ROWcount{%
   \afterassignment\subROWcount\let\next= %
}%
\def\subROWcount{%
   \ifx\next\endcount %
      \let\next=\relax%
   \else%
      \ncase=0%
      \ifx\next\cr %
         \global\advance\countREGISTER by 1%
         \ncase=0%
      \fi%
      \ifx\next\endrow %
         \global\advance\countREGISTER by 1%
         \ncase=0%
      \fi%
      \ifx\next\crthick %
         \global\advance\countREGISTER by 1%
         \ncase=0%
      \fi%
      \ifx\next\crnorule %
         \global\advance\countREGISTER by 1%
         \ncase=0%
      \fi%
      \ifx\next\header %
         \ncase=1%
      \fi%
      \relax%
      \ifcase\ncase %
         \let\next\ROWcount%
      \or %
         \let\next\argROWskip%
      \else %
      \fi%
   \fi%
   \next%
}
\def\counthdROWS#1\into#2{%
\dvr{10}%
   \let\countREGISTER=#2%
   \countREGISTER=0%
\dvr{11}%
\dvr{13}%
   \expandafter\hdROWcount\the#1\endcount%
\dvr{12}%
}%
\def\hdROWcount{%
   \afterassignment\subhdROWcount\let\next= %
}%
\def\subhdROWcount{%
   \ifx\next\endcount %
      \let\next=\relax%
   \else%
      \ncase=0%
      \ifx\next\cr %
         \global\advance\countREGISTER by 1%
         \ncase=0%
      \fi%
      \ifx\next\endrow %
         \global\advance\countREGISTER by 1%
         \ncase=0%
      \fi%
      \ifx\next\crthick %
         \global\advance\countREGISTER by 1%
         \ncase=0%
      \fi%
      \ifx\next\crnorule %
         \global\advance\countREGISTER by 1%
         \ncase=0%
      \fi%
      \ifx\next\header %
         \ncase=1%
      \fi%
\relax%
      \ifcase\ncase %
         \let\next\hdROWcount%
      \or%
         \let\next\arghdROWskip%
      \else %
      \fi%
   \fi%
   \next%
}%
{\catcode`\|=13\letbartab
\gdef\countCOLS#1\into#2{%
   \let\countREGISTER=#2%
   \global\countREGISTER=0%
   \global\multispancount=0%
   \global\firstrowtrue
   \expandafter\COLcount\the#1\endcount%
   \global\advance\countREGISTER by 3%
   \global\advance\countREGISTER by -\multispancount
}%
\gdef\COLcount{%
   \afterassignment\subCOLcount\let\next= %
}%
{\catcode`\&=13%
\gdef\subCOLcount{%
   \ifx\next\endcount %
      \let\next=\relax%
   \else%
      \ncase=0%
      \iffirstrow
         \ifx\next& %
            \global\advance\countREGISTER by 2%
            \ncase=0%
         \fi%
         \ifx\next\span %
            \global\advance\countREGISTER by 1%
            \ncase=0%
         \fi%
         \ifx\next| %
            \global\advance\countREGISTER by 2%
            \ncase=0%
         \fi
         \ifx\next\|
            \global\advance\countREGISTER by 2%
            \ncase=0%
         \fi
         \ifx\next\multispan
            \ncase=1%
            \global\advance\multispancount by 1%
         \fi
         \ifx\next\header
            \ncase=2%
         \fi
         \ifx\next\cr       \global\firstrowfalse \fi
         \ifx\next\endrow   \global\firstrowfalse \fi
         \ifx\next\crthick  \global\firstrowfalse \fi
         \ifx\next\crnorule \global\firstrowfalse \fi
      \fi
\relax
      \ifcase\ncase %
         \let\next\COLcount%
      \or %
         \let\next\spancount%
      \or %
         \let\next\argCOLskip%
      \else %
      \fi %
   \fi%
   \next%
}%
\gdef\argROWskip#1{%
   \let\next\ROWcount \next%
}
\gdef\arghdROWskip#1{%
   \let\next\ROWcount \next%
}
\gdef\argCOLskip#1{%
   \let\next\COLcount \next%
}
}
}
\def\spancount#1{
   \nspan=#1\multiply\nspan by 2\advance\nspan by -1%
   \global\advance \countREGISTER by \nspan
   \let\next\COLcount \next}%
\def\dvr#1{\relax}%
\def\header#1{%
\dvr{1}{\let\cr=\@mpersand%
\hdtks={#1}%
\counthdROWS\hdtks\into\hdrows%
\advance\hdrows by 1%
\ifnum\hdrows=0 \hdrows=1 \fi%
\dvr{5}\makehdPREAMBLE{\the\hdrows}%
\dvr{6}\getHDdimen{#1}%
{\parindent=0pt\hsize=\hdsize{\let\ifmath0%
\xdef\next{\valign{\headerpreamble #1\crnorm}}}\dvr{7}\next\dvr{8}%
}%
}\dvr{2}}
\def\makehdPREAMBLE#1{
\dvr{3}%
\hdrows=#1
{
\let\headerARGS=0%
\let\cr=\crnorm%
\edef\xtp{\vfil\hfil\hbox{\headerARGS}\hfil\vfil}%
\advance\hdrows by -1
\loop
\ifnum\hdrows>0%
\advance\hdrows by -1%
\edef\xtp{\xtp&\vfil\hfil\hbox{\headerARGS}\hfil\vfil}%
\repeat%
\xdef\headerpreamble{\xtp\crcr}%
}
\dvr{4}}
\def\getHDdimen#1{%
\hdsize=0pt%
\getsize#1\cr\end\cr%
}
\def\getsize#1\cr{%
\endsizefalse\savetks={#1}%
\expandafter\lookend\the\savetks\cr%
\relax \ifendsize \let\next\relax \else%
\setbox\hdbox=\hbox{#1}\newhdsize=1.0\wd\hdbox%
\ifdim\newhdsize>\hdsize \hdsize=\newhdsize \fi%
\let\next\getsize \fi%
\next%
}%
\def\lookend{\afterassignment\sublookend\let\looknext= }%
\def\sublookend{\relax%
\ifx\looknext\cr %
\let\looknext\relax \else %
   \relax
   \ifx\looknext\end \global\endsizetrue \fi%
   \let\looknext=\lookend%
    \fi \looknext%
}%
%
%
\def\tablelet#1{%
   \tableLETtokens=\expandafter{\the\tableLETtokens #1}%
}%
\catcode`\@=12

\newcount\sectionct
\newcount\eqct
\sectionct=0 
\def\Equ#1{\global\advance\eqct by1\relax{\number\sectionct.\number\eqct}%
{\expandafter\xdef\csname x1#1\endcsname%
{\number\sectionct.\number\eqct}}%
\saywhatitis{#1}}

\def\resetcts{\global\eqct=0}

\def\Secno{\global\advance\sectionct by1\relax%
\S\number\sectionct.%
\resetcts}

\def\newsection#1#2{\heading{\bf\Secno \saywhatitis{#1} {#2}}\endheading%
{\expandafter\xdef\csname x1#1\endcsname%
{\noexpand\S\number\sectionct}}
}

\def\refer#1{\csname x1#1\endcsname\saywhatitis{#1}}

\def\saywhatitis#1{}

\newcount\biblioct
\global\biblioct=0
\def\Biblio#1{\global\advance\biblioct by1\relax%
{\expandafter\xdef\csname x1#1\endcsname{\number\biblioct}}}


\newcount\hour
\hour=\time
\define\today{ \divide\hour by 60\relax\number\hour:%
\multiply \hour by -60\relax\advance\time by\the\hour\relax\number\time,
\number\month--\number\day--\number\year }

\documentstyle{amsppt}
\magnification = \magstep1
\hsize=6.5truein
\define\Z{{\Bbb Z}}
\define\Q{{\Bbb Q}}
\define\R{{\Bbb R}}
\define\C{{\Bbb C}}
\define\F{{\Bbb F}}
\define\N{{\Bbb N}}
\define\Pj{{\Bbb P}}
\define\9{\,{}^{t}\,}
\define\8#1{{\Cal#1}}
\define\s{\sigma}
\define\Half{{\Cal H}}
\define\inv{^{-1}}

\define\Jac{\operatorname{Jac}}
\define\SH{\operatorname{Semihull}}
\define\tr{\operatorname{tr}}
\define\IM{\operatorname{Im}}
\define\RE{\operatorname{Re}}
\define\PR{{\Cal P}_n(\R)}
\define\Pn{{\Cal P}_n}
\define\PQ{{\Cal P}_n(\Q)}
\define\PZ{{\Cal P}_n(\Z)}
\define\PfourZ{{\Cal P}_4(\Z)}

\define\PsR{{\Cal P}^{\text{semi}}_n(\R)}
\define\PtwosR{{\Cal P}^{\text{semi}}_2(\R)}
\define\PsQ{{\Cal P}^{\text{semi}}_n(\Q)}

\define\classone{\text{type one}}
\define\Classone{\text{Type one}}
\define\<{\langle}
\define\>{\rangle}

\define\SpnZ{\operatorname{Sp}_n(\Z)}

\define\SptwoQ{\operatorname{Sp}_2(\Q)}
\define\SptwoZ{\operatorname{Sp}_2(\Z)}

\define\GLnZ{\operatorname{GL}_n(\Z)}
\define\GLn{\operatorname{GL}_n}
\define\SLn{\operatorname{SL}_n}
\define\SLtwon{\operatorname{SL}_{2n}}
\define\GLtwon{\operatorname{GL}_{2n}}
\define\Spn{\operatorname{Sp}_n}
\define\GSpnQ{\operatorname{GSp}^+_n(\Q)}
\define\GSpn{\operatorname{GSp}^+_n}
%
%
\define\Xn{{\Cal X}_n}
\define\FS{\operatorname{FS}}

\define\twomat#1#2#3#4{\pmatrix#1&#2\\#3&#4\endpmatrix}
\define\smtwomat#1#2#3#4{{\tsize{\bigl({#1\atop#3}{#2\atop #4}\bigr)}}}
\define\smtwobmat#1#2#3#4{{\tsize{\bigl[{#1\atop#3}{#2\atop #4}\bigr]}}}

\define\mm#1{{B_#1}}
\define\rep#1#2#3{{\Cal V}(#1,#2,#3)}
\define\boxop#1#2{#1\boxtimes#2}
\define\mmm#1{{$B_#1$}}

%
\define\Psn{{\Cal P}^{\text{semi}}_n}
\define\supp{\operatorname{supp}}

\define\({\left(}
\define\){\right)}
\define\redtr{\tilde{tr}}

\define\Tr{\operatorname{Tr}}
\define\X2{{\Cal X}_2}
\define\Xsn{{\Cal X}_n^{\text{semi}}}
\define\Xtwos{{\Cal X}_2^{\text{semi}}}
\define\GLtwoZ{\operatorname{GL}_2(\Z)}
\define\GL{\operatorname{GL}}
\define\SLtwoZ{\operatorname{SL}_2(\Z)}

\define\PtwoQ{{\Cal P}_2(\Q)}
\define\PtwoR{{\Cal P}_2(\R)}

\define\mut{\tilde{\mu}}
\define\kap{\kappa}
\define\Sym{\operatorname{Sym}}
\define\mc{\varpi}

\define\TB{\operatorname{THBK}}
\define\Grit{\operatorname{Grit}}
\define\eps{\operatorname{u}}
\define\Diag{\operatorname{Diag}}
\define\Twin{\operatorname{Twin}}
\define\fav{19}
\define\JJ{J}
\define\Span{\operatorname{Span}}
\define\Aut{\operatorname{Aut}}
\define\spin{\operatorname{spin}}
\define\smfourmatNUMBERTWO#1#2#3#4#5#6#7#8#9{{\tsize{\left(%
{{2\atop#4}\atop{#5\atop #7}}%
{{#4\atop#1}\atop{#6\atop #8}}%
{{#5\atop#6}\atop{#2\atop #9}}%
{{#7\atop#8}\atop{#9\atop #3}}%
\right)}}}
\define\smfourmatNUMBERFOUR#1#2#3#4#5#6#7#8#9{{\tsize{\left(%
{{4\atop#4}\atop{#5\atop #7}}%
{{#4\atop#1}\atop{#6\atop #8}}%
{{#5\atop#6}\atop{#2\atop #9}}%
{{#7\atop#8}\atop{#9\atop #3}}%
\right)}}}
\define\smfourmatNUMBERTEN#1#2#3#4#5#6#7#8#9{{\tsize{\left(%
{{10\atop#4}\atop{#5\atop #7}}%
{{#4\atop#1}\atop{#6\atop #8}}%
{{#5\atop#6}\atop{#2\atop #9}}%
{{#7\atop#8}\atop{#9\atop #3}}%
\right)}}}
\define\smfourmatNUMBERTWELVE#1#2#3#4#5#6#7#8#9{{\tsize{\left(%
{{12\atop#4}\atop{#5\atop #7}}%
{{#4\atop#1}\atop{#6\atop #8}}%
{{#5\atop#6}\atop{#2\atop #9}}%
{{#7\atop#8}\atop{#9\atop #3}}%
\right)}}}
\define\smfourmatNUMBERFOURTEEN#1#2#3#4#5#6#7#8#9{{\tsize{\left(%
{{14\atop#4}\atop{#5\atop #7}}%
{{#4\atop#1}\atop{#6\atop #8}}%
{{#5\atop#6}\atop{#2\atop #9}}%
{{#7\atop#8}\atop{#9\atop #3}}%
\right)}}}
\define\smfourmatNUMBERSIXTEEN#1#2#3#4#5#6#7#8#9{{\tsize{\left(%
{{16\atop#4}\atop{#5\atop #7}}%
{{#4\atop#1}\atop{#6\atop #8}}%
{{#5\atop#6}\atop{#2\atop #9}}%
{{#7\atop#8}\atop{#9\atop #3}}%
\right)}}}
\define\smfourmatNUMBEREIGHTEEN#1#2#3#4#5#6#7#8#9{{\tsize{\left(%
{{18\atop#4}\atop{#5\atop #7}}%
{{#4\atop#1}\atop{#6\atop #8}}%
{{#5\atop#6}\atop{#2\atop #9}}%
{{#7\atop#8}\atop{#9\atop #3}}%
\right)}}}
\define\smfourmatNUMBERTWENTY#1#2#3#4#5#6#7#8#9{{\tsize{\left(%
{{20\atop#4}\atop{#5\atop #7}}%
{{#4\atop#1}\atop{#6\atop #8}}%
{{#5\atop#6}\atop{#2\atop #9}}%
{{#7\atop#8}\atop{#9\atop #3}}%
\right)}}}
\define\lastprime{600}
\define\ideal{\bold p\/}
\define\Xtwo{{\Cal X}_2}

\Biblio{Andrianov}
\Biblio{AGM}
\Biblio{BKold}
\Biblio{BK}
\Biblio{Dern}
\Biblio{EZ}
\Biblio{Grit1}
\Biblio{GritHulek}
\Biblio{IbukHash}
\Biblio{Ibuk1}
\Biblio{Ibuk2}
\Biblio{Ibuk3}
\Biblio{Ibuk4}
\Biblio{Ibuk5}
\Biblio{IgusaMFZ}
\Biblio{Iwaniec}
\Biblio{Marschner}
\Biblio{Nipp}
\Biblio{OMeara}
\Biblio{PoorYuenRel}
\Biblio{PoorYuenExt}
\Biblio{PoorYuenComp}
\Biblio{PoorYuenS22}
\Biblio{URL}
\Biblio{RS}
\Biblio{Shimura}
\Biblio{Skoruppa}
\Biblio{Sturm}
\Biblio{Tate}
\Biblio{Witt}
\Biblio{Yoshida}
\Biblio{Yoshida2}


\topmatter
\title  Paramodular cusp forms  \endtitle
\author  Cris Poor,  poor\@fordham.edu \\
David S\. Yuen, yuen\@lakeforest.edu \\
\endauthor
\address{
Department of Mathematics,
Fordham University, Bronx, NY 10458
\vskip0pt {\it Email:}\/ {\rm poor\@fordham.edu}}\endaddress
\address{
Math/CS Department,
Lake Forest College,
555 N. Sheridan Rd.,
Lake Forest, IL 60045
\vskip0pt{\it Email:}\/  {\rm yuen\@lakeforest.edu}}\endaddress
\date\today\enddate 
\abstract
We classify Siegel modular cusp forms of weight two for the paramodular group $K(p)$ 
for primes $p< 600$.  We find that  weight two Hecke 
eigenforms beyond the Gritsenko lifts correspond to certain 
abelian varieties defined over $\Q$ of conductor $p$.    
The arithmetic classification is in the companion article by 
A\. Brumer and K\. Kramer \cite{\refer{BK}}. 
The Paramodular Conjecture, supported by these computations and  consistent with the Langlands
philosophy and the work of H\. Yoshida \cite{\refer{Yoshida}}, is a partial extension to degree 2 of the
Shimura-Taniyama Conjecture. 
These nonlift Hecke eigenforms share 
Euler factors with the corresponding abelian variety $A$ 
and satisfy congruences modulo $\ell$ with
Gritsenko lifts, whenever $A$ has rational $\ell$-torsion.
\endabstract
\endtopmatter

\document
\nopagenumbers
\pageheight{7.5in}  
\rightheadtext{Paramodular cusp forms}
\leftheadtext{Cris Poor,  David S\.  Yuen}

{\bf{\newsection{sec0}{   Introduction.}} }   

In $1980$, H\. Yoshida conjectured that for every abelian surface defined 
over $\Q$, there exists a discrete group $\Gamma \subseteq \SptwoQ$ and 
a degree two Siegel modular form of weight two for $\Gamma$  
with the same $L$-function.  He supported this conjecture by constructing 
lifts and giving specific examples.  A broader context for this conjecture may be found in the 
recent article \cite{\refer{Yoshida2}} of H. Yoshida.  

Systematic computational evidence for nonlifts required   
specification of the discrete group of the 
putative Siegel modular form. The Paramodular Conjecture posits 
the paramodular group $K(N)$ as the group
corresponding to certain rational abelian surfaces of conductor $N$.
Accordingly, this
article studies spaces of Siegel paramodular cusp forms.   
We believe that the examples given here are the
first nonlifts of weight~{two} found.  
Although we have verified the 
equality of some Euler factors in our examples, we have
not proven the equality of any $L$-functions.
%
For a natural number $N$,  the paramodular group $K(N)$ is defined by:
$$
K(N)=\SptwoQ  \cap    
\pmatrix 
* & * & */N & *  \\
N* & * & * & *  \\
N* & N* & * & N*  \\
N* & * & * & *  
\endpmatrix, 
\quad\text{ for $* \in \Z$.}
$$
Let $S_2^k\left( K(N) \right)$ denote the $\C$-vector space of Siegel modular cusp
forms of weight $k$ and degree two with respect to the group $K(N)$.  
%
A general statement of the Paramodular Conjecture, 
a degree two version of the Shimura-Taniyama
Conjecture,  may be found in the companion article \cite{\refer{BK}} by A. Brumer and K.
Kramer.  
%
%
%
Here, however, we focus on the simplest case:  
\proclaim{ \Equ{A1b}  Paramodular Conjecture for rational abelian surfaces of prime conductor}
Let $p$ be a prime.  
There is a bijection between lines of 
Hecke eigenforms $f \in S_2^2\left( K(p) \right) $ that have rational eigenvalues   
and   
are not Gritsenko lifts 
and  isogeny classes of  rational abelian surfaces ${\Cal A\/}$ 
of conductor $p$.  
In this correspondence,  we have 
$$
L\left({\Cal A\/},s,\text{\rm Hasse-Weil}\right)=L\left(f,s, \spin \right).  
$$
\endproclaim

In \cite{\refer{BK}}, the authors obtain a list of the odd numbers $N<1000$ for which an abelian surface of conductor
$N$ could exist; examples of such surfaces are given for most of the primes on this list. The first is $p = 277$ and the known
surfaces of that conductor are all isogenous to the Jacobian $\Cal{A}_{277}$ of the curve
$y^2+y=x^5+5x^4+8x^3+6x^2+2x$.
Our results cover primes $p < \lastprime$  and are consistent with the Conjecture~\refer{A1b}.  
In particular, there are rational nonlift Hecke eigenforms
where the Paramodular Conjecture indicates there should be and there are none 
where it indicates that there should be none. 
\proclaim{ \Equ{A2} Theorem}
For primes $p < \lastprime $ and not in the set 
$
\{277,349,353, 389, 461, 523, 587\}  
$, 
$ S_2^2\left( K(p) \right)$ is spanned by Gritsenko lifts.  
\endproclaim

\proclaim{ \Equ{A3} Theorem}
The subspace of Gritsenko lifts in  
$ S_2^2\left( K(277) \right)$ has dimension $10$ whereas 
$ S_2^2\left( K(277) \right)$ has dimension $11$.  
There is a rational Hecke eigenform $f$ that is not a Gritsenko lift.  
The Euler factors of $L(f,s,\spin)$ for $q=2$, $3$, $5$ and the 
linear coefficients of the Euler factors for $q=7,11$, $13$ agree with those 
of $L({\Cal A}_{277},s,\text{\rm Hasse-Weil})$.  
\endproclaim

The abelian surface ${\Cal A}_{277}$ has rational $15$-torsion.  
We have the following congruence on the modular side.
\proclaim{ \Equ{A4} Theorem }
Let $f$ be as in Theorem~\refer{A3} and be chosen so that   
$f\in S_2^2\left( K(277) \right)(\Z)$ has 
Fourier coefficients of unit content.    
Let the first Fourier Jacobi coefficient of $f$ be $\phi \in J_{2,277}$
and let $R= \Grit(\phi)\in S_2^2\left( K(277) \right)(\Z)$. 
We have 
$
f \equiv R \mod 15  
$.  
\endproclaim
The above theorems answer, in the case of the paramodular group, 
the challenge posed in \cite{\refer{BKold}, pg.~16} to show that the first nonlift of weight two for prime level  
occurs at $ 277$.  
Further examples, all currently conjectural, of weight two paramodular nonlifts and  their  
Hecke eigenvalues and congruences may be found in section~$7$, see Table~5 in particular.   
Please see our website \cite{\refer{URL}} for further data and details.     

Computations in weight $k=2$ pose a special challenge since no dimension formula is known. 
T\. Ibukiyama used trace formula techniques \cite{\refer{Ibuk2}} to 
give $\dim S_2^k\left( K(p) \right)$  for $k \ge 5$ and any prime $p$.
He has recently proven dimension formulae for $k=3$ and $4$, see \cite{\refer{Ibuk5}}.
The ring structure of $ M_2\left( K(p) \right)$, 
the graded ring of Siegel modular forms for $K(p)$, 
has been studied  
for $p=2$, $3$ and $5$ in \cite{\refer{Ibuk3}}, \cite{\refer{Dern}} and \cite{\refer{Marschner}},
respectively.   T\. Ibukiyama has also proven \cite{\refer{Ibuk4}} that $ S_2^2\left( K(p)
\right)=\{0\}$  for primes $p \le 23$.  
These were hitherto the 
only systematic computations concerning paramodular
cusp forms of weight two.

In weight two, standard constructions of Siegel modular forms 
leave something to be desired.  There are no Eisenstein series 
of weight two.  The useful construction of tracing theta series 
from $M_2^k\left( \Gamma_0(p) \right)$ to $M_2^k\left( K(p) \right)$ 
always vanishes in weight two.  See Theorem~{3.5} for the proof of this fact.  
The Gritsenko lift,   
$\Grit: J_{2,p}^{\text{cusp}} \to S_2^2(K(p))$,  
which constructs paramodular forms for $K(p)$ from 
Jacobi forms of index~$p$, is a nontrivial construction; however, the 
subspace of lifts in weight two is precisely the uninteresting space 
for arithmetic geometry.  In order to construct nonlifts of weight two, 
we used a method of integral closure.  

Given two linearly independent Gritsenko lifts $g_1,g_2 \in S_2^2( K(N))$, 
define the space ${\Cal H}(g_1,g_2)=\{(H_1,H_2)\in S_2^4( K(N)) \times S_2^4( K(N)): H_1 g_2=H_2 g_1\}$.  
The map $\imath_{g_1,g_2}: S_2^2(K(N)) \to {\Cal H}(g_1,g_2)$ given by 
$\imath_{g_1,g_2}(f)=(g_1 f, g_2 f)$ injects.  
If $ \dim {\Cal H}(g_1,g_2) \le \dim J_{2,N}^{\text{cusp}}$ for 
some choice of $g_1,g_2$, then $S_2^2(K(N))$ consists entirely of lifts.  
Suppose, on the other hand, that all choices of $g_1,g_2$ yield 
$ \dim {\Cal H}(g_1,g_2) > \dim J_{2,N}^{\text{cusp}}$.  
For $(H_1,H_2) \in {\Cal H}(g_1,g_2)$ but not in 
$ \imath_{g_1,g_2} \Grit(J_{2,N}^{\text{cusp}})   $, 
we might hope that the meromorphic $f=H_1/g_1=H_2/g_2$ is actually holomorphic.  
In this case we use the initial Fourier expansion of $f$ to 
search for an $F \in  S_2^4(K(N))$ with $f^2=F$, or 
$H_1^2=g_1^2F$.  The validity of a weight~$8$ identity 
$H_1^2=g_1^2F$ is proof that $H_1/g_1$ is in the integral closure 
of $M_2(K(N))$ and is holomorphic.  

Thus, to rule out nonlifts of weight two, 
one spans  $S_2^4(K(p))$, while to construct nonlifts one must span 
$S_2^8(K(p))$ as well.  For primes $p < 600$ and $p \ne 499$, 
$S_2^4(K(p))$ was spanned by tracing theta series, by 
multiplying weight two Gritsenko lifts and by smearing with 
Hecke operators.  The ring structure of $M_2(K(p))$ plays a 
crucial role in spanning $S_2^4(K(p))$ because products of lifts do not in general 
span a Hecke stable subspace.  This same multiplication of Fourier series, 
however, is computationally expensive.  
It was difficult to span $S_2^4(K(479))$ of dimension $440$, for example, 
because the number of Gritsenko lifts in $S_2^2(K(479))$ is relatively small, just $8$.    
Filling $S_2^4(K(479))^{+}$, which turned out to have dimension~$341$, required smearing $S_2^2(K(479)) S_2^2(K(479))$ 
with Hecke operators nine times.  This required computing $1.9$~million 
Fourier coefficients for each of the $8$ weight two Gritsenko lifts.  
The Fourier coefficients of the Gritsenko lifts were computed using the 
method of {\sl theta blocks\/} due to 
V. Gritsenko, N\. Skoruppa and D. Zagier, see section~$4$.  
Finding $99$ linearly independent elements in $S_2^4(K(479))^{-}$ was achieved by theta tracing.  

Finite sets of Fourier coefficients that determine vanishing and congruence 
in $S_2^k(K(p))$ for any $k \in \Z^{+}$ are given in section~$5$.  
In particular, a generalization of Sturm's Theorem \cite{\refer{Sturm}} to $n=2$ may 
be of independent interest.  For $f \in M_n^k$, denote the 
Fourier coefficients by $a(T;f)$.  
The Fourier coefficients of a level one elliptic modular 
form $f\in M_1^k$ are in the $\Z$-module spanned by the first $k/12$:  
$\forall n, a(n;f) \in \Z\< a(j;f): j \le k/12\>$.  
In Theorem~{5.15} we prove that for    
level one Siegel forms $f\in M_2^k$,
all the Fourier coefficients $a(T;f)$ are in the $\Z$-module spanned by those 
whose index $T$ has dyadic trace less than or equal to $k/6$: 
$$
\forall f \in M_2^k,\,\forall T,\, a(T;f) \in \Z\< a(S;f): w(S) \le k/6\>.  
$$ 
Section~$8$ contains examples of nonlifts of higher weight, in particular, 
weight~$3$ cusp forms whose construction was requested by    
A\. Ash, P\. Gunnells and M\. McConnell \cite{\refer{AGM}}.  
We plan a sequel that studies $S_2^k(K(N))$ for composite $N$ and 
makes more use of Fourier-Jacobi expansions.

We thank T\. Ibukiyama and N\. Skoruppa for new
results which improved our exposition and the extent of our
computations.   First, T\. Ibukiyama \cite{\refer{Ibuk5}} 
proved new dimension formulas for prime levels for Siegel forms of weights $3$ and $4$, 
including the paramodular case.  Secondly, N\. Skoruppa explained to us  his
joint work with V. Gritsenko and D. Zagier on theta blocks \cite{\refer{Skoruppa}}.  
We thank Fordham University for the more powerful computer acquired through the Office 
of Internal Grants.  We thank N\. Skoruppa and Siegen University for the use of their 
cluster.  The computations of this paper were mostly performed with Mathematica and $C^{++}$ but 
also with Maple, Perl, Pari and Fermat.

 A\. Brumer began promoting this project in 1997 and we thank
him for finally convincing us to undertake it.  Without his unflagging encouragement and
assistance,   this work would not have been attempted, nor once attempted,  
as satisfactorily concluded.  In particular, 
A\. Brumer suggested studying Gritsenko lifts and 
computed the coefficients 
of thousands of the Jacobi forms used here.

{\bf{\newsection{sec1}{   Notation}} }

For a commutative ring $R$, let $M_{m\times n}(R)$ denote the $R$-module of
$m$-by-$n$ matrices with coefficients in $R$.  For $x\in M_{m\times n}(R)$, 
let $x'\in M_{n\times m}(R)$ denote the transpose.  Let
$V_n(R)=\{x\in M_{n\times n}(R): x'=x\}$ be the symmetric $n$-by-$n$ matrices
over $R$; $V_n(\R)$ is a euclidean vector space under the inner product
$\langle x,y\rangle=\tr(xy)$.
For $R\subseteq \R$, an element $x\in V_n(R)$ is called
positive definite, written $x>0$, when $v'xv >0$ for all $v \in \R^n\setminus\{0\}$;
we denote the set of these by   $\Pn(R)$. 
When $x\in V_n(R)$ and $v'xv \ge 0$ for all $v \in \R^n$, 
we write  $x\in\Psn(R)$. 
The half-integral matrices are
$\Xsn=\{ T\in \PsQ: \forall\, v\in\Z^n,\,v'Tv\in\Z\}$ and 
$\Xn=\Xsn \cap \PQ$.

Let $\GLn(R)=\{x\in M_{n \times n}(R): \det(x) \text{ is a unit in $R$} \}$
be the general linear group
and
$\SLn(R)=\{x\in \GLn(R): \det(x) =1 \}$
the special linear group.
For $x \in \GLn(R)$ let $x^*$ denote the inverse transpose.
Let $I_n\in \GLn(R)$ be the identity
matrix and set $J_n=  \smtwomat{0}{I_n}{-I_n}{0}\in  \SLtwon(R)$.                         
The symplectic group is defined by
$\Spn(R)=\{x\in\GLtwon(R): x'J_nx=J_n\}$. 
For $T,u\in\text{GL}_n(\R)$, we define $T[u]=u'Tu$ and  
denote the $\GLnZ$ equivalence class of $T$ by 
$[T]= \cup_{u} T[u]$ for $ u \in \GLnZ$.  
We write $\Gamma_n=\SpnZ$ and for $R \subseteq \R$
define the group of positive $R$-similitudes  by $\GSpn(R)=\{
x\in M_{2n \times 2n}(R): \exists \mu\in {\R}^{+}: g'J_ng=\mu
J_n \}$. Each $\gamma \in \GSpn(R)$ has a unique
$\mu=\mu(\gamma)=\det(\gamma)^{1/n}$.  
For $S\in V_n(R)$, let $t(S) = \smtwomat IS0I  $ 
define a homomorphism $t: V_n(R) \to \Spn(R)$.  
For $U\in \GLn(R)$, let $u(U) = \smtwomat U00{U^*} $ 
define a homomorphism $u: \GLn(R) \to \Spn(R)$.  
We let $\Gamma_0(N) = \{\smtwomat ABCD\in\SpnZ : C\equiv 0 \mod N\}$, 
$G\Delta_n^{+}(R) = \{\smtwomat AB0D\in \GSpn(R)\}$
and $\Delta_n(R) = \{\smtwomat AB0D\in\Spn(R)\}$.  
The group $\Gamma_0(N)$ is normalized by the Fricke involution 
$F_N = \frac{1}{\sqrt{N}} \smtwomat{0}{I_n}{-NI_n}{0}$.

Define the Siegel upper half space
$\Half_n=\{ \Omega\in V_n(\C): \IM \Omega \in \PR\}$.
The group $\GSpn(\R)$ acts on $\Half_n$ as $\gamma\langle
\Omega\rangle=(A\Omega+B)(C\Omega+D)\inv$ for $\gamma=
\smtwomat{A}{B}{C}{D}$. 
For any function $f:\Half_n\to \C$ and any $k \in \Z$,  
we follow Andrianov and, 
letting $\langle n \rangle=n(n+1)/2$,  
define the group action for $\gamma \in \GSpn(R)$ by
$$
(f|_k \gamma)(\Omega)=\mu(\gamma)^{kn-\langle n\rangle }\det(C\Omega+D)^{-k}f(\gamma\langle \Omega\rangle ).
$$
Let $\Gamma$ be a group commensurable with $\Gamma_n$.  
The complex vector space of Siegel modular forms of degree $n$ and weight $k$ 
automorphic with respect to $\Gamma$ 
is denoted by $M_n^k(\Gamma)$ and is
defined as the set of holomorphic $f:\Half_n\to \C$ such that
$f|_k \gamma=f$ for all $\gamma\in \Gamma$ and such that for all $Y_0\in \PR$ 
and for all $\gamma\in \Gamma_n$, 
$f| \gamma$ is bounded on $\{\Omega\in \Half_n: \IM \Omega > Y_0 \}$.
For $f \in M_n^k(\Gamma)$ the Siegel $\Phi$-map
is defined by
$(\Phi f)(\Omega)=\lim_{\lambda\to+\infty} 
f(\smtwomat{i\lambda}{0}{0}{\Omega})$     
and the space of cusp forms is defined by
$S_n^k(\Gamma)=\{ f\in M_n^k(\Gamma): \forall \gamma\in\Gamma_n, \Phi(f| \gamma)=0\}$.  
The graded ring of Siegel modular forms is 
$ M_n(\Gamma) = \oplus_k M_n^k(\Gamma)$ and the graded ideal of cusp forms is 
$ S_n(\Gamma) = \oplus_k S_n^k(\Gamma)$.

Let $e(z)=e^{2\pi i z}$.
By the Koecher principle, an
$f \in S_n^k(\Gamma)$
has a Fourier expansion
$$
f(\Omega)    =\sum 
 a(T;f) e\left( \langle T,\Omega \rangle \right),   
\text{   also written }  
\FS_n(f)    =\sum 
 a(T;f) q^T,   
$$
where the summation is over $T\in \Pn(\Q)$.  
The  $a(T;f)$ satisfy
$a(T[U];f)=\det(U)^k a(T;f)$ for all
$U\in \GLn(\Q)$ such that $u(U)\in\Gamma$.  
We let $\supp(f)=\{T\in \Pn(\Q): a(T;f) \ne 0\}$.  
The set $\supp(f)$ is contained in the lattice dual to 
$\{S\in V_n(\Q): t(S) \in \Gamma\}$.  The integrality 
properties of  $\supp(f)$ are sometimes important.   
For $f \in S_n^k\left(\Gamma_0(N) \right)$, 
we have $\supp(f) \subseteq \Xn$ and 
$a(T[U];f)=\det(U)^k a(T;f)$ for all
$U\in \GLn(\Z)$.  
Let $ {}^N\X2=\{\smtwomat{a}bbc\in \X2: N|a \}$ and 
$ {}^N\X2^{\text{semi}}=\{\smtwomat{a}bbc\in \X2^{\text{semi}}: N|a \}$.  
For $f \in S_2^k\left(K(N) \right)$, 
we have $\supp(f) \subseteq {}^N\X2$ and 
$a(T[U];f)=\det(U)^k a(T;f)$ for all
$U\in \hat{\Gamma}_0(N)$, where 
$ \hat{\Gamma}_0(N)=\< {\Gamma}_0(N), \smtwomat{1}00{-1}\>$.  
The paramodular groups satisfy $\SptwoQ=K(p)\Delta_2(\Q)$ 
so that there is essentially only one Fourier expansion to consider.  
For primes $p$, we have 
$\SptwoQ=\Gamma_0(p)\Delta_2(\Q) \cup \Gamma_0(p)E_1\Delta_2(\Q) \cup \Gamma_0(p)J_2\Delta_2(\Q)$ 
with $E_1 = I_2 \oplus J_1$, so that there are $3$ basic 
Fourier expansions to consider for $\Gamma_0(p)$.

{\bf{\newsection{sec2}{   Paramodular Forms.}} }

For weights $k \ge 5$, the dimensions $\dim S_2^k\left(K(p)\right)$ 
have been given in \cite{\refer{Ibuk2}}; 
T\. Ibukiyama \cite{\refer{Ibuk5}} has recently proven
the dimensions for weights $k=3$ and $4$.  
There are many ways to construct modular forms.  
The difficulty is knowing when to stop and thus these 
dimension formulae are powerful.

\proclaim{ \Equ{C1} Theorem}{\rm (T\. Ibukiyama)}
Let $p\ge 5$ be a prime number.  
$$
\dim S_2^4\left( K(p) \right) =
\frac{p^2}{576}+\frac{p}{8}-\frac{143}{576} +
\left(\frac{p}{96}-\frac{1}{8}\right)\left(\dfrac{-1}{p}\right)
+\frac{1}{8}\left(\dfrac{2}{p}\right)
+\frac{1}{12}\left(\dfrac{3}{p}\right)
+\frac{p}{36}\left(\dfrac{-3}{p}\right)
$$
\endproclaim
\smallskip
\centerline{Table 1. \ Dimensions for weight $4$ paramodular cusp forms.  }
\smallskip
\begintable
 p\| 2 | 3 |5 |7 |11 |13 |17 |19 |23 |29 
|31 |37 |41 |43 |47 |53 | 59 | 61 | 67 \cr
$\dim S_2^4\left(K(p)\right)$
 \| 0 | 0 | 0 | 1 | 1 | 2 | 2 | 3 | 3 | 4 
| 6 | 8 | 7 | 9 | 8 | 10 | 11 | 16 | 17
\endtable\smallskip

Since weight one paramodular forms are trivial \cite{\refer{RS}}, 
only the dimensions for weight two remain unknown for prime level.  
The Paramodular Conjecture~\refer{A1} is not a dimension formula 
but it does provide an accounting for 
$S_2^2(K(p))$ in terms of isogeny classes of rational abelian 
varieties.  

For a group $G$, let $G_1$ and $G_2$ be subgroups that satisfy the 
following finiteness condition:  
$\forall g \in G, \vert G_1\backslash G_1 g G_2 \vert < +\infty$.  
Let $L(G_1, G)$ be the $\C$-vector space with a basis given by the left cosets 
$\cup_{g\in G} \{G_1 g\}$.  The subgroup $G_2$ has a right action 
$L(G_1, G) \times G_2 \to L(G_1, G)$ given on basis elements by 
$(G_1 g, g_2) \mapsto G_1 gg_2$ and extended linearly.  Denote  the 
fixed subspace of $G_2$ by 
$H(G_1,G,G_2)=\{x\in L(G_1, G):\forall g_2\in G_2, xg_2=x\}$,  
this is the space of Hecke operators for the triple $  (G_1,G,G_2)$.  
For a disjoint union $G_1gG_2=\cup_i G_1g_i$, 
set $[G_1gG_2]= \sum_i G_1g_i \in H(G_1,G,G_2)$;  
then $H(G_1,G,G_2)$ is generated by these double cosets.  
We may check that the multiplication of Hecke operators 
$ H(G_1,G,G_2)\times H(G_2,G,G_3) \to H(G_1,G,G_3)$
given by
$(\sum_i G_1 g_i, \sum_j G_2 h_j ) \mapsto \sum_{i,j} G_1 g_i h_j$ 
is well-defined.  
For $G=\GSpnQ$ and subgroups $G_i$ commensurable with $\SpnZ$, 
the necessary finiteness condition is satisfied.  We have an action of the 
Hecke operators on Siegel modular forms 
$M_n^k(G_1) \times H(G_1,G,G_2) \to M_n^k(G_2)$ 
given by 
$(f, \sum_i G_1 g_i) \mapsto \sum_i f|_k g_i$.  
Hecke operators send cusp forms to cusp forms because of the 
factorization $\GSpnQ=\SpnZ \,G\Delta_n^+(\Q)$.  
For $G_1=G_2$, we have the Hecke algebra $H(G_1,G)=H(G_1,G,G_1)$ 
acting on $M_n^k(G_1)$.  If $w \in G$ normalizes $G_1$, then the single 
coset $G_1 w=[G_1 w G_1]$ is a useful Hecke operator 
that is often just abbreviated by $w$.
We use $B(N)$ to denote the Iwahori subgroup of $K(N)$, 
$$
B(N)=\SptwoZ  \cap    
\pmatrix 
* & * & * & *  \\
N* & * & * & *  \\
N* & N* & * & N*  \\
N* & N* & * & *  
\endpmatrix, 
\quad\text{ for $* \in \Z$.}
$$
For a description of the Hecke operators for the group $B(p)$, 
we refer to \cite{\refer{Ibuk1}} and list the results here for the reader's convenience.   
The Hecke operator $T_{m}$ is defined by the double coset 
$\{\gamma \in \GSpn(\Z): \mu(\gamma)=m\}$.  For $f \in M_2^k( B(p) )$ and 
$\smtwomat{a}{b/2}{b/2}{c} \in \Xtwo$ we have: 
$$
\align
&a\left( \smtwomat{a}{b/2}{b/2}{c}; f|T_{q^{\delta}} \right)=
\sum_{\alpha, \beta, \gamma \in \Z: \alpha+\beta+\gamma=\delta;\, \alpha, \beta, \gamma \ge 0}     \\
&\sum_{u \in R(q^{\beta}): a_u \equiv 0 \mod q^{\beta+\gamma};\, b_u \equiv c_u \equiv 0 \mod q^{\gamma} }
a\left(q^{\alpha} \smtwomat{a_u q^{-\beta-\gamma}}{b_u q^{-\gamma}/2}{b_u q^{-\gamma}/2}{c_u q^{ \beta-\gamma}}; f  \right),   
\endalign
$$ 
where $\smtwomat{a_u  }{b_u/2  }{b_u/2  }{c_u  }=u' \smtwomat{a}{b/2}{b/2}{c} u$ 
and $R(q^{\beta}) \subseteq \Gamma_0(p)$ 
is any lift of $\Pj( \Z / q^{\beta}\Z)$ 
under the map $\smtwomat{u_1}{v_1}{u_2}{v_2} \mapsto (u_1,u_2)$.  
The operator $T_q$, for $(q,p)=1$, is then given by $a\left( T; f | T_q \right)=$ 
$$
a\left( qT;f \right) + 
q^{2k-3}\, a\left( \tfrac1{q}T;f \right) + 
q^{k-2}\sum_{j \mod q} a\left( \tfrac1{q}T \smtwobmat{1}{0}{jp}{q};f \right) + 
q^{k-2}\, a\left( \tfrac1{q}T \smtwobmat{q}{0}{0}{1};f \right).   
$$
The same formulas apply to $f \in M_2^k( K(p) )$ because the number of single cosets 
in the double coset is the same.  
For an eigenform $f$, with eigenvalues 
$f|T(q^{\delta})=\lambda_{q^{\delta}} f$, we use the spinor Euler factor 
$Q_q(f;x)$ given for a $q$ prime to the level 
by (this is the palandrome of the factor in \cite{\refer{Ibuk1}}):
$$
Q_q(f,x)=
1-\lambda_{q}x+(\lambda_{q}^2-\lambda_{q^2}-q^{2k-4})x^2
-\lambda_{q}\,q^{2k-3}x^3+q^{4k-6}x^4.  
$$
Following the work of Andrianov \cite{\refer{Andrianov}}, the spinor $L$-function is 
given, for $\RE(s)>>0$, by 
$$
L(f,s,\text{spin})= \prod_{\text{primes $q$}}
Q_q(f,q^{-s})\inv.  
$$

Theta series give us modular forms on $\Gamma_0(p) $ 
and we can use the Hecke operator $\Tr$, given below, to obtain modular forms on $K(p)$.  
In weight two these theta series trace to zero but for weights $k>2$ we can use this
method to construct paramodular forms.  For even weights, this is the only method we have that 
constructs paramodular forms in the minus space.  
The plus and minus spaces of $S_2^k(K(p))$ are defined as follows.  
Define elements $\mu$ and $\mut$ as below.    
We note that $\mu$ normalizes $K(p)$; 
since $\mu^2=-I_4$  the space $S_2^k(K(p))$ decomposes into $\mu$-eigenspaces 
with eigenvalues $\pm 1$ and we set 
$ S_2^k(K(p))^{\pm}= \{ f \in S_2^k(K(p)): f| \mu = \pm f \}$.  
$$
\mu=\frac{1}{\sqrt{p}}
\pmatrix 
0 & 1 & 0 & 0  \\
-p & 0 & 0 & 0  \\
0 & 0 & 0 & p  \\
0 & 0 & -1 & 0  
\endpmatrix;  \quad
\mut=
\pmatrix 
0 & 1 & 0 & 0  \\
-p & 0 & 0 & 0  \\
0 & 0 & 0 & 1  \\
0 & 0 & -\frac1p & 0  
\endpmatrix.  
$$  
Note that $\Gamma_0(N)$ is not a subgroup of $K(N)$; 
we let $\Gamma_0'(N)=K(N) \cap \SptwoZ$.  

\proclaim{ \Equ{B1} Theorem}
Let the following double cosets define Hecke operators:   
$$
\align
\left[B(p)\Gamma_0'(p)\right]  &: M_2^k\left( B(p)\right)  \to M_2^k\left(\Gamma_0'(p)\right) \text{ and }\\
\left[\Gamma_0'(p)K(p)\right]  &: M_2^k\left( \Gamma_0'(p)\right)  \to M_2^k\left(K(p)\right) .  
\endalign
$$
We define  $\Tr: M_2^k\left(\Gamma_0(p)\right)  \to M_2^k\left(K(p)\right) $ by 
$\Tr= \vert_{M_2^k\left(\Gamma_0(p)\right)}  \left[B(p)\Gamma_0'(p)\right] \left[\Gamma_0'(p)K(p)\right]$.  
For all $f\in  M_2^k\left(\Gamma_0(p)\right)$, 
the Hecke operator $\Tr$ 
satisfies $f|\Tr = $  
$$
 \sum_{\beta \mod p} f\vert t \smtwomat{\frac{\beta}{p}}000 
+\left(f\vert E_1\right)\vert \mut 
+ \sum_{\alpha,\beta \mod p} \left(f\vert E_1\right)\vert t \smtwomat{\frac{\beta}{p}}00\alpha
+\sum_{\alpha \mod p}\left(f\vert F_p\right)\vert \mu\, t \smtwomat{0}00\alpha.  
$$ 
Furthermore, we have $ F_p\Tr=\Tr \mu  $ as Hecke operators in 
$H \left( \Gamma_0(p), \GSpnQ, K(p) \right)$.   
\endproclaim
\demo{Proof}
From \cite{\refer{Ibuk1}} we have 
$\Gamma_0'(p)= B(p) \cup B(p)s_2B(p)$   
where 
$$
s_2=\pmatrix 
1 & 0 & 0 & 0 \\
0 & 0 & 0 & -1 \\
0 & 0 & 1 & 0  \\
0 & 1 & 0 & 0 \endpmatrix.   
$$
A calculation shows that 
$$
s_2\inv B(p) s_2 \cap B(p)=\SptwoZ  \cap    
\pmatrix 
* & * & * & *  \\
p* & * & * & p*  \\
p* & p* & * & p*  \\
p* & p* & * & *  
\endpmatrix, 
\quad\text{ for $* \in \Z$.}
$$
Thus we have 
$$
\align
B(p)  &=\bigcup_{\alpha \mod p}
\left( s_2\inv B(p) s_2 \cap B(p) \right) t \smtwomat{0}00\alpha 
\text{ and hence } \\
B(p)s_2B(p) &=\bigcup_{\alpha \mod p}
B(p) s_2\, t \smtwomat{0}00\alpha .   
\endalign
$$
Thus $\left[B(p)\Gamma_0'(p)\right]=$ 
$B(p)+ \sum_{\alpha \mod p}
B(p) s_2\, t \smtwomat{0}00\alpha $ as a Hecke operator.  
For the coset $\Gamma_0'(p) \backslash K(p)$, 
it follows close upon the definitions that 
$$
K(p)  = \Gamma_0'(p)\JJ(p) \cup 
\bigcup_{\beta \mod p}
\Gamma_0'(p) t \smtwomat{\frac{\beta}{p}}000 
\text{ with the notation } 
\JJ(p)=\pmatrix 
0 & 0 & \frac1p & 0 \\
0 & 0 & 0 & 1 \\
-p & 0 & 0 & 0  \\
0 & -1 & 0 & 0 \endpmatrix,    
$$
and hence that 
$\left[\Gamma_0'(p)K(p)\right]=$ 
$\Gamma_0'(p)\JJ(p)+ \sum_{\beta \mod p}
\Gamma_0'(p) t \smtwomat{\frac{\beta}{p}}000 $ as a Hecke operator.    
By the definition of multiplication of Hecke operators we have 
$$
\align
\Tr  &=\Gamma_0(p)
\left( I+ \sum_{\alpha \mod p} s_2\, t \smtwomat{0}00\alpha \right)
\left( \JJ(p)+ \sum_{\beta \mod p} t \smtwomat{\frac{\beta}{p}}000 \right)   \\
 &=\Gamma_0(p)
\left( \JJ(p)+ \sum_{\beta \mod p} t \smtwomat{\frac{\beta}{p}}000  
+\sum_{\alpha \mod p} s_2\, t \smtwomat{0}00\alpha \JJ(p)
+\sum_{\alpha,\beta \mod p} s_2\, t \smtwomat{\frac{\beta}{p}}00\alpha 
\right) .  
\endalign
$$
$\Gamma_0(p)$ has three cusps, which we represent by 
$I$, $E_1=s_2\inv$ and $F_p$.  
After determining the cusp, we select a simple upper triangular element for the 
$\Gamma_0(p)$-coset.   
We have 
$$
\JJ(p)=u\smtwomat{0}{-1}10 F_p\mu \in \Gamma_0(p)\, F_p  \mu.  
$$
For $\alpha \not\equiv 0 \mod p$, there is a $\tau\in\Z$ such that $\alpha\tau+1\equiv 0 \mod p$ 
and we have
$$
s_2\, t \smtwomat{0}00\alpha \JJ(p) \in  
\Gamma_0(p)\, F_p  \mu\, t \smtwomat{0}00\tau  
\text{  for $\alpha\tau+1\equiv 0 \mod p $.} 
$$
For $\alpha \equiv 0 \mod p$, we have 
$
s_2 \JJ(p) \in  
\Gamma_0(p)\, E_1 \mut\,   
$.  
If we note
$s_2 \in \Gamma_0(p) E_1$, then 
the formula given for $\Tr$ follows.  

Now we show $ F_p\Tr=\Tr \mu  $.   
Using the identity 
$
\mu \, t \smtwomat{\frac{\beta}{p}}00\alpha = t \smtwomat{\frac{\alpha}{p}}00\beta  \mu
$
we can see that the $I_4$ and $F_p$ cusps swap:  
$$
F_p \left( \Gamma_0(p) \sum_{\beta \mod p} t \smtwomat{\frac{\beta}{p}}000  \right) =
\Gamma_0(p) F_p \mu  \sum_{\beta \mod p}\mu t \smtwomat{\frac{\beta}{p}}000 =
 \left( \Gamma_0(p)  F_p \mu \sum_{\beta \mod p} t \smtwomat{0}00\beta  \right) \mu.
$$
To see that the $E_1$ cusp is stabilized, note that 
$F_p E_1=-u\smtwomat{0}{1}10 E_1\mu \in \Gamma_0(p)\, E_1\mu$ 
and therefore we have 
$F_p\, \Gamma_0(p) E_1 t\smtwomat{\frac{\beta}{p}}00\alpha
=\Gamma_0(p) E_1\mu t\smtwomat{\frac{\beta}{p}}00\alpha$ 
$=\Gamma_0(p) E_1 t \smtwomat{\frac{\alpha}{p}}00\beta \, \mu$.  
Furthermore, we have 
$F_p\, \Gamma_0(p) E_1 \mut= \Gamma_0(p)\, E_1\mu\mut =\Gamma_0(p)\, E_1\mut\,\mu.  \qed$  
\enddemo

In order to use the formula for $\Tr$ in 
Theorem~\refer{B1}  to provide Fourier expansions of paramodular forms, 
we must be able to expand a theta series at each of the three cusps: 
$I_4$, $E_1$ and $J_2$.   
General formulas from 
Andrianov \cite{\refer{Andrianov}}, Prop~{3.14} may be simplified to give the following results.

\proclaim{ \Equ{B2} Theorem} 
Let $k,q\in\N$.  
Let $Q$ be a $2k$-by-$2k$ even quadratic form with $qQ\inv$ even.  
Let $N(Q)=\{h\in \Z^{2k}/q\Z^{2k}: Qh=0 \mod q \}$ be the $\Z/q\Z$-nullspace of $Q$.  
Define $\vartheta^{Q}[T]:\Half_g \to \C$ 
for $T\in \Z^{2k\times g}/q\Z^{2k\times g}$ 
by 
$$
\vartheta^{Q}[T](\Omega)= \sum_{N\in\Z^{{2k}\times g} }
e\left( \frac12 \< Q[N+\tfrac{1}{q}T],  \Omega \> \right). 
$$
For $g=2$, we have the following expansions at the other two cusps, 
$$
\align
\vartheta^{Q}\vert E_1  &= i^k\, \det(Q)^{-1/2} 
\sum_{h\in N(Q)} \vartheta^{Q}[0,h]  \\
\vartheta^{Q}\vert F_q &= (-1)^{k} \det(Q)^{-1} q^{k}\,
\vartheta^{qQ^{*}}.  
\endalign
$$
Alternatively, we may use
$$
\vartheta^{Q}\vert J_2  = \det(Q)^{-1} 
\sum_{a,b\in N(Q)} \vartheta^{Q}[a,b]  .  
$$
\endproclaim
For small weights, in order to find the linear combinations of the $ \vartheta^{Q}| \Tr$ 
that are cusp forms, it suffices to cancel the constant term.  
To explain this we need the following Lemma.  
\proclaim{ \Equ{B3b}  Lemma} 
The Witt map 
$W: M_2\left( K(N) \right) \to M_1|\smtwomat{N}001 \otimes M_1$, 
given by $(Wf)(\tau_1,\tau_2)= f\smtwomat{\tau_1}00{\tau_2}$, 
is a weight preserving homomorphism of graded rings.  
\endproclaim
\demo{Proof} 
This follows from the fact that 
$ \left( \smtwomat{N}001\inv \SLtwoZ \smtwomat{N}001 \right) 
\oplus \SLtwoZ \subseteq K(N)$.  \qed
\enddemo
Note that for $ f\in M_2^k\left( K(N) \right) $, 
the vanishing of $Wf$ immediately implies the vanishing of $\Phi f$.  
Furthermore, for $k < 12$, the nontrivial $M_1^k$ are spanned by a single 
Eisenstein series and so $Wf$ vanishes precisely when the 
Fourier expansion of $f$ has a constant term of zero.  
We will use this fact to construct paramodular cusp forms of weight four.

The trace $\Tr$ from $M_2\left(\Gamma_0(p)\right)$ to 
$M_2\left(K(p)\right)$ cannot be used on theta series to 
construct weight two forms.  Indeed, $\Tr$  is identically 
zero on theta series in $M_2^2\left(\Gamma_0(p)\right)$.  

\proclaim{ \Equ{B4}  Theorem} 
Let $p$ be an odd prime.  
Let $Q$ be a $4$-by-$4$ even quadratic form of level $p$ and 
square determinant.  In degree two we have 
$  \vartheta^{Q}| \Tr =0$.  
\endproclaim

Since $M_2^2\left(K(2)\right)=\{0\}$, the above Theorem also holds for
$p=2$, see \cite{\refer{Ibuk3}}. 
The remainder of this section is devoted to proving Theorem~\refer{B4} 
and its consequences.  

\proclaim{ \Equ{B5} Lemma} 
Let $p$ be an odd prime.  
Let $Q$ be a $4$-by-$4$ even quadratic form of level $p$ and 
square determinant. 
Then $\det(Q)=p^2$, the Hasse invariant of $Q$ is $(-1,-1)$ 
and $Q$ is equivalent over $\Q$ to $\Diag(1,\eps,p,\eps  p)$, 
where the resdiue class of $\eps$ in $\F_p$ is a nonsquare unit.  
Furthermore, $Q$ has the property
$$
\forall b \in \Z^4,\quad  b'Qb \equiv 0 \mod p^2 
\iff  b \equiv 0 \mod p .  
\tag \Equ{B5}
$$
\endproclaim
\demo{Proof} 
This is elementary, compare page~{152} of \cite{\refer{OMeara}}. 
\enddemo

The notion of twinning is used in the proof of Theorem~\refer{B4}.  
Recall the Fricke involution in degree one: $F_p=\frac1{\sqrt{p}}\smtwomat{0}{1}{-p}{0}$.  
\proclaim{ \Equ{B6} Definition} 
For $T=\smtwomat{a}bbc \in\PtwoQ$, 
define $\Twin(T)=F_p' T F_p=
\smtwomat{pc}{-b}{-b}{{a\over p}} \in\PtwoQ$.  
\endproclaim

Notice that twinning stabilizes ${}^p\X2$ and respects the 
$\Gamma_0(p)$-equivalence class of $T$.  Twins do have the 
same determinant but the $\GLtwoZ$-equivalence classes may differ.  
For example, 
$T=\smtwomat{10}554 \in{}^5\X2$ has 
$\Twin(T)=\smtwomat{20}552 \in{}^5\X2$ but 
$T=\smtwomat{10}554\in \left[ \smtwomat{4}114 \right]$ 
whereas 
$\Twin(T)= \smtwomat{20}552 \in \left[ \smtwomat{2}118 \right]$.  
The following elements twin and rescale the Fourier coefficients of a 
Siegel modular form:    

$$
\mu=\frac{1}{\sqrt{p}}
\pmatrix 
0 & 1 & 0 & 0  \\
-p & 0 & 0 & 0  \\
0 & 0 & 0 & p  \\
0 & 0 & -1 & 0  
\endpmatrix;  \quad
\mut=
\pmatrix 
0 & 1 & 0 & 0  \\
-p & 0 & 0 & 0  \\
0 & 0 & 0 & 1  \\
0 & 0 & -\frac1p & 0  
\endpmatrix;  \quad
\kap=
\pmatrix 
0 & \frac1p & 0 & 0  \\
-1 & 0 & 0 & 0  \\
0 & 0 & 0 & p  \\
0 & 0 & -1 & 0  
\endpmatrix.  
$$

\proclaim{ \Equ{B8} Lemma} 
Let $f\in M_2^k(\Gamma)$ have Fourier coefficients $a(T;f)$.
$$
\alignat2
a\left( T; f|\mu\right)  &= a\left( \Twin(T); f\right), & 
\supp(f|\mu) &= \Twin\left( \supp(f) \right), \\
a\left( T; f|\mut\right)  &= p^k a\left(\frac1p \Twin(T); f\right), & \qquad
\supp(f|\mut) &= p\Twin\left( \supp(f) \right), \\
a\left( T; f|\kap\right)  &= p^{-k} a\left( p \Twin(T); f\right), &
\supp(f|\kap) &= \frac1p\Twin\left( \supp(f) \right) .  
\endalignat
$$ 
\endproclaim
If $f$ has $\mu$-eigenvalue $\pm 1$, then 
$a\left( \Twin(T); f\right)=\pm a\left( T; f\right)$, so that 
the $\mu$-eigenspace is determined by the twins.  
It is useful to define a type of projection to ${}^p\X2$.  

\proclaim{ \Equ{B9} Definition} 
Define $\pi: M_2(\Gamma(p))\to {\Cal O}(\Half_2)$ by: 
If 
$f(\Omega)= \sum_{T\in \frac1p \Xtwos}
a(T) e\left( \< \Omega, T\> \right)$ 
then 
$(\pi f)(\Omega)= \sum_{T\in {}^p \Xtwos}
a(T) e\left( \< \Omega, T\> \right)$.   
\endproclaim

For $f\in M_2(\Gamma_0(p))$, the four terms in $f | \Tr$,   
$$
 \sum_{\beta \mod p} f\vert t \smtwomat{\frac{\beta}{p}}000 
+\left(f\vert E_1\right)\vert \mut 
+ \sum_{\alpha,\beta \mod p} \left(f\vert E_1\right)\vert t \smtwomat{\frac{\beta}{p}}00\alpha
 +\sum_{\alpha \mod p}\left(f\vert F_p\right)\vert \mu\, t \smtwomat{0}00\alpha,  
$$ 
have the following images under $\pi$:  
$$
\pi\left(\sum_{\beta \mod p} f\vert t\smtwomat{\frac{\beta}{p}}000 \right) 
=p\, \pi(f), \qquad
\pi\left(f\vert E_1\mut \right) 
=
f\vert E_1 \mut, 
$$
$$
\pi\left(
\sum_{\alpha,\beta \mod p} \left(f\vert E_1\right)\vert t \smtwomat{\frac{\beta}{p}}00\alpha
\right)
=
p^2  \pi(f|E_1), \quad
\pi\left(
\sum_{\alpha \mod p}\left(f\vert F_p\right)\vert \mu\, t \smtwomat{0}00\alpha 
\right)
=
p \, \pi(f|F_p \mu).  
$$
We will see that, when applied to weight two theta series, 
the sum of the first and third terms, and the the sum of the second and fourth terms, 
cancel separately.   

\proclaim{ \Equ{B10} Lemma} 
Let $p$ be an odd prime.  
Let $Q$ be a $4$-by-$4$ even quadratic form of level $p$ and 
square determinant. 
For $b \in N(Q)$, 
$\pi \left( \vartheta^{Q}[0,b] \right)=0$ 
unless $b \equiv 0 \mod p$.  
We have 
$$
\pi \left( \vartheta^{Q}| E_1 \right)=
-\frac1p 
\pi \left( \vartheta^{Q}\right) \text{ and } 
\sum_{\beta \mod p} \vartheta^{Q}\vert t \smtwomat{\frac{\beta}{p}}000 +
\sum_{\alpha,\beta \mod p} \left(\vartheta^{Q}\vert E_1\right)\vert t
\smtwomat{\frac{\beta}{p}}00\alpha =0.  
$$
\endproclaim
\demo{Proof}
We have 
$ \vartheta^{Q}[0,b](\Omega)= 
\sum_{c,d\in\Z^2} e \left( \frac12\<\pmatrix c \\ d+b/p \endpmatrix' Q
\pmatrix c \\ d+b/p \endpmatrix, \Omega   \> \right)$.  
We may write this as 
$ \vartheta^{Q}[0,b](\Omega)= 
\sum_{c,d\in\Z^2} e \left( \frac12\<T, \Omega   \> \right)$ for 
$$
T=
\pmatrix
 c'Qc & c'Qd +\frac1p b'Qc \\c'Qd +\frac1p b'Qc & 
d'Qd +\frac2p b'Qd+\frac1{p^2} b'Qb 
\endpmatrix .
$$
If $b \in N(Q)$ and $T\in {}^p\X2$ then $b'Qb\in p^2\Z$ and so 
by property~\refer{B5} we have $b \equiv 0 \mod p$.  

For the second part, notice 
$\vartheta^{Q}\vert E_1  = - \frac1p 
\sum_{b\in N(Q)} \vartheta^{Q}[0,b]$ by Theorem~\refer{B2} so that we have 
$\pi \left( \vartheta^{Q}| E_1 \right)=
-\frac1p 
\pi \left( \vartheta^{Q}\right)$.  
To prove the third part, notice that 
$\sum_{\beta \mod p} \vartheta^{Q}\vert t \smtwomat{\frac{\beta}{p}}000 +
\sum_{\alpha,\beta \mod p} \left(\vartheta^{Q}\vert E_1\right)\vert t
\smtwomat{\frac{\beta}{p}}00\alpha$ is already supported on $ {}^p\X2$. 
Therefore, it it equal to its projection, which is 
$ p\, \pi \left( \vartheta^{Q} \right) + p^2 \pi \left( \vartheta^{Q} |E_1\right)$ 
$=  p\, \pi \left( \vartheta^{Q} \right) + p^2  \left( -\frac1p 
\pi \left( \vartheta^{Q}\right) \right)=0$.  \qed 
\enddemo
 
\proclaim{ \Equ{B11} Lemma} 
Let $p$ be an odd prime.  
Let $Q$ be a $4$-by-$4$ even quadratic form of level $p$ and 
square determinant. 
For $a,b \in N(Q)$, 
$\pi \left( \vartheta^{Q}[a,b]| \mut \right)=0$ 
unless $a \equiv 0 \mod p$.  
We have 
$$
\pi \left( \vartheta^{Q}| F_p\mu \right)=
-\frac1p 
\pi \left( \vartheta^{Q}| E_1 \mut\right) \text{ and } 
\left(\vartheta^{Q}\vert E_1\right)\vert \mut +
\sum_{\alpha \mod p}\left(\vartheta^{Q}\vert F_p\right)\vert \mu\, t \smtwomat{0}00\alpha
  =0.  
$$
\endproclaim
\demo{Proof}
We have 
$ \vartheta^{Q}[a,b](\Omega)= 
\sum_{c,d\in\Z^2} e \left( \frac12\<\pmatrix c+a/p \\ d+b/p \endpmatrix' Q
\pmatrix c+a/p \\ d+b/p \endpmatrix, \Omega   \> \right)$.  
We may write this as 
$ \vartheta^{Q}[a,b](\Omega)= 
\sum_{c,d\in\Z^2} e \left( \frac12\<T, \Omega   \> \right)$ for 
$$
T=
\pmatrix
 c'Qc +\frac2p a'Qc+\frac1{p^2} a'Qa & * \\ * & * 
\endpmatrix .
$$
We have $\pi \left( \vartheta^{Q}[a,b]| \mut \right)=0$ unless 
$ \supp(\vartheta^{Q}[a,b]|\mut) \cap {}^p\X2$ is nonempty.  
However, since $ \supp(\vartheta^{Q}[a,b]|\mut)=p\Twin(\supp(\vartheta^{Q}[a,b]))$, 
this is equivalent to
$T\in \supp(\vartheta^{Q}[a,b])\cap \frac1p ({}^p\X2)$.  
From $a \in N(Q)$ and $T\in  \frac1p ({}^p\X2)$ we derive 
$a'Qa \equiv 0 \mod p^2$.  
By property~\refer{B5} we have $a \equiv 0 \mod p$.  

For the second part, notice 
$\vartheta^{Q}| F_p\mu=
\vartheta^{Q}| J_2\mut=
 \frac{1}{p^2}\sum_{a,b\in N(Q)} \vartheta^{Q}[a,b]| \mut$ 
by  Theorem~\refer{B2}.  
We have 
$\pi \left( \vartheta^{Q}| F_p\mu \right)=
\frac{1}{p^2}\pi \left( \sum_{a,b\in N(Q)} \vartheta^{Q}[a,b]| \mut\right)$ 
$=\frac{1}{p^2}\pi \left( \sum_{b\in N(Q)} \vartheta^{Q}[0,b]| \mut\right)$  
by the first part.  Recognizing the formula for $\vartheta^{Q}\vert E_1$ 
from Theorem~\refer{B2}, we may conclude that  
$\pi \left( \vartheta^{Q}| F_p\mu \right)
=\frac{1}{p^2}\pi \left(-p  \vartheta^{Q}\vert E_1|
\mut\right)$ 
$=-\frac1p 
\pi \left( \vartheta^{Q}| E_1 \mut\right)$.  

To prove the third part, notice that 
$\left(\vartheta^{Q}\vert E_1\right)\vert \mut +
\sum_{\alpha \mod p}\left(\vartheta^{Q}\vert F_p\right)\vert \mu\, t \smtwomat{0}00\alpha$ 
is already supported on $ {}^p\X2$. 
Therefore, it it equal to its projection, which is 
$ \pi \left( \vartheta^{Q}|E_1\mut \right) +
p\, \pi \left( \vartheta^{Q}|F_p\mu \right)$
$=\pi \left( \vartheta^{Q}|E_1\mut \right) +
p  \left( -\frac1p 
\pi \left( \vartheta^{Q}| E_1 \mut\right) \right)=0$.  
 \qed 
\enddemo
 
\demo{Proof of Theorem~\refer{B4}} 
We add the final conclusions from Lemmas~\refer{B10} and  \refer{B11}.  \qed
\enddemo
Here is a consequence of Theorem~\refer{B4} that is useful for proving the integrality 
of modular forms constructed by theta tracing.  

\proclaim{Theorem  \Equ{B12}}
Let $p$ be an odd prime.
Let $A,B\in\PfourZ$ be even with determinant $p^2$ and 
$pA^{-1}$,  $pB^{-1}$ both even.
The Fourier coefficients of $\Tr(\vartheta^Q)$ are  multiples of 4.
\endproclaim
\demo{Proof} 
Let $Q=A\oplus B$, an $8\times8$ matrix.  
First, by the formulas in \refer{B1}, \refer{B2} and \refer{B8}, 
the constant term of $\Tr(\vartheta^Q)$ is
$p+p^4/p^2+p^2/p^2+p=(p+1)^2$ which is a multiple of 4.
We now show that $a(T; \Tr(\vartheta^Q))\equiv0\mod4$ for the nonzero $T\in\,^p\Xtwo$. 
By Theorem~\refer{B4}, we have that $\Tr(\vartheta^A)=0=\Tr(\vartheta^B)$.  
For $f=\vartheta^Q-\vartheta^A-\vartheta^B$ and $g=\Tr(f)$, 
it suffices to prove that $a(T; g)\equiv0\mod4$.  We have 
$$
a(T; g)=a(T; \Tr(f))=
p\,a(T; f) + a(T; f|E_1\tilde\mu) + p^2\,a(T; f|E_1)+p\,a(T; f|F_p\mu).
$$
We will show that the Fourier coefficients
of each of the four terms above, treated individually, are multiples of 4. 
When using the slash operator, it must be remembered
that $f\in M_2^2(K(p)) \oplus M_2^4(K(p))$ is a sum of forms of different weights.  
Considering the first term, we wish to show that
$$
p\, a(T; \vartheta^Q) - p\, a(T; \vartheta^A) - p\, a(T; \vartheta^B) \equiv 0\mod 4.
$$
We have
$$
\align
a(T; \vartheta^Q)&=\#\{w\in\Z^{8\times2}:w'Qw=2T\}=\#\{(u,v):u,v\in\Z^{4\times2},u'Au+v'Bv=2T\},  \\
a(T; \vartheta^A)&=\#\{u\in\Z^{4\times2}:u'Au=2T\}, \\
a(T; \vartheta^B)&=\#\{v\in\Z^{4\times2}:v'Bv=2T\}.
\endalign
$$
Then, using that $T$ is nonzero, 
$$
a(T; \vartheta^Q) {-} a(T; \vartheta^A) {-} a(T; \vartheta^B)=
\#\{(u,v):u,v\in\Z^{4\times2},u\ne0,v\ne0,u'Au+v'Bv=2T\}.  
$$

This set can be partitioned into subsets of 4 elements each
of the form $\{\pm u, \pm v\}$ because the $u,v$ are nonzero.
This proves that the above number is a multiple of 4.
Now consider the second term.
We wish to show that
$$
a(T; \vartheta^Q|E_1\tilde\mu) - a(T; \vartheta^A|E_1\tilde\mu) - a(T; \vartheta^B|E_1\tilde\mu) \equiv 0\mod 4.
$$
Denoting $S={1\over p}\text{Twin}(T)$,
the lefthand side is
$$
p^4\,a(S; \vartheta^Q|E_1) - p^2\,a(S; \vartheta^A|E_1) - p^2\,a(S; \vartheta^B|E_1).
$$
Applying formula~\refer{B2} for slashing with $E_1$, this becomes 
$$
\frac{p^4}{p^2}\,a(S; \sum_{\beta\in N(Q)} \vartheta^Q[0,\beta]) + \frac{p^2}{p}\,a(S; \sum_{h\in N(A)} \vartheta^A[0,h])
+ \frac{p^2}{p}\,a(S; \sum_{\ell\in N(B)} \vartheta^B[0,\ell]).
$$
Consider $\beta=({h\atop\ell})\in N(Q)=N(A\oplus B)$.
We break the above sum into four parts: 
$$
\align
& \sum \{p^2\, a(S;\vartheta^Q[0,(\tsize{h\atop\ell})]):
(\tsize{h\atop\ell})\in N(Q), h\ne 0, \ell \ne 0\}\tag a\\
+&\sum\{p^2\, a(S; \vartheta^Q[0,(\tsize{h\atop0})])+p\,a(S; \vartheta^A[0,h]) :
h\in N(A), h \ne 0\}\tag b\\
+&\sum\{p^2\, a(S; \vartheta^Q[0,(\tsize{0\atop\ell})])+p\,a(S; \vartheta^B[0,\ell]) :
\ell\in N(B), \ell \ne 0\}\tag c\\
+&p^2 \,a(S; \vartheta^Q[0,(\tsize{0\atop 0})])+p\,a(S; \vartheta^A[0,0])
+p\,a(S; \vartheta^B[0,0]). \tag d
\endalign
$$

We show parts (a) to (d) individually are $0$ modulo $4$.
Part (a) is
$$
\align
p^2\, \#\{(u,v,h,\ell):\ &u,v\in \Z^4,(\tsize{h\atop\ell})\in N(Q), h\ne0, \ell\ne 0,\\
&(u+h/p)'A(u+h/p)+(v+\ell/p)'B(v+\ell/p)=2S\}.
\endalign
$$
We can partition these $(u,v,h,\ell)$ into subsets of the form
$$
\{(u,v,h,\ell), (u,-v,h,-\ell), (-u,v,-h,\ell), (-u,-v,-h,-\ell)\}.
$$
Because $p$ is odd, we have $h\not\equiv-h\mod p$ and $\ell\not\equiv-\ell\mod p$
for nonzero $h,\ell$.
Thus these subsets always have 4 distinct elements and this proves that Part (a) is a multiple of 4.
To analyze Part (b), note that 
$(\tsize{h\atop0}) \in N(Q)$ if and only if $h \in N(A)$, so that 
$$
\align
&\sum \{a(S; \vartheta^Q[0,(\tsize{h\atop0})]):h\in N(A), h\ne0\}  
= \\
 \#  &\{(u,v,h): u,v\in\Z^4,h\in N(A), h\ne0,
(u+h/p)'A(u+h/p)+v'Bv=2S\}.
\endalign
$$
We can put these $(u,v,h)$ into equivalence classes of the form
$$
\{(u,v,h),(u,-v,h),(-u,v,-h), (-u,-v,-h)\}.
$$
These classes will have four elements unless $v=0$.
So modulo 4, we can ignore all but the case where $v=0$, and
thus modulo 4, Part (b) is equivalent to
$$
\align
&\sum p^2\, \#\{(u,h): u\in\Z^4,
(u+h/p)'A(u+h/p)=2S\}+ \\
&\sum\{p\,a(S; \vartheta^A[0,h]) :
h\in N(A), h \ne 0\}
\endalign
$$
But this is equal to
$$
\sum\{(p^2+p)\, \#\{(u,h): u\in\Z^4,
(u+h/p)'A(u+h/p)=2S\} :
h\in N(A), h \ne 0\}.  
$$
Because we can pair $\pm(u,h)$,  
$\#\{(u,h): u\in\Z^4,(u+h/p)'A(u+h/p)=2S\}$ 
is an even number;  
since $p^2+p$ is also even, the above is a multiple of 4.
Thus Part (b) is a multiple of 4.
Similarly, Part (c) is a multiple of 4.
Finally, Part (d) can be rewritten as
$$
\align
& p^2 \,\left(a(S; \vartheta^Q[0,(\tsize{0\atop 0})])- a(S; \vartheta^A[0,0])
-a(S; \vartheta^B[0,0]) \right) \\
+ &(p^2+p)\, a(S; \vartheta^A[0,0])
+ (p^2+p)\, a(S; \vartheta^B[0,0]).
\endalign
$$

Here the first line is a multiple of 4 by the exact argument as for the first term.
The second line is a multiple of 4 because $a(S; \vartheta^A[0,0])$ and
$a(S; \vartheta^B[0,0])$ are even because $S$ is
nonzero.
This proves the second term is $0$ modulo $4$.
The third term is
$$
p^2\,a(T; \vartheta^Q|E_1) - p^2\,a(T; \vartheta^A|E_1) - p^2\,a(T; \vartheta^B|E_1),
$$
and the same argument used for the second term will show that this
is also $0$ modulo 4.
The fourth term can be rewritten as
$$
p\,a(T; \vartheta^{pQ^*}|\mu)-p\,a(T; \vartheta^{pA^*}|\mu)-p\,a(T; \vartheta^{pB^*}|\mu).
$$
Noting that $pQ^*=(pA^*)\oplus(pB^*)$,  this is
$$
p\,a(S; \vartheta^{(pA^*)\oplus(pB^*)})
-p\,a(S; \vartheta^{pA^*})-p\,a(S; \vartheta^{pB^*}).
$$
for $S=\text{Twin}(T)$.
The same argument used for the first term 
shows that this is a multiple of 4.
Having shown all four terms are $0$ modulo $4$, we conclude that $a(T; g)\equiv0\mod4$.  \qed
\enddemo

\proclaim{Theorem \Equ{B14}}
Let $g\in S_2^k(K(N))(\Z)$ for $k \ge 2$.  
Let $T_q$ be the Hecke operator for a fixed prime $q$.  
Then the following congruence holds:  
$$
T_q(g^q)\equiv  g \mod q.
$$
Furthermore, if $\phi\in J_{qk,N}^{\text{cusp}}$ is the first Fourier Jacobi coefficient of
$T_q(g^q)$ and $\psi\in J_{k,N}^{\text{cusp}}(\Z)$ that of $g$, then 
$\Grit(\psi)\equiv \Grit(\phi)\mod q$. 
\endproclaim
\demo{Proof}
Take any $T\in\,{}^N\Xtwo$.
Then we have 
$$
a(qT; g^q) = \sum_{s_i \in  {}^N\Xtwo: s_1+\cdots+s_q=qT}a(s_1;g)\cdots a(s_q; g).
$$
Since $q$ is prime,
unless $s_1=\cdots=s_q$, then there are a multiple of $q$ nontrivial ways
to permute the $s_1,\ldots,s_q$.
Thus modulo $q$, we have
$
a(qT; g^q) \equiv a(T; g)^q\mod q 
$.  
From Fermat's congruence $x^q\equiv x\mod q$ for any integer $x$,
we have
$
a(qT; g^q)\equiv a(T; g)\mod q  
$.  
Next, in the formula for any cusp form given in this section 
$$
a(T; T_q(f)) = a(qT; f) + \text{terms with coefficients that are positive powers of }q,
$$
as long as $f$ has weight at least $3$.
Here, the weight of $g^q$ is at least $4$,
so we can conclude
$
a(T; T_q(g^q))\equiv a(qT; g^q)\mod q  
$.  
Thus $a(T; T_q(g^q))\equiv a(T; g)\mod q$ for all $T\in\,{}^N\Xtwo$.
Hence we have the first assertion 
$T_q(g^q)\equiv g\mod q$.  

Letting $\phi\in J_{qk,N}^{\text{cusp}}$ be the first Fourier-Jacobi coefficient of $T_q(g^q)$,
then $\phi$ is congruent modulo $q$ to the first Fourier Jacobi coefficient of $g$;  
that is, 
$\phi\equiv\psi\mod q$.
Now, if $T=\smtwomat{mN}{r/2}{r/2}{n} \in\,{}^N\Xtwo$, then
$$\align
a(T; \Grit(\phi))&=\sum_{\delta|(n,rm)}\delta^{qk-1} c({mn\over\delta^2},{r\over\delta}; \phi),\\
a(T; \Grit(\psi))&=\sum_{\delta|(n,rm)}\delta^{k-1} c({mn\over\delta^2},{r\over\delta}; \psi).  
\endalign
$$
Since $\delta^{qk}\equiv\delta^k\mod q$
and since $c({mn\over\delta^2},{r\over\delta}; \phi)\equiv c({mn\over\delta^2},{r\over\delta}; \psi)\mod q$,
we may conclude 
$a(T; \Grit(\phi))\equiv a(T; \Grit(\psi))\mod q$; 
that is,
$\Grit(\phi)\equiv \Grit(\psi)\mod q$.  \qed
\enddemo

{\bf{\newsection{sec3}{   Jacobi Forms.}} }

The Gritsenko lift  constructs paramodular 
forms from Jacobi forms and so we need to compute spaces of 
Jacobi forms.
The basic reference for Jacobi forms is the book of Eichler and Zagier \cite{\refer{EZ}}.  
The weight~$2$ Jacobi forms in this article were originally computed from 
weight~$3/2$ modular forms on $\Gamma_0(4p)$.  
A more appealling technique, however, is the method
of theta blocks, which seems to work very well for low weight.  
We thank N\. Skoruppa, see \cite{\refer{Skoruppa}}, for explaining 
to us his joint work  with 
V\. Gritsenko and D\. Zagier on theta blocks.  

The following  subgroup $\Gamma_{\infty}(\Z)$ of $\SptwoZ $ stabilizes the 
Fourier-Jacobi expansion of a level one Siegel modular form  term by term 
and this gives some motivation for the definition of Jacobi forms.    
$$
\Gamma_{\infty}(\Z)=\SptwoZ  \cap    
\pmatrix 
* & 0 & * & *  \\
* & * & * & *  \\
* & 0 & * & *  \\
0 & 0 & 0 & *  
\endpmatrix, 
\quad\text{ for $* \in \Z$.}
$$
To admit half-integral weights in the following definition, 
we require that $\mu(\gamma,\Omega)=\chi(\gamma) \det(C\Omega+D)^k$ 
be a factor of automorphy on $\Gamma_{\infty}(\Z) \times \Half_2$, 
compare \cite{\refer{GritHulek}}.  
We write $q=e(\tau)$ and $\zeta=e(z)$. 
%
\proclaim{Definition \Equ{b1}}
A level one Jacobi form of weight $k\in \frac12 \Z$, index $m\in \Q$ and 
multiplier $\chi:\Gamma_{\infty}(\Z)\to \C$, 
denoted $\phi \in J_{k,m}(\chi)$, 
is a holomorphic map  
$\phi:\Half_1 \times \C  \to\C $ given by 
$(\tau,z)   \mapsto \phi(\tau,z)$  
such that, if ${\tilde \phi}:\Half_2\to\C$ is defined by 
${\tilde \phi}(\Omega)= \phi(\tau,z) e( m \omega  )$ for 
$\,\Omega =\pmatrix \tau & z \\ z & \omega \endpmatrix \in \Half_2$,
then we have \newline 
1.)  $\forall \gamma\in \Gamma_{\infty}(\Z),\ {\tilde \phi}\vert_{k} \gamma =\chi(\gamma) {\tilde \phi}$, and
\newline 2.) $\phi(\tau,z)= \sum_{n\ge 0,r\in\Z} c(n,r)q^n\zeta^r$, where $c(n,r)=0$ unless $4mn-r^2
\ge 0$.  
\endproclaim
If $c(n,r)=0$ unless $4mn-r^2 > 0$ then $\phi$ is called a {\it cusp\/} form and we write 
$\phi \in J_{k,m}^{\text{cusp}}(\chi)$.   
In \cite{\refer{EZ}, pp\. 121, 131-132} we can find dimension formulae for Jacobi forms.  
We thank N\. Skoruppa for rewriting these for us.
\proclaim{Theorem \Equ{c0}}  
For $k \in \Z_{ \ge 0}$, let $\{\{k\}\}=\dim S_1^k$.  
For $m \in \N$, let $\s_0(m)$ be the number of positive divisors of $m$.  
Let $\delta(k,m)$ be zero unless $k=2$ and let 
$\delta(2,m)= \frac12\s_0(m)-1 $ for nonsquare $m$ and 
$\delta(2,m)= \frac12\s_0(m)-\frac12 $ for square $m$.  
$$
\align
\text{For even $k \ge 2$,   }\quad
\dim J_{k,m}^{\text{cusp}}  &= \delta(k,m)+ 
\sum_{j=1}^m \left(  \{\{ k+2j\}\}- [[ \frac{j^2}{4m} ]]  \right).  \\
\text{For odd $k \ge 3$,   }\quad
\dim J_{k,m}^{\text{cusp}}  &= 
\sum_{j=1}^{m-1} \left(  \{\{ k+2j-1\}\}- [[ \frac{j^2}{4m} ]]  \right). 
\endalign
$$
\endproclaim

\smallskip
\centerline{Table 2. \ Dimensions for Jacobi cusp forms of weight $2$ and index $p$.  }
\smallskip
\begintable
 p\| 2 | 3 |5 |7 |11 |13 |17 |19 |23 |29 
|31 |37 |41 |43 |47 |53 | 59 | 61 | 67 \cr
$\dim J_{2,p}^{\text{cusp}}$
 \| 0 | 0 | 0 | 0 | 0 | 0 | 0 | 0 | 0 | 0 
| 0 | 1 | 0 | 1 | 0 | 1 | 0 | 1 | 2
\endtable\smallskip

The Fourier expansions of theta blocks can be computed from 
the Dedekind eta function 
$
\eta(\tau)=q^{\frac{1}{24}} 
\prod_{n \in \N} (1-q^n) 
$
and the Jacobi theta function:  
$$
\align
\vartheta(\tau, z) &=-i \theta[{\tsize \matrix 1/2 \\ 1/2\endmatrix}](z,\tau)  = -i 
\sum_{n \in \Z} e\left( \frac12 (n +\frac12)^2\tau + (n +\frac12)(z+\frac12) \right)\\
&=\sum_{n \in \Z} (-1)^n q^{\frac{(2n+1)^2}{8}}\zeta^{\frac{2n+1}{2}}  \\
&=q^{\frac18}\left( \zeta^{\frac12}- \zeta^{-\frac12} \right) 
\sum_{n \in \N}(-1)^{n+1}\,q^{\binom{n}{2}}\,\sum_{j\in\Z: \vert j \vert \le n-1 }
\zeta^j \\
&=q^{\frac18}\left( \zeta^{\frac12}- \zeta^{-\frac12} \right) 
\left(1-q(\zeta+1+\zeta^{-1})+q^{3} (\zeta^{2}+\zeta+1+\zeta^{-1}+\zeta^{-2})+\dots \right)
\endalign
$$

A modular form may be viewed as a Jacobi form of index zero.  
For example, 
$
\eta \in J_{\frac12,0}(\epsilon) 
$
with the multiplier $\epsilon:\SLtwoZ\to \C$ defined as in \cite{\refer{Iwaniec}, pg\. 45}.  
Since $\eta^{24}=\Delta \in S_1^{12}$ we know that $\epsilon^{24}=1$.  
Let  $H(\Z)=\{ u\smtwomat{1}0{\lambda}1 t\smtwomat0{\mu}{\mu}{\kappa} 
\in \Gamma_{\infty}(\Z): \lambda, \mu,\kappa \in \Z \}$ 
be the integral Heisenberg group.  
Define the character 
$v_H:H(\Z)\to\C$ by $v_H\left( u\smtwomat{1}0{\lambda}1 t\smtwomat0{\mu}{\mu}{\kappa} \right)
=(-1)^{\lambda+\mu +\kappa}$.   Embed $i_{\infty}:\SLtwoZ\to \Gamma_{\infty}(\Z)$ 
via $i_{\infty}(\sigma)=\sigma \oplus I_2$.  
The isomorphism $\Gamma_{\infty}(\Z) \backslash \{\pm I_4\} \cong$ 
$i_{\infty}\left( \SLtwoZ \right) \ltimes H(\Z)$ allows us to view 
$\epsilon^a\,v_H^b$ as a map $\epsilon^a\,v_H^b:\Gamma_{\infty}(\Z)\to \C$ trivial on $\pm I_4$ 
for any integers $a,b\in\Z$ and we use this shorthand.      
For the theta function we have:
$$
\vartheta \in J_{\frac12,\frac12}(\epsilon^3 v_H).    
$$
%
For $\vartheta_d$ defined by $\vartheta_d(\tau,z)=\vartheta(\tau,dz)$ we have
$\vartheta_d \in J_{\frac12,\frac{d^2}2}(\epsilon^3 v_H^d)$, see 
\cite{\refer{GritHulek}}.

\proclaim{Theorem \Equ{c1}} {\rm (Gritsenko, Skoruppa, Zagier) } 
Let $\ell\in\N$ and $t\in \Z$.  
Let $ {\bold{n}}=(n_1, \dots,n_{\ell})\in\Z^{\ell}$ and 
$ {\bold{d}}=(d_1, \dots,d_{\ell})\in\N^{\ell}$.  
Let $n=\sum_{i=1}^{\ell} n_i$ for brevity.  
Define a meromorphic function 
$\TB(t, {\bold{n}}, {\bold{d}}):\Half_1\times\C\to\C$ by 
$$
\TB(t, {\bold{n}}, {\bold{d}})(\tau,z)=
\eta(\tau)^t\,
\prod_{i=1}^{\ell}\, \vartheta(\tau, d_iz)^{n_i}.  
$$
We have $\TB(t, {\bold{n}}, {\bold{d}})\in J_{k,m}^{\text{cusp}}$ if and only if
\roster
\item $2k=t+n$,  
\item $2m= \sum_{i=1}^{\ell} n_i d_i^2$,
\item $t+3n \equiv 0 \mod 24$,    
\item $\forall d\in\N$, $\sum_{i: d | d_i} n_i \ge 0$, 
\item The function $\frac{k}{12} +  \sum_{i=1}^{\ell} n_i {\bar B}_2(d_ix)$ has a positive minimum on $[0,1]$.  
Here ${ B}_2(x)=\frac12x^2-\frac12x +\frac1{12}$ and 
${\bar B}_2(x)=B_2(x-[[x]])$ is the periodic 
extension of its restriction to $[0,1]$.  
\endroster
\endproclaim

We will only avail ourselves of the simplest cases.  
For $k=2$, we always take $t=-6$ and all $n_i=1$ so that 
$n=10$.  With this in mind, let: 
$$
\TB_2(d_1,d_2,\dots,d_{10})(\tau,z)=
\eta(\tau)^{-6}\,
\prod_{i=1}^{10}\, \vartheta(\tau, d_iz).  
$$
For $k=4$, a case used only  incidentally, 
we always take $t=0$ and all $n_i=1$ so that 
$n=8$.  With this in mind, let: 
$$
\TB_4(d_1,d_2,\dots,d_{8})(\tau,z)=
\prod_{i=1}^{8}\, \vartheta(\tau, d_iz).  
$$

In the articles of Gritsenko, see \cite{\refer{Grit1}}, his lift is proved for the 
group $\Gamma[N]=UK(N)U$ where 
$U=u  \smtwomat0110  $.  
We restate his results for cusp forms on the paramodular group $K(N)$.  

\proclaim{Theorem \Equ{c2}} {\rm (Gritsenko) }
Let $\phi \in J_{k,N}^{\text{cusp}}$ and let  
$\phi(\tau,z)= \sum_{n> 0,r\in\Z} c(n,r)q^n\zeta^r$
be the Fourier expansion.   
There is a form $\Grit(\phi)\in S_2^k\left(K(N)\right)$  
given by
$$
\Grit(\phi) \pmatrix \tau & z \\ z & \omega \endpmatrix = 
\sum_{n,r,m}\left( \sum_{\delta \vert (n,r,m)} 
\delta^{k-1}c\left( \dfrac{mn}{\delta^2}, \dfrac{r}{\delta} \right)  \right) 
q^{mN} \zeta^r  e(n\omega).  
$$
For $k$ even, $\Grit(\phi)$ is in the $\mu$-plus space;  
for $k$ odd,  $\Grit(\phi)$ is in the $\mu$-minus space.  
\endproclaim
Thus, the Fourier coefficients of $\Grit(\phi)$ are: 
$$
a\pmatrix mN & r/2 \\ r/2 & n \endpmatrix = 
\sum_{\delta \vert (n,r,m)} \delta^{k-1}c\left( \dfrac{mn}{\delta^2}, \dfrac{r}{\delta} \right) .  
$$

Additionally, V\. Gritsenko and K\. Hulek \cite{\refer{GritHulek}} have developed a lift for 
Jacobi forms with certain characters.  The case of a quadratic character 
will sometimes be useful here.  For odd, squarefree $N \in \N$, 
consider a representation of $N$ as the sum of four squares: 
$N=a^2+b^2+c^2+d^2$, for $a$, $b$, $c$, $d \in \N$.  
We have $\vartheta_a \vartheta_b \vartheta_c \vartheta_d 
\in J_{2, \frac12 N}^{\text{cusp}}(\epsilon^{12}v_H)$.  
The following Theorem shows how to lift such forms to $K(N)$.  
\proclaim{Theorem \Equ{c3}} {\rm (Gritsenko and Hulek) }  
Let $N \in \N$.  
There is a character  $\chi_2^{(N)}: K(N) \to \{ \pm 1\}$.  
Let $\phi \in J_{k,\frac12 N}^{\text{cusp}}$ and let  
$\phi(\tau,z)= \sum_{n,r\in \frac12+\Z: 2Nn>r^2, n>0} c(n,r)q^n\zeta^r$
be the Fourier expansion.   
There is a form $\Grit(\phi)\in S_2^k\left(K(N), \chi_2^{(N)}  \right)$  
given by
$$
\Grit(\phi) \pmatrix \tau & z \\ z & \omega \endpmatrix = 
\sum_{n,r,m \text{ odd } \in \N}\left( \sum_{\delta \vert (n,r,m)} 
\delta^{k-1}c\left( \dfrac{mn}{2\delta^2}, \dfrac{r}{2\delta} \right)  \right) 
q^{\frac{mN}{2}} \zeta^{\frac{r}2}  e(\frac{n}2 \omega).  
$$
\endproclaim
The point is that for $i \in \{1,2\}$ and 
$\phi_i \in J_{k_i,\frac12 N}^{\text{cusp}}$, 
we have a paramodular form with trivial character 
$\Grit(\phi_1)\Grit(\phi_2) \in S_2^{k_1+k_2}\left( K(N) \right)$.

{\bf{\newsection{sec4}{   Vanishing Theorems and Congruences.}} }

Computations often require apriori sets of Fourier coefficients 
that determine linear dependence among Siegel modular forms.  
We discuss such sets in degree two and prove that some also 
determine congruences among Fourier coefficients.  
The forms that index 
Fourier coefficients are a partially ordered set with no 
natural linear order.  Indeed, the intrinsic measure of the
vanishing order of a Fourier series is the closure in $\PsR$ 
of the convex ray hull of its support.  
For computational purposes however, ordering the support with 
a convex function $\phi$ is a versatile expedient.  We review the
results of \cite{\refer{PoorYuenExt}} and \cite{\refer{PoorYuenComp}}.  

\definition{Definition \Equ{D1}}
A function $\phi:\PsR\to\R_{\ge 0}$
is called
{\it{\classone}\/} if
\roster
\item For all $s \in \PR$, $\phi(s) > 0$,  
\item for all
$\lambda\in\R_{\ge0}$ and $s\in\PsR$, $\phi(\lambda s)=\lambda\phi(s)$, 
\item for all $s_1,s_2\in\PsR$,  $\phi(s_1+s_2)\ge\phi(s_1)+\phi(s_2)$.
\endroster
\enddefinition
\Classone\ functions are continuous on $\PR$ and 
respect the partial order on $\PsR$.  
Basic examples are:  
For $s\in\PsR$, define
\roster
\item $m(s) = \inf_{u\in\Z_n\backslash\{0\}}u'su$, the Minimum function, 
\item $\redtr(s) = \inf_{u\in\GLnZ}\tr(u'su)$, the reduced trace, 
\item $\delta(s) = \det(s)^{1/n}$, the reduced determinant, 
\item $w(s) = \inf_{u\in\PR}{\<u,s\>\over m(u)}$, the dyadic trace.
\endroster
For $n=2$, the dyadic trace of a 
Minkowski reduced $s=\smtwomat{a}{b}{b}{c}\in\PtwoR$ 
is given by 
$w(s)=a+c-|b|$, see \cite{\refer{PoorYuenRel}}.  For $n=2$, 
Minkowski reduced means $2\vert b \vert \le a \le c$.  

\proclaim{ \Equ{D5} Vanishing Theorem}
Let $\phi$ be \classone.  
For all $n \in \Z^{+}$ there exists a $c_n(\phi)\in \R_{>0}$ such that: 
For any subgroup $\Gamma \subseteq \Gamma_n$ with finite index $I$ and 
coset decomposition $\Gamma_n=\cup_{i=1}^I \Gamma M_i$, we have  
$$
\forall\,k\in\Z^{+},\,\forall\, f\in S_n^k(\Gamma), \quad
\frac1{I} \sum_{i=1}^I 
\inf \phi\left(\supp(f| M_i)\right) > c_n(\phi)\,k \implies f \equiv 0. 
\tag \Equ{D6f}
$$
For $n=2$, we may take $c_n(\phi)= 
\inf\phi\left([\frac1{30}\smtwomat{3}{1}{1}{3} ] \right)$.  
\endproclaim
\demo{Proof}
This is Theorem~{2.5} from \cite{\refer{PoorYuenS22}}, except for the last comment, 
which is Corollary~{5.8} from \cite{\refer{PoorYuenExt}}.   
\enddemo

We obtain our Vanishing Theorem for automorphic forms on 
$K(p)$ by viewing them as automorphic with respect to 
$\Gamma_0'(p)=K(p) \cap \Gamma_2$.  
To use Theorem~\refer{D5} we need coset representatives
$\Gamma_2 = \cup_{i=1}^{(1+p)(1+p^2)}\Gamma_0'(p)Y_i$;  
these may be found in \cite{\refer{IbukHash}, pg\. 71}.  

\proclaim{ \Equ{D6} Theorem}{\rm (Hashimoto, Ibukiyama)\/}
As a complete set of $1+p+p^2+p^3$ representatives of the 
coset space $\SptwoZ/ \Gamma_0'(p)$ we may take
$$
X_1(a,b,c)=
\pmatrix 
1 & 0 & 0 & 0 \\ 
a & 1 & 0 & 0 \\ 
b & c & 1 & -a \\
c & 0 & 0 & 1 \endpmatrix;\quad 
X_2(a,b)=
\pmatrix 
0 & 1 & 0 & 0 \\ 
1 & 0 & 0 & 0 \\ 
a & 0 & 0 & 1 \\
b & a & 1 & 0 \endpmatrix;
$$
$$
X_3(a)=
\pmatrix 
0 & -a & -1 & 0 \\ 
0 & 1 & 0 & 0 \\ 
1 & 0 & 0 & 0 \\
a & 0 & 0 & 1 \endpmatrix ;\quad 
X_4=
\pmatrix 
0 & 0 & 0 & -1 \\ 
0 & 0 & -1 & 0 \\ 
0 & 1 & 0 & 0 \\
1 & 0 & 0 & 0 \endpmatrix
$$
where $a$, $b$, $c$ run over the integers modulo $p$.  
\endproclaim
%
%
It is helpful to abbreviate
$$
\kap=u\left( \matrix 0 & \frac{1}{p} \\ -1 & 0 \endmatrix \right)=
\pmatrix 
0 & \frac{1}{p} & 0 & 0 \\ 
-1 & 0 & 0 & 0 \\ 
0 & 0 & 0 & p \\
0 & 0 & -1 & 0 \endpmatrix.  
$$ 
\proclaim{ \Equ{D7} Corollary}
A complete set of $1+p+p^2+p^3$ coset representatives $Y$ for 
$\Gamma_0'(p)Y\in \Gamma_0'(p) \backslash \SptwoZ $ 
and a representative from $\Delta_2(\Q)$ for 
$K(p)Y $ and a representative, $I_4$ or $\kap$, for $K(p)Y \Delta_2(\Z)$
is given by 
$$
X_4\inv=
\pmatrix 
0 & 0 & 0 & 1 \\ 
0 & 0 & 1 & 0 \\ 
0 & -1 & 0 & 0 \\
-1 & 0 & 0 & 0 \endpmatrix 
\in K(p) 
\pmatrix 
0 & \frac{1}{p} & 0 & 0 \\ 
1 & 0 & 0 & 0 \\ 
0 & 0 & 0 & p \\
0 & 0 & 1 & 0 \endpmatrix
\subseteq K(p) \kap\, \Delta_2(\Z);    
$$
$$
X_3(a)\inv=
\pmatrix 
0 & 0 & 1 & 0 \\ 
0 & 1 & 0 & 0 \\ 
-1 & -a & 0 & 0 \\
0 & 0 & -a & 1 \endpmatrix 
\in K(p) 
\pmatrix 
-\frac{1}{p} & -\frac{a}{p} & 0 & 0 \\ 
0 & -1 & 0 & 0 \\ 
0 & 0 & -p & 0 \\
0 & 0 & a & -1 \endpmatrix
\subseteq K(p) \kap\, \Delta_2(\Z),  
$$
where $a$ runs over the integers modulo $p$;    
$$
X_2(a,b)\inv=
\pmatrix 
0 & 1 & 0 & 0 \\ 
1 & 0 & 0 & 0 \\ 
-a & -b & 0 & 1 \\
0 & -a & 1 & 0 \endpmatrix 
\in K(p) 
\pmatrix 
-{\hat b}a & -\frac{1}{p} & 0 & \frac{\hat b}{p} \\ 
-1 & 0 & 0 & 0 \\ 
0 & 0 & 0 & -p \\
0 & 0 & -1 & -{\hat b}a \endpmatrix
\subseteq K(p) \kap\, \Delta_2(\Z),  
$$
where $a$, $b$, $\hat b$ run over the integers modulo $p$ 
with $b{\hat b} \equiv 1 \mod p$;     
$$
X_2(a,0)\inv=
\pmatrix 
0 & 1 & 0 & 0 \\ 
1 & 0 & 0 & 0 \\ 
-a & 0 & 0 & 1 \\
0 & -a & 1 & 0 \endpmatrix 
\in K(p) 
\pmatrix 
-\frac{1}{p}  & {\hat a} & 0 & \frac{\hat a}{p} \\ 
0 & 1 & -{\hat a} & 0 \\ 
0 & 0 & -p & 0 \\
0 & 0 & -{\hat a}p & 1 \endpmatrix
\subseteq K(p) \kap\, \Delta_2(\Z),  
$$
where $a$, $\hat a$ run over the integers modulo $p$ 
with $a{\hat a} \equiv 1 \mod p$; 
$$
X_2(0,0)\inv=
\pmatrix 
0 & 1 & 0 & 0 \\ 
1 & 0 & 0 & 0 \\ 
0 & 0 & 0 & 1 \\
0 & 0 & 1 & 0 \endpmatrix 
\in  \Delta_2(\Z);    
$$
$$
X_1(a,b,c)\inv=
\pmatrix 
1 & 0 & 0 & 0 \\ 
-a & 1 & 0 & 0 \\ 
-b & -c & 1 & a \\
-c & 0 & 0 & 1 \endpmatrix 
\in K(p) 
\pmatrix 
-\frac{1}{p} & -\frac{{\hat b}c}{p} & \frac{\hat b}{p} & \frac{{\hat b}a}{p} \\
0 & 1 & -a{\hat b} & 0 \\ 
0 & 0 & -p & 0 \\
0 & 0 & -{\hat b}c & 1 \endpmatrix
\subseteq K(p) \kap\, \Delta_2(\Z),  
$$
where $a$, $b$, $\hat b$, $c$ run over the integers modulo $p$ 
with $b{\hat b} \equiv 1 \mod p$;  
$$
X_1(a,0,c)\inv=
\pmatrix 
1 & 0 & 0 & 0 \\ 
-a & 1 & 0 & 0 \\ 
0 & -c & 1 & a \\
-c & 0 & 0 & 1 \endpmatrix 
\in K(p) 
\pmatrix 
0 & -\frac{1}{p} & \frac{\hat c}{p} & \frac{{\hat c}a}{p} \\
1 & 0 & 0 & -{\hat c} \\ 
0 & 0 & 0 & -p \\
0 & 0 & 1 & 0 \endpmatrix
\subseteq K(p) \kap\, \Delta_2(\Z),  
$$
where $a$, $c$, $\hat c$  run over the integers modulo $p$ 
with $c{\hat c} \equiv 1 \mod p$; 
$$
X_1(a,0,0)\inv=
\pmatrix 
1 & 0 & 0 & 0 \\ 
-a & 1 & 0 & 0 \\ 
0 & 0 & 1 & a \\
0 & 0 & 0 & 1 \endpmatrix 
\in  \Delta_2(\Z),  
$$
where $a$  runs over the integers modulo $p$.  
\endproclaim
\demo{Proof}
First of all, it is the $Y=X_i\inv$ from Theorem~\refer{D6} 
that give left coset representatives. 
Finally, 
one can check these assertions directly by taking inverses
and multiplying.  In the case of $X_2(a,b)\inv$, 
for example, the following element is in $K(p)$:  
$$
\pmatrix 
0 & 1 & 0 & 0 \\ 
1 & 0 & 0 & 0 \\ 
-a & -b & 0 & 1 \\
0 & -a & 1 & 0 \endpmatrix
\pmatrix 
-{\hat b}a & -\frac{1}{p} & 0 & \frac{\hat b}{p} \\ 
-1 & 0 & 0 & 0 \\ 
0 & 0 & 0 & -p \\
0 & 0 & -1 & -{\hat b}a \endpmatrix\inv=
\pmatrix 
-p & a{\hat b} & -\frac{\hat b}{p}  & 0 \\ 
0 & -1 & 0 & 0 \\ 
bp & a(1-b{\hat b}) & \frac{b{\hat b}-1}{p} & 0 \\
ap & -a^2{\hat b} & 0 & -1 \endpmatrix.  
$$
Furthermore, the following element is in $\Delta_2(\Z)$:
$$
\pmatrix 
0 & \frac{1}{p} & 0 & 0 \\ 
-1 & 0 & 0 & 0 \\ 
0 & 0 & 0 & p \\
0 & 0 & -1 & 0 \endpmatrix\inv 
\pmatrix 
-{\hat b}a & -\frac{1}{p} & 0 & \frac{\hat b}{p} \\ 
-1 & 0 & 0 & 0 \\ 
0 & 0 & 0 & -p \\
0 & 0 & -1 & -{\hat b}a \endpmatrix=
\pmatrix 
1 & 0 & 0 & 0 \\ 
-a{\hat b} & -1 & 0 & {\hat b} \\ 
0 & 0 & 1 & -a{\hat b} \\
0 & 0 & 0 & -1 \endpmatrix.  \qed
$$
\enddemo

With these coset representatives, 
we may use Theorem~\refer{D5} to prove: 
\proclaim{ \Equ{D8} Theorem}
Let $f\in S_2^k\left( K(p) \right)$.  
Let $\phi$ be a \classone\ $\GLtwoZ$-class function.  
Unless $f \equiv 0$, we have 
$$
\inf \phi\left(\supp(f) \right)
+p \inf \phi\left(\supp(f|\mu) \right)\le 
\phi\left(\tfrac1{30}\smtwomat{3}{1}{1}{3}  \right)k(1+p^2).  
$$
\endproclaim
\proclaim{ \Equ{D9} Corollary}
For nontrivial $f\in S_2^k\left( K(p) \right)$ we have   
$$
\min \delta\left(\supp(f) \right)\le \frac{\sqrt{2}}{15}
k\,\frac{1+p^2}{1+p}.   
$$
If $f$ is additionally a $\mu$-eigenform then we have 
$$
\align
\min w\left(\supp(f) \right)
&\le \frac{k}{6}
\frac{1+p^2}{1+p},  \\ 
\min \redtr\left(\supp(f) \right)
&\le \frac{k}{5}
\frac{1+p^2}{1+p},  \\ 
\min m\left(\supp(f) \right)
&\le \frac{k}{10}
\frac{1+p^2}{1+p}.   
\endalign
$$
\endproclaim

\demo{Proof of Theorem~\refer{D8} and Corollary~\refer{D9}}
Let 
$\Gamma_2 = \cup_{i=1}^{I}\Gamma_0'(p)Y_i$ 
with $I=(1+p)(1+p^2)$.  For nontrivial $f$, Theorem~\refer{D5} gives 
$$
\frac1{I} \sum_{i=1}^I 
\inf \phi\left(\supp(f| Y_i)\right) \le  c_2(\phi)\,k . 
\tag \Equ{D6e}
$$
For each $i$ we have 
$Y_i\in K(p) \kap^{\epsilon_i} 
\pmatrix u_i & * \\ 0 & u_i^* \endpmatrix$ 
with $u_i\in \GLtwoZ$, with $\epsilon_i=0$ in $1+p$ cases 
and with $\epsilon_i=1$ in $p^2+p^3$ cases by 
Corollary~\refer{D7}.  
Since $\phi$ is a class function and 
$\supp\left( f|Y_i \right)=
u_i'\kap'{}^{\epsilon_i}\supp(f)\kap^{\epsilon_i} u_i$, 
we have 
$\phi\left(\supp\left( f|Y_i \right)\right)=
 \phi\left(\supp\left( f|\kap^{\epsilon_i}  \right)\right) $.  
%
When $\epsilon_i=0$ this is 
$\phi\left(\supp\left( f\right)\right)$ 
and when $\epsilon_i=1$ this is
$ \phi\left(\supp\left( f|\kap  \right)\right) 
= \phi\left(\frac{1}{p}\supp\left( f|\mu \right)\right)$ by Lemma~\refer{B8}.  
In view of these two cases, 
condition~\refer{D6e} becomes 
$$
\frac1{I} \left(  
(1+p)\inf \phi\left(\supp(f)\right)+
(p^2+p^3)
\inf  \phi\left(\frac{1}{p}
\supp\left( f|\mu \right)\right)\right)\le  c_2(\phi)\,k, 
$$
which is equivalent to the conclusion of Theorem~\refer{D8}.  

In the Corollary, we have replaced $\inf$ by $\min$ 
because
the infimum of $\phi\left(\supp\left( f|Y_i \right)\right)$
is attained whenever $\phi$ has finite shells or takes values in a 
lattice.  If $f$ is a $\mu$-eigenform then 
$\supp(f|\mu)=\supp(f)$ and the conclusions for $w$, $m$ and $\redtr$
follow upon evaluating   
$w\left(1/{30}\smtwomat{3}{1}{1}{3}  \right)=1/6$, 
$\redtr\left(1/{30}\smtwomat{3}{1}{1}{3}  \right)=1/5$ and
$m\left(1/{30}\smtwomat{3}{1}{1}{3}  \right)=1/10$.   
In general, twinning does not change the determinant, 
so $\delta\left(\supp(f|\mu)\right)=\delta\left(\supp(f)\right)$ 
and we note that 
$\delta\left(1/{30}\smtwomat{3}{1}{1}{3}  \right)=\sqrt{2}/15$.   
  \qed
\enddemo

The goal of this section is to give similar sets of 
Fourier coefficients that determine congruences.    
The essential difficulty is at the primes two and three.  
Results of this type depend upon the work of J\. I\. Igusa  
\cite{\refer{IgusaMFZ}}.  
Let $R$ be a ring.  For $A \subseteq \C$, let 
$$
M_n^k\left( \Gamma\right)(A) = 
\{f\in M_n^k\left( \Gamma\right): 
\forall\, T \in \PsQ, a(T;f)\in A\}.  
$$
We see that $M_n^k\left( \Gamma\right)(A)$ is an $R$-module, 
or $M_n\left( \Gamma\right)(A)$ a ring, whenever $A$ is.   
Further, for a prime $\ell$, let 
$M_n^k\left( \Gamma\right)(\F_{\ell}) $ denote the 
reduction modulo $\ell$ of the coefficients of the 
Fourier series from $M_n^k\left( \Gamma\right)(\Z) $.  
Let $R_{\ell}: M_n^k\left( \Gamma\right)(\Z)\to 
M_n^k\left( \Gamma\right)(\F_{\ell}) $ be the natural reduction map.  

Using Igusa's work we can prove the following: 
%
\proclaim{ \Equ{D10} Theorem} 
Let $K\subseteq \C$ be a number field, 
${\Cal O}$ its ring of integers and 
$\ideal$ a prime ideal in ${\Cal O}$.  
Let $ f \in S_2^k \left( K(p)\right)({\Cal O})$
be a $\mu$-eigenform with $a(T;f) \in \ideal $ for all $  T\in {}^p\X2$ satisfying 
$  w(T) \le \frac{k}{6} \frac{p^2+1}{p+1}$.  
Then we have $a(T;f) \in \ideal$ for all $  T\in {}^p\X2$.  
\endproclaim
The proof we will give is valid only for primes $p$ because we rely on the coset 
decomposition of Corollary~{\refer{D7}}.  
The proof of Theorem~\refer{D10} is at the end of this section.  
These results partially generalize Sturm's Theorem \cite{\refer{Sturm}} on elliptic 
modular forms for $\Gamma_0(N)$, which we state for level one in the next Theorem.    
For a ring $R \subseteq \C$ and a set $S  \subseteq \C$, 
let $R\<S\>$ denote the $R$-module generated by finite 
$R$-linear combinations of elements from $S$.  
Recall the notation $ \{\{k\}\}=\dim S_1^k$.   
\proclaim{ \Equ{D13} Theorem  }
If $f \in M_1^k\left( \C \right)$,  
then 
$f \in M_1^k\left( \Z\<a(j;f):j \le \{\{k\}\} \>\right)$.  
\endproclaim
Theorem~\refer{D10} is proven by reduction to a congruence criterion for 
integral forms of level one 
which is of independent interest. 
For $g=2$, we will show in Theorem~{5.15} that
$$
\text{ If $f \in M_2^k\left( \C \right)$,   
then 
$f \in M_2^k\left( \Z\<a(T;f):w(T) \le \frac{k}{6} \>\right)$. }
\tag\Equ{D18}
$$
For $k \ne 2$, each $M_1^k(\Z)$ has a basis 
whose $i^{\underline{\text{th}}}$ element has a Fourier expansion 
$q^i +0(q^{i+1})$ for $0 \le i \le \dim S_1^k$.   
An immediate consequence is
\proclaim{ \Equ{D14} Theorem }
If $f \in \Sym\left(M_1\otimes M_1 \right)^k(\C)$, 
then we have 
$$f \in \Sym\left(M_1\otimes M_1 \right)^k\left( \Z\<a(i,j;f):i,j \le \{\{k\}\} \>\right).$$  
\endproclaim
For a ring $R \subseteq\C$, this  shows that 
$\Sym\left(M_1\otimes M_1 \right)^k(R)=
\Sym\left(M_1^k(R)\otimes M_1^k(R) \right)$.  
Interest in the intermediate ring $\Sym\left(M_1\otimes M_1 \right)^k(R)$ 
stems from the Witt map \cite{\refer{Witt}} 
$$
\align
W: M_2^k(R) &\to \Sym\left(M_1\otimes M_1 \right)^k(R)  \\
\left( \Omega \mapsto f(\Omega) \right) &\mapsto 
\left( (\tau_1,\tau_2) \mapsto f\left(\matrix \tau_1 & 0 \\ 0 & \tau_2 \endmatrix \right) \right) 
\endalign
$$
for which the Fourier coefficients obey 
$
a(i,j;W(f))=\sum_{b} a\left(\pmatrix i & b \\ b & j \endpmatrix ;f \right). 
$
The following exact sequence is often the basis of computing in the ring 
$M_2(\C)$.  
$$
0 \to M_2^{k-10}(\C)\, {\overset {\cdot X_{10}}\to {\to}} M_2^{k}(\C)
\,{\overset {W }\to {\to}} \Sym\left(M_1\otimes M_1 \right)^k(\C) \to 0.  
$$
Igusa showed that the following sequence is exact 
\cite{\refer{IgusaMFZ}, Lemma~7, p\. 163}.   
$$
0 \to M_2^{k-10}(\Z)\, {\overset {\cdot X_{10}}\to {\to}} M_2^{k}(\Z)
\,{\overset {W }\to {\to}} \Sym\left(M_1\otimes M_1 \right)^k(\Z) 
$$
but that the Witt map from 
$M_2^{k}(\Z)$ to $\Sym\left(M_1\otimes M_1 \right)^k(\Z) $ does {\sl not\/} surject 
when twelve divides $k$.  
The following Lemma has greater applicability if we note that 
$M_2^k(R)= M_2^k(\Z) \otimes_{\Z} R$ for $\Z$-modules $R$, cf\. \cite{\refer{IgusaMFZ}, pg\. 150}.  
%
\proclaim{ \Equ{D15} Lemma {\rm ( \cite{\refer{IgusaMFZ}}, Lemma~{13}, pg\. 171)} }
We have 
$W\left( M_2^{k}(\Z) \right) = \Sym\left(M_1\otimes M_1 \right)^k(\Z) $ 
if and only if  $k \not\equiv 0 \mod 12$.  
\endproclaim

We need a lemma that records the reach of each Fourier expansion.  
Some familiarity with the geometry of numbers, as in \cite{\refer{PoorYuenExt}}, 
is required for the proof of the next Lemma.  In particular, let 
$\mc(A)=\R_{\ge 0}\< x\in \Z^2\setminus\{0\}: x'Ax=m(A)\>$ be the cone in $\PtwosR$ generated 
by the minimal vectors of $A$.  
Also, for $X \subseteq \PtwosR$, $\SH(X)$ denotes the closure in $\PtwosR$ 
of the convex ray hull of $X$.  Finally, we need to mention the technique of 
restriction to modular curves \cite{\refer{PoorYuenComp}} in the simplest case.  
For $s \in \PZ$, define $\phi_s:     \Half_1 \to \Half_n$ by 
$\tau \mapsto s\tau$ and note that, for any $\ell \in \N$ with $\ell s\inv\in \PZ$, 
the pullback is a ring homomorphism  
$\phi_{s}^*: M_n \to M_1\left( \Gamma_0(\ell) \right)$ 
that multiplies weights by~$n$.   
If $f \in M_n^k$ has a Fourier series 
$\FS_n(f)=\sum_{T} a(T) q^T$, then 
$\FS_1\left( \phi_{s}^{*}f\right)=\sum_{j}
\(  \sum_{T: \<s,T\>=j  }  a(T)  \)  q^j$.  
\proclaim{ \Equ{D17}  Lemma }
Let $R \subseteq \C$ be a $\Z$-module.  Let $b{>0}$.  
Let $f \in M_2^k(\C)$ and $f' \in M_2^{k-10}(\C)$ with $f=X_{10}f'$.  
Then we have 
$$
\left(  \forall\, T\in \Xtwos: w(T) < b, \,a(T;f) \in R \right)
\implies 
\forall\, T\in \Xtwos: w(T) < b-\frac{3}{2},\, a(T;f') \in R.   
$$
\endproclaim
\demo{Proof}
Let $\supp_R(f')=\{T\in \Xtwos:  a(T;f') \not\in R \}$.  
Our goal is to prove the inequality
$\min w\left( \supp_R(f') \right) \ge b-3/2$.  
Letting 
$K=\SH(\supp_R(f'))$, it is equivalent to show that 
$\min w\left(K \right) \ge b-3/2$.   
Suppose, by way of contradiction, that $T\in \supp_R(f')$  
is a vertex of $K$ with minimal dyadic trace $w(T) < b-3/2$.  
By changing representatives within the $\GLtwoZ$-equivalence class 
of $T$, we may assume $T\in \mc(A)$ for $A=\smtwomat{1}{-1/2}{-1/2}{1}$.  

In the dual 
$K^\sqcup= \{y\in \PtwosR: \forall x\in K, \<x,y\> \ge 1 \}$,  
the vertex $T$ corresponds to the convex two dimensional face 
$F=\{y \in K^\sqcup: \<T,y\>=1\}$.  The point $Y=A/w(T)$ is in this 
face and also in the interior of the cone $\mc({\hat A})$ for 
${\hat A} =\smtwomat{1}{1/2}{1/2}{1}$.  
Let $N $ be a neighborhood of $Y$ with $N \subseteq \mc({\hat A})$.  
For any $H \in N \cap F^{\circ}$, we have both 
\therosteritem1  
$\{x\in K: \<x,H\>=1\}=\{T\}$ 
because $H\in F^{\circ}$ and 
\therosteritem2   
$\{x\in [{\hat A}]: \<x,H\>=w(H)\}=\{{\hat A}\}$ 
because $H\in \mc({\hat A})^{\circ}$.  
Consider the continuous function 
$\< {\hat A} +T, \cdot \>/m(\cdot)$ evaluated at $Y$:  
$$
\frac{\< {\hat A} +T, Y \>}{m(Y)}= 
\frac{\< {\hat A} +T, A \>}{m(A)}= 
 \frac{\< {\hat A} , A \>}{m(A)}+\frac{\<  T, A \>}{m(A)}=
\frac32 +w(T) < 
\frac32 +\left(b-\frac32\right)=b.  
$$
As $Y$ is an accumulation point of $N \cap F^{\circ}$, 
there is $H \in N \cap F^{\circ}$ with 
${\< {\hat A} +T, H \>}/{m(H)} < b$.  

The leading term of $\phi_H^* X_{10}$ is 
$a({\hat A};X_{10})q^{\<{\hat A} , H\>}=q^{\<{\hat A} , H\>}$ by \therosteritem2.  
The leading term of $\phi_H^* f'$ 
modulo $R[[q]]$ is 
$ a(T;f')q^{\<T , H\>}$ by \therosteritem1.  
Thus $\phi_H^* f=\left(\phi_H^* X_{10} \right)\left(\phi_H^* f' \right)$ 
has a leading term modulo $R[[q]]$ given by 
$ a(T;f')q^{\<{\hat A} +T , H\>}$ and hence 
$$
\sum_{X \in \supp(f): \<X,H\>=\<{\hat A} +T , H\>} 
a(X;f) \equiv a(T;f') \quad\text{ modulo  } R.  
$$
For $X \in \supp(f)$ with $ \<X,H\>=\<{\hat A} +T , H\>$ we have 
$$
w(X) \le \frac{\<X,H\>}{m(H)}= 
\frac{\<{\hat A} +T,H\>}{m(H)}< b,  
$$
so that $a(X;f)\in R$ by hypothesis and we obtain 
$a(T;f')\in R$ and the contradiction $T \not\in \supp_R(f')$.  \qed
\enddemo
\proclaim{ \Equ{D19} Theorem }
Let $f \in M_2^k\left( \C \right)$,   
then 
$f \in M_2^k\left( \Z\<a(T;f):w(T) \le \frac{k}{6} \>\right)$.   
Furthermore, 
$f \in M_2^k\left( \Z\<a(T;f):w(T) \le \dfrac{k-\nu(k)}{6} \>\right)$ 
for $\nu(k)=0,1,2,3,4,5$ for $k \equiv 0,10,8,6,4,2 \mod 12$, respectively.  
\endproclaim
\demo{Proof}
It suffices to prove the second assertion 
and we do this by induction on~$k$.  
Let $A=\Z\<a(T;f):w(T) \le  (k-\nu(k))/{6} \>$.  
Let $f \in M_2^k\left( \C \right)$ and consider 
${\hat f}=f$ for $k \not\equiv 0\mod 12$ and 
${\hat f}=E_4f$ for $k \equiv 0\mod 12$.  
Thus, in all cases, the weight ${\hat k}$ of ${\hat f}$ is not divisible by $12$.  
The Witt image 
$W({\hat f}) \in \Sym\left(M_1\otimes M_1 \right)^{\hat k}(\C) $ has Fourier coefficients 
$$
a(i,j;W({\hat f}))=\sum_{b} a\left(\pmatrix i & b \\ b & j \endpmatrix ;{\hat f} \right). 
$$
By Theorem~\refer{D14}, we have 
$W({\hat f})\in \Sym\left(M_1\otimes M_1 \right)^{\hat k}
\left( \Z\<a(i,j;W({\hat f})):i,j \le \{\{{\hat k}\}\} \>\right)$. 
For $i,j \le \{\{{\hat k}\}\}$, we have 
$w\smtwomat{i}bbj  \le \tr\smtwomat{i}bbj = i+j \le 2\{\{{\hat k}\}\}$, 
so that $a(i,j;W({\hat f}))\in A$ and 
$W({\hat f}) \in \Sym\left(M_1\otimes M_1 \right)^{\hat k}(A)$ 
if $   2\{\{{\hat k}\}\} \le (k-\nu(k))/6$.  
This inequality holds in all six cases, see Table~3.  
By Lemma~\refer{D15}, 
there is an $F \in M_2^{{\hat k}}(A)$ 
such that  $W(F)=W({\hat f})$. Therefore, $F-{\hat f}=X_{10}f'$ for some 
$f'\in M_2^{k'}\left( \C \right)$ with $k'={\hat k}-10$.   For all $T$ with 
$w(T)\le (k-\nu(k))/6$ we have $a(T;F-{\hat f})\in A$.  
Since $ (k-\nu(k))/6 \in \frac12\Z$ and since $w\left(\Xtwos\right)\subseteq \frac12\Z$, 
the strict inequality $w(T)< (k-\nu(k))/6+1/2$ is equivalent to the 
inequality $w(T)\le (k-\nu(k))/6$.  
By Lemma~\refer{D17}, we have $a(T;f')\in A$ for 
all $T$ with $w(T) < (k-\nu(k))/6-1$, or equivalently, 
$w(T) \le (k-\nu(k))/6-3/2$.    

By the induction hypothesis we have 
$f'\in M_2^{k'}
\left( \Z\<a(T;f'):w(T) \le (k'-\nu(k'))/6 \>\right)$.  
Thus we have $f'\in M_2^{k'}
\left(  A \right)$ since 
$ (k'-\nu(k'))/6 \le(k-\nu(k))/6-3/2$, an inequality 
that holds with equality in all six cases.   
Since $X_{10}\in S_2^{10}(\Z)$, we have 
$F-{\hat f}=X_{10}f'\in  S_2^{\hat k}
\left( A \right)$ and we also have 
${\hat f}\in  M_2^{\hat k}\left( A \right)$.  
This is either ${  f}\in  M_2^{  k}\left( A \right)$ or 
${E_4  f}\in  M_2^{  k+4}\left( A \right)$, from which 
${  f}\in  M_2^{  k}\left( A \right)$ follows.  
It remains to check the base case of the induction.  
It suffices to note that for $k<10$, nontrivial $M_2^k(\C)$ are spanned by 
one Eisenstein series.   \qed
\enddemo
The following Table is an aid to checking the proof of Theorem~\refer{D19}.  
\smallskip
\centerline{Table 3.  Functions from the proof of Theorem~\refer{D19}.} \nopagebreak
\smallskip
\begintable
 1\| $k \mod 12$ | 0 | 10 | 8 | 6 | 4 | 2  \cr
2
\| $\nu(k)$ | 0 | 1 | 2 | 3 | 4 | 5   \cr
3
 \| ${\hat k}$ | $k+4$ | $k$ | $k$ | $k$ | $k$ | $k$ \cr
4
 \| ${ k'}$ | $k-6$ | $k-10$ | $k-10$ | $k-10$ | $k-10$ | $k-10$ \cr
5
 \| $\nu(k')$ | $3$ | $0$ | $1$ | $2$ | $3$ | $4$ \cr
6
 \| $12\{\{k\}\}$ | $k$ | $k-10$ | $k-8$ | $k-6$ | $k-4$ | $k-14$ \cr
7
 \| $k-\nu(k)$ | $k$ | $k-1$ | $k-2$ | $k-3$ | $k-4$ | $k-5$ \cr
8
 \| $k'-\nu(k')$ | $k-9$ | $k-10$ | $k-11$ | $k-12$ | $k-13$ | $k-14$ 
\endtable\smallskip
Theorem~\refer{D19} is a module criterion for 
integral forms of level one.  
To prove a congruence criterion for $K(p)$, we 
use the explicit coset representatives from Corollary~\refer{D7}.  
%
\demo{Proof of Theorem~\refer{D10}  }
Write the Fourier expansion of $f$ as 
$\FS(f)=\sum_T a(T)q^T$.  
For each coset $\Gamma_0'(p)Y\in \Gamma_0'(p)\backslash\Gamma_2$ 
as in Corollary~\refer{D7}, we consider the Fourier expansion of $f|Y$.  
There are $I=1+p+p^2+p^3$ of these cosets.  
These $Y$ break down into $p+1$ cases of one type where 
$K(p)Y\Delta_2(\Z)=K(p)\Delta_2(\Z)$ and 
$p^2+p^3$ cases of another type where  $K(p)Y\Delta_2(\Z)=K(p)\kap\Delta_2(\Z)$.  
For the first type we have 
$$
Y=\delta=
\pmatrix U & X U^* \\ 0 & U^* \endpmatrix \in \Delta_2(\Z)
$$
for $U\in \GLtwoZ$ and $X\in M_{2\times 2}^{\text{sym}}(\Z)$,  
so that for $f_Y= \det(U)^k f$ 
we have the Fourier expansion
$$
\FS(f_Y)=\sum_{T\in {}^p\X2 } a(T)q^{T[U]}.    
$$
For the second type we have 
$Y=\kap\delta$.   
Using $\kap=\mu \frac1{\sqrt{p}} \smtwomat{I}00{pI}$ along with 
our assumption $f|\mu = \pm f$ and letting $f_Y= \pm\det(U)^k p^k f$ 
and $\zeta_p=e(\frac1p)$, 
we have the Fourier expansion
$$
\FS(f_Y)=\sum_{T\in {}^p\X2 } a(T)\zeta_{p}^{\<X,T \>}q^{\frac1p T[U]}.    
$$
Each of these Fourier expansions has coefficients in ${\Cal O}[\zeta_{p}]$  
and furthermore in $\ideal[\zeta_{p}]$ for $T$ such that 
$  w(T) \le \frac{k}{6} \frac{p^2+1}{p+1}$.  
We consider the product of these series 
$$
F=\left(
\prod_{\Gamma_0'(p)Y\in \Gamma_0'(p)\backslash\Gamma_2}f_Y  \right)
\in M_2^{kI}\left( {\Cal O}[\zeta_{p}]
\right). 
$$
From the product for $F$ we see that 
$\supp_{\ideal[\zeta_{p}]}(F) \subseteq 
\sum_Y \supp_{\ideal[\zeta_{p}]}(f_Y)$ and so 
$$
\align
&\min\  w\left( \supp_{\ideal[\zeta_{p}]}(F) \right)  
\ge 
\min\ w\left( \sum_Y \supp_{\ideal[\zeta_{p}]}(f_Y) \right) 
\ge 
\min \left( \sum_Y w\left(\supp_{\ideal[\zeta_{p}]}(f_Y) \right) \right) \\
&\ge 
 \sum_Y \min\ w\left(\supp_{\ideal[\zeta_{p}]}(f_Y) \right) 
> (p+1) \frac{k}{6} \frac{p^2+1}{p+1} + \frac1p(p^2+p^3) \frac{k}{6} \frac{p^2+1}{p+1}
=\frac{k}{6}.  
\endalign
$$
By Theorem~\refer{D18}, we have 
$F\in S_2^{kI}\left( \ideal[\zeta_{p}]\right)$.  

Now we will show that 
$\FS(f)=\sum_T a(T)q^T$ has $a(T) \in \ideal$.  
We proceed by contradiction; if not, 
there is a $T_0\in {}^p\X2$ with $a(T_0)\not\in \ideal$.  
Pick a prime ideal $L$ containing 
$\ideal {\Cal O}_{K(\zeta_p)}$ in $ {\Cal O}_{K(\zeta_p)}$.  
Denote by 
$$
R_{L}:M_2\left( {\Cal O}_{K(\zeta_p)}\right)\to 
M_2\left( {\Cal O}_{K(\zeta_p)}/L\right)
$$
the reduction of Fourier coefficients modulo $L$.  
Then $R_{L}(F)=0$ but we can show each $R_{L}(f_Y)$ is nonzero 
at $q^{T_0}$ or at $q^{\frac1p T_0[U]}$:    
for if $a(T_0)\zeta_p^j\in L$ for some $j\in\Z$ and $a(T_0) \in {\Cal O}$, 
then $a(T_0) \in L \cap {\Cal O}=\ideal$.  
We now have a contradiction using the fact that power series over the domain 
$ {\Cal O}_{K(\zeta_p)}/L$ cannot be zero divisors.  \qed
\enddemo

We conclude this section by presenting vanishing and congruence conditions 
for Jacobi forms.  We thank O\. Richter for suggesting this.  
\proclaim{ \Equ{D21} Corollary}
Let $p \in \N$ be prime or one.  
Let $\phi \in J_{k,p}^{\text{cusp}}$ have the Fourier expansion 
$$
\phi(\tau,z)= \sum_{n \in \N, r \in \Z:\, 4np>r^2}
c(n,r)q^n \zeta^r= 
 \sum_{D \in \N, r \in \Z:\, -D \equiv r^2 \mod 4p} 
c(D) q^{\frac{D+r^2}{4p}} \zeta^r.  
$$
The form $\phi$ is trivial if and only if 
$c(D)=0$ whenever $D \le \frac{8}{225}\left( k \frac{1+p^2}{1+p} \right)^2$.  

Let $K$ be a number field, 
${\Cal O}$ its integers and 
$\ideal$ a prime ideal in ${\Cal O}$. 
Let $\phi \in J_{k,p}^{\text{cusp}}$ have all $c(D) \in {\Cal O}$.  
We have all $c(D) \in \ideal$ if and only if $c(D) \in \ideal$ 
whenever $D \le \frac{1}{27}\left( k \frac{1+p^2}{1+p} \right)^2$.
\endproclaim
\demo{Proof}
By Theorem~{2.2} of \cite{\refer{EZ}, page 23}, for index $p$ prime or $1$, 
the Fourier coefficients $c(n,r)$ depend only upon $4np-r^2$ so that we may write 
$c(n,r)=c(4np-r^2)$.  First examine the vanishing condition.  
The Fourier coefficients of $\Grit(\phi) \in S_2^k(K(p))$ are  
$$
a\pmatrix mp & r/2 \\ r/2 & n \endpmatrix = 
\sum_{\delta \vert (n,r,m)} \delta^{k-1}c\left( \dfrac{mn}{\delta^2}, \dfrac{r}{\delta} \right)
= \sum_{\delta } \delta^{k-1} c\left(\frac{4nm-r^2}{\delta}\right),  
$$ 
compare Theorem~\refer{c2}.  
From $c(D)=0$ whenever $D \le \frac{8}{225}\left( k \frac{1+p^2}{1+p} \right)^2$, 
we see that $a(T)=0$ whenever 
$4 \det(T) \le \frac{8}{225}\left( k \frac{1+p^2}{1+p} \right)^2$.  
Equivalently, $a(T)=0$ for $\delta(T) \le  \frac{\sqrt{2}}{15} k \frac{1+p^2}{1+p} $, 
which proves the vanishing of $\Grit(\phi)$ by Corollary~\refer{D9}.  
Hence $\phi$ also vanishes.  

Now examine the congruence condition for 
$\Grit(\phi) \in S_2^k(K(p))^{\epsilon}({\Cal O})$ for $\epsilon=(-1)^k$.  
If $a\left( T; \Grit(\phi) \right) \in \ideal$ for 
$T=\pmatrix mp & r/2 \\ r/2 & n \endpmatrix \in {}^p \Xtwo$ satisfying 
$w(T) \le \frac{k}{6}\frac{1+p^2}{1+p} $, 
then by Theorem~\refer{D10} we have $a\left( T; \Grit(\phi) \right) \in \ideal$ for 
all $T$.  It would follow that 
$c(D)=c(n,r;\phi)=a\left( \pmatrix p & r/2 \\ r/2 & n \endpmatrix;\Grit(\phi) \right) \in \ideal$ 
for all $D$.   
So suppose $T$ satisfies $w(T) \le \frac{k}{6}\frac{1+p^2}{1+p} $;  
we will show that $a\left( T; \Grit(\phi) \right) \in \ideal$.  
From $w(T) \ge \frac{2}{\mu_2} \delta(T)$, we see 
$\delta(T) \le \frac{1}{\sqrt{3}}w(T) \le \frac{k}{6\sqrt{3}} \frac{1+p^2}{1+p} $; 
thus $D=4mnp-r^2= 4\det(T) \le \frac{1}{27} \left( k \frac{1+p^2}{1+p} \right)^2$ 
and $c(D) \in \ideal$.  Therefore, by the above formula for the Fourier coefficients of $\Grit(\phi)$, 
we have $a\left( T; \Grit(\phi) \right) \in \ideal$.  \qed
\enddemo

{\bf{\newsection{sec5}{   Integral Closure.}} }

The constructions considered so far generate a large subring 
$R \subseteq M_2(K(p))$.  Let $R$ be the Hecke stable subring generated 
by Gritsenko lifts and  traces of theta series.  For weights $k \ge 3$, 
the dimension formulae of Ibukiyama reveal when $R$ contains 
$S_2^k(K(p))$; indeed, we usually have containment in our examples for $k \ge 4$.  
We  may use the following Lemmas to construct nonlifts in 
$S_2^2(K(p))$ by studying the integral closure of $S_2(K(p))$. 
This technique, in the case of elliptic modular forms, was used by 
J\. Tate \cite{\refer{Tate}} to construct nonbanal examples 
of weight one cusp forms for which 
the Artin Conjecture could be tested.  
The present article, in much the same spirit, 
aims at nonbanal examples of weight two paramodular cusp forms 
for which the Paramodular Conjecture can be tested.

Recall that $\Xn=\{ T\in \PQ: \forall\, v\in\Z^n,\,v'Tv\in\Z\}$ and 
$ {}^N\X2=\{\smtwomat{a}bbc\in \X2: N|a \}$.    

\proclaim{ \Equ{E5} Definition} 
Set  
$
{\Cal H}_N(2)=\{ H\in  S_2^4(K(N)): \supp(H) \subseteq {}^N\X2+{}^N\X2 \}  
$.  
Also define 
$
{\Cal H}_N(2)^{\pm}=\{ H\in  S_2^4(K(N))^{\pm}: \supp(H) \subseteq {}^N\X2+{}^N\X2 \}  
$.
\endproclaim
%
The next lemma is useful when $S_2^2(K(N))$ has linearly independent Gritsenko lifts.

\proclaim{ \Equ{E1} Lemma}
Let $g_1,g_2\in S_2^2(K(N))$ be nontrivial.  
Define a linear map 
$$
\align
\imath_{g_1,g_2}: S_2^2(K(N)) &\to 
\{(H_1,H_2)\in {\Cal H}_N(2)\times {\Cal H}_N(2): H_1\, g_2=H_2\, g_1 \}  \\
f &\mapsto (fg_1,fg_2).  
\endalign
$$
The map $\imath_{g_1,g_2}$ is injective.  
\endproclaim
\demo{Proof}
It suffices to point out that 
the image of $\imath_{g_1,g_2}$ is contained in ${\Cal H}_N(2)\times {\Cal H}_N(2)$.  \qed
\enddemo

\proclaim{ \Equ{E2} Corollary}
Let $g_1,g_2\in S_2^2(K(N))$ be nontrivial.  
For primes  $\ell$,   
we have the inequality  
$
\dim_{\C} S_2^2(K(N)) \le  
\dim_{\F_{\ell}}
\{(H_1,H_2)\in {\Cal H}_N(2)(\F_{\ell})\times 
{\Cal H}_N(2)(\F_{\ell}): H_1\, R_{\ell}(g_2)=H_2\, R_{\ell}(g_1) \}.
$
\endproclaim
This Corollary finds an upper bound on $\dim S_2^2(K(N))$ 
from a basis of $S_2^4(K(N))(\Z)$ modulo $ \ell$.   
When this upper bound equals $\dim J_{2,N}^{\text{cusp}}$  
then $S_2^2(K(N))$ is spanned by lifts.  
It is also useful to have versions of Lemma~\refer{E1} that treat the 
plus and minus spaces separately.  
We define a set, ${\Cal F}_N$, of indices that 
are ${\hat \Gamma}_0(N)$-equivalent to those that appear in the first Fourier-Jacobi coefficient  
of forms from $S_2(K(N))$. 
\proclaim{ \Equ{E5B} Definition}  
Set ${\Cal F}_N=\{ T \in {}^N\X2: \exists \smtwomat{N}bbc \in {}^N\X2, 
\exists U \in {\hat \Gamma}_0(N): T[U]= \smtwomat{N}bbc \}$ for $N \in \N$.  
Set ${}^N\X2' = {}^N\X2 \setminus {\Cal F}_N $.  
Define 
$ {\Cal H}'_N(2) =\{ H\in  {\Cal H}_N(2): \supp(H) \subseteq {}^N\X2+{}^N\X2' \}$.  
Define 
$ {\Cal H}''_N(2) =\{ H\in  {\Cal H}_N(2): \supp(H) \subseteq {}^N\X2'+{}^N\X2' \}$.  
Set $ {\Cal H}'_N(2)^{\pm}= {\Cal H}'_N(2) \cap S_2^4(K(N))^{\pm}$ and 
$ {\Cal H}''_N(2)^{\pm}= {\Cal H}''_N(2) \cap S_2^4(K(N))^{\pm}$.  
\endproclaim
\proclaim{ \Equ{E2B} Corollary}  
Let $N \in \N$.  
Let $g_1,g_2\in S_2^2(K(N))^{+}$ be nontrivial.  
We have the inequalities  
$
\dim_{\C} S_2^2(K(N))^{-} \le  
\dim 
\{(H_1,H_2)\in {\Cal H}_N(2)^{-} \times 
{\Cal H}_N(2)^{-} : H_1\, g_2=H_2\, g_1 \}  
$  and 
$
\dim_{\C}\left( S_2^2(K(N))^{+}/ \Grit(J_{2,N}^{\text{cusp}}) \right) \le  
\dim 
\{(H_1,H_2)\in {\Cal H}_N'(2)^{+} \times 
{\Cal H}_N'(2)^{+} : H_1\, g_2=H_2\, g_1 \}  
$.  
\endproclaim
\demo{Proof}
When $g_1$ and $g_2$ are plus forms, 
the restriction of $\imath_{g_1,g_2}$ to $S_2^2(K(N))^{-}$ has an image in 
${\Cal H}_N(2)^{-} \times {\Cal H}_N(2)^{-} $.  The injectivity of $\imath_{g_1,g_2}$ 
then proves the first inequality.  
For the second inequality, we note that each $ f \in S_2^2(K(N))^{+}$ has a 
unique represenative ${\hat f} \in f + \Grit(J_{2,N}^{\text{cusp}})$ 
whose first Fourier Jacobi coefficient vanishes.  
Thus, we have $\supp({\hat f}) \subseteq {}^N\X2'$.  
The map ${\hat\imath}_{g_1,g_2}$ defined on $S_2^2(K(N))^{+}$ by 
${\hat\imath}_{g_1,g_2}(f)=(g_1{\hat f}, g_2{\hat f})$ has kernel 
$\Grit(J_{2,N}^{\text{cusp}})$ and an image in 
${\Cal H}_N'(2)^{+} \times {\Cal H}_N'(2)^{+}$.  
\qed
\enddemo

This next Lemma is useful when $S_2^2(K(N))$ has 
a nontrivial Gritsenko lift.  

\proclaim{ \Equ{E3} Lemma}
Let $g\in S_2^2(K(N))$ be nontrivial.  
Define a (nonlinear) map 
$$
\align
\jmath_{g}: S_2^2(K(N)) &\to 
\{(F,H)\in {\Cal H}_N(2)\times {\Cal H}_N(2): H^2=Fg^2 \}  \\
f &\mapsto (f^2,fg ).  
\endalign
$$
The map $\jmath_{g}$ is bijective.  
\endproclaim
\demo{Proof}
The map $\jmath_{g}$ is clearly injective.  
On the other hand, suppose we have $F$, $H\in {\Cal H}_N(2)\subseteq S_2^4(K(N))$ 
with $H^2=Fg^2$.  The function $f=H/g $ is a weight $2$ 
meromorphic form whose square, $f^2={H^2}/{g^2}=F$, is holomorphic.  
Hence $f$ is holomorphic and 
$\jmath_{g}(f)=(f^2,fg )=(F,H)$.  \qed
\enddemo

Our strategy to construct nonlifts 
is to use Lemma~\refer{E1} to find a meromorphic function $f$
with distinct representations $f= {H_1}/{g_1}= {H_2}/{g_2}$.  
One then gains a pretty good idea of whether or not 
$f$ is holomorphic by applying the 
formulae for the action of the 
Hecke operators on Fourier coefficients 
to the initial Fourier expansion of $f$. To prove $f$ is holomorphic
we use Lemma~\refer{E3}.  This requires demonstrating the vanishing of 
$H^2-Fg^2\in S_2^8(K(p))^{+}$.      
One way to verify these identities in weight $8$ is to 
span  $S_2^8(K(p))^{+}$ but this is not always computationally feasible. 
Another path to proving holomorphicity would be to study the divisors of 
$g_1$ and $g_2$; note however, that 
$K(p) \backslash K(p)(\Half_1\oplus\Half_1)$ will always be a common divisor 
of any weight two paramodular cusp forms.

For small levels, when the dimension of the Gritsenko lifts is less than $2$, 
nonlifts were eliminated by the Restriction Technique, compare~\cite{\refer{PoorYuenComp}}.  
In order to avoid a lengthy description of the Restriction Technique here, we 
provide the following Lemmas.  The first was used for levels $37$, $43$ and $53$.  
%
%
\proclaim{ \Equ{E7B} Lemma} 
If $ \dim J_{2,N}^{\text{cusp}} \ge 1 $ and 
$\dim S_2^4(K(N))= 1+ \dim J_{4,N}^{\text{cusp}} $ then 
${\Cal H}''_N(2)=\{0\}$.  Also, 
if $\dim {\Cal H}_N(2) \le 1 \le \dim J_{2,N}^{\text{cusp}} $ then 
${\Cal H}''_N(2)=\{0\}$.  
\endproclaim
\demo{Proof}
Let $ f = \Grit(\phi) \in S_2(K(N))^{+}$ for a nontrivial $ \phi \in J_{2,N}^{\text{cusp}}$.  
Since the first Fourier-Jacobi coefficient of $f^2$ is zero, we have 
$ \Grit(J_{4,N}^{\text{cusp}}) \subsetneq \C f^2+\Grit(J_{4,N}^{\text{cusp}}) \subseteq S_2^4(K(N))^{+}$.  
Applying the hypothesis $\dim S_2^4(K(N))= 1+ \dim J_{4,N}^{\text{cusp}} $, we have 
$S_2^4(K(N)) = S_2^4(K(N))^{+} = \C f^2+\Grit(J_{4,N}^{\text{cusp}})$.  
From this we can show ${\Cal H}''_N(2)=\{0\}$.  

Any $H \in {\Cal H}''_N(2) \subseteq S_2^4(K(N))$ can be expressed as 
$H= \alpha f^2 + \Grit(\psi)$ for some $\alpha \in \C$ and $\psi \in J_{4,N}^{\text{cusp}}$.  
The first Fourier-Jacobi coefficient of $H$ is $\psi$ and so $\psi=0$ and $H = \alpha f^2$.  
The Fourier-Jacobi expansion of $H$ begins 
$
H\smtwomat{\tau}{z}{z}{\omega}= 
\alpha\, \phi(w,z)^2 e\left( 2N\omega \right)+ \dots
$
and because $\phi$ is nontrivial, 
$\supp( \phi^2 e(\cdot)^{2N} )$ contains a definite index 
$T=\smtwomat{2N}{r/2}{r/2}{m}$ with $r \in \Z$ and $m \in \N$.  
However, no element of this form is in ${}^N\X2'+{}^N\X2'$, 
so that  $H \in {\Cal H}''_N(2) $ implies that $\alpha=0$ and $H=0$.  
This proves the first assertion. 
The final assertion follows from the same argument because any 
$H \in {\Cal H}''_N(2)$ can be written $H = \alpha f^2$ for $f$ as above.  
\qed
\enddemo

Items \therosteritem4 and \therosteritem5 of the next Lemma are the most commonly used 
on the website~\cite{\refer{URL}} to prove that, for most primes $p<600$,  
$S_2^2(K(p))$ consist entirely of lifts.  
\proclaim{ \Equ{E7} Lemma} Let $N \in \N$.  
\roster
\item If $S_2^4(K(N))=\Grit( J_{4,N}^{\text{cusp}} )$ then ${\Cal H}_N(2)=\{0\}$.  
\item  If ${\Cal H}_N(2)^{+}=\{0\}$ then $S_2^2(K(N))=\{0\}$. 
\item  If ${\Cal H}''_N(2)=\{0\}$ then $S_2^2(K(N))=\Grit( J_{2,N}^{\text{cusp}} )$.   
\item  If ${\Cal H}''_N(2)^{+}=\{0\}$ then $S_2^2(K(N))^{+}=\Grit( J_{2,N}^{\text{cusp}} )$.   
\item  If $ \dim {\Cal H}_N(2)^{-} < \dim J_{2,N}^{\text{cusp}} $ then $S_2^2(K(N))^{-}=\{0\}$.  
\item  If $\dim{\Cal H}''_N(2)\le 2$, then $\Grit  ( J_{2,N}^{\text{cusp}}  )$ has codimension 
at most $1$ in $S_2^2(K(N))$. 
\item  If $\dim{\Cal H}''_N(2)^{+}\le 2$, then $\Grit  ( J_{2,N}^{\text{cusp}}  )$ has codimension 
at most $1$ in $S_2^2(K(N))^{+}$. 
\endroster     
\endproclaim
\demo{Proof}
For \therosteritem1 it is enough to show that 
${\Cal H}_N(2) \cap \Grit  ( J_{4,N}^{\text{cusp}}  )=\{0\}$.   
For $ f \in {\Cal H}_N(2)$ the first Fourier-Jacobi coefficient is $0$, 
hence $f \in \Grit  ( J_{4,N}^{\text{cusp}}  )$ further implies that $f=\Grit(0)=0$.  
For \therosteritem2, if $S_2^2(K(N)) \ne \{0\}$ then there is a nontrivial $f$ in 
$S_2^2(K(N))^{+}$ or $S_2^2(K(N))^{-}$.  In either case, $f^2 \in {\Cal H}_N(2)^{+}$ 
is nontrivial as well.  
For \therosteritem3, suppose that $f\in S_2^2(K(N))$ is not a Gritsenko lift.  
Let $\phi\in J_{2,N}^{\text{cusp}}$ be the first Fourier-Jacobi coefficient
of $f$.  Then ${\hat f}= f - \Grit(\phi)$ has a trivial 
first Fourier-Jacobi coefficient but is itself nontrivial.  
From
$\supp({\hat f}) \subseteq {}^N\X2' $ we see that   
${\hat f}^2\in {\Cal H}''_N(2)$ and ${\Cal H}''_N(2)\ne \{0\}$.  
Item~\therosteritem4 follows from the same argument.  
For \therosteritem5: 
when $f \in S_2^2(K(N))^{-}$ is nontrivial then 
$f \Grit(J_{2,N}^{\text{cusp}}) \subseteq {\Cal H}_N(2)^{-}$ 
so that $\dim J_{2,N}^{\text{cusp}} \le \dim {\Cal H}_N(2)^{-}$.  
For \therosteritem6, if $f,g\in S_2^2(K(N))$ are linearly independent 
modulo Gritsenko lifts, then ${\hat f}^2$, ${\hat f}{\hat g}$, ${\hat g}^2 \in {\Cal H}''_N(2)$ 
are linearly independent, noting that any quadratic has linear factors over $\C$.   
Item~\therosteritem7 follows from the same argument. 
\qed
\enddemo

{\bf{\newsection{sec6}{   Examples of weight two.}} }

We explain how the Theorems stated in the Introduction were proven.  
We first construct initial Fourier expansions of cusp forms 
in $S_2^4\left(K(p)\right)$ by multiplying Gritsenko lifts 
from  $S_2^2\left(K(p)\right)$, by applying Hecke operators 
and by tracing theta series from  $S_2^4\left(\Gamma_0(p)\right)$.  
The 
dimension formula of Ibukiyama in Theorem~\refer{C1} tells us if 
and when we have spanned $S_2^4\left(K(p)\right)$.   In this manner 
we were able to span $S_2^4\left(K(p)\right)$ for $p < 600$,   
with the exception of $p=499$.   
A nontrivial minus form of weight four first appears for $p=83$.  

\smallskip
\centerline{Table 4. \ Dimensions of the $\mu$-plus and $\mu$-minus subspaces in $S_2^4\left(K(p)\right)$} \nopagebreak
\smallskip
\begintable
 p\| 83 | 89 | 97 | 101 | 103 | 107 | 109 | 113 | 127 | 131  
| 137 | 139 | 149 | 151 | 157 | 163  \cr
$\dim S_2^4\left(K(p)\right)^{+}$
 \| 18 | 23 | 32 | 27 | 32 | 27 | 38 | 33 | 44 | 38  
| 45 | 51 | 51 | 59 | 65 | 65 \cr
$\dim S_2^4\left(K(p)\right)^{-}$
 \| 1 | 0 | 0 | 1 | 1 | 2 | 0 | 1 | 2 | 3  
| 2 | 2 | 3 | 2 | 3 | 4 
\endtable\smallskip  
\smallskip
\begintable
 p\| 167 | 173 | 179 | 181 | 191 | 193 | 197 | 199 | 211 | 223  
| 227 | 229 | 233 | 239 | 241 \cr
$\dim S_2^4\left(K(p)\right)^{+}$
 \| 55 | 62 | 65 | 83 | 73 | 92 | 78 | 91 | 100 | 106 | 91  
| 121 | 106 | 105 | 133    \cr
$\dim S_2^4\left(K(p)\right)^{-}$
 \| 8 | 8 | 6 | 3 | 7 | 4 | 10 | 6 | 7 | 12 | 18  
| 7 | 13 | 15 | 7  
\endtable\smallskip
\smallskip
\begintable
 p\| 251 | 257 | 263 | 269 | 271 | 277 | 281 | 283 | 293 | 307  
| 311 | 313 | 317 | 331 \cr
$\dim S_2^4\left(K(p)\right)^{+}$
 \| 113 | 124 | 120 | 134 | 149 | 161 | 149 | 155 | 149 | 177  
| 163 | 200 | 174 | 211    \cr
$\dim S_2^4\left(K(p)\right)^{-}$
 \| 18 | 18 | 23 | 20 | 17 | 17 | 18 | 24 | 31 | 30  
| 32 | 21 | 34 | 26  
\endtable\smallskip
\smallskip
\begintable
 p\| 337 | 347 | 349 | 353 | 359 | 367 | 373 | 379 | 383 | 389  
| 397 | 401 | 409 | 419 \cr
$\dim S_2^4\left(K(p)\right)^{+}$
 \| 227 | 192 | 239 | 212 | 210 | 241 | 263 | 264 | 226 | 256  
| 289 | 274  | 318  | 272     \cr
$\dim S_2^4\left(K(p)\right)^{-}$
 \| 25 | 47 | 29 | 42 | 45 | 45 | 39 | 39 | 62 | 48  
| 49 |  48 |  39 |  69  
\endtable\smallskip
\smallskip
\begintable
 p\| 421 | 431 | 433 | 439 | 443 | 449 | 457 | 461 | 463 | 467  
| 479 | 487 | 491 | 499 \cr
$\dim S_2^4\left(K(p)\right)^{+}$
 \| 333 | 287 | 343 | 333 | 297 | 333 | 378 | 335 | 362 | 321  
| 341 | 393  | 363  | 422 or 423     \cr
$\dim S_2^4\left(K(p)\right)^{-}$
 \| 43 | 73 | 53 | 64 | 82 | 65 | 59 | 83 | 76 | 98  
| 99 |  88 |  98 |  81 or 80  
\endtable\smallskip
\smallskip
\begintable
 p\| 503 | 509 | 521 | 523 | 541 | 547 | 557 | 563 | 569 | 571  
| 577 | 587 | 593 | 599 \cr
$\dim S_2^4\left(K(p)\right)^{+}$
 \| 363 | 400 | 426 | 437 | 506 | 478 | 460 | 443 | 502 | 530  
| 558 | 480  | 518  | 519     \cr
$\dim S_2^4\left(K(p)\right)^{-}$
 \| 120 | 104 | 101 | 112 | 90 | 119 | 138 | 156 | 121 | 117  
| 114 |  169 |  156 |  156  
\endtable\smallskip

For example, $\dim S_2^4\left(K(229)\right) =128$ and products of
the $7$ weight two Gritsenko  lifts  give $28$ linearly independent
cusp  forms of weight~$4$.  Applying the Hecke operators
$T_2$, $T_3$, $T_2^2$ and $T_5$,  we span spaces of dimension  
$56$, $84$, $112$ and at least $121$, respectively. 
When the action of the Hecke operators seemed to stabilize,
we used Theorems~\refer{B1} and \refer{B2} to compute initial
Fourier expansions of
$\Tr\left( \vartheta_P \vartheta_Q\right)\in M_2^4\left(K(229)\right) $
for $P,Q\in {\Cal A}$ where
$$ 
{\Cal A}=\{ 
\smfourmatNUMBERTEN{12}{24}{24}{1}{-3}{0}{-1}{3}{10}, 
\smfourmatNUMBERTEN{12}{12}{40}{1}{2}{1}{-1}{4}{2}, 
\smfourmatNUMBERTWELVE{14}{16}{28}{2}{-1}{5}{-1}{6}{8}, 
\smfourmatNUMBERTWELVE{18}{18}{20}{3}{5}{8}{-1}{0}{3}
 \}. 
$$      
Any linear combination of these  that cancels the constant term 
in the Fourier expansion 
gives an element of  $S_2^4\left(K(229)\right) $ as explained after Lemma~\refer{B3b}.  
The addition of these linear combinations of theta traces 
increased the dimension of the constructed subspace to at least $128$ and 
hence spanned $S_2^4\left(K(229)\right) $.    
The involution $\mu$ twins the indices of the Fourier coefficients so that 
the dimensions of the 
plus and minus subspaces are easily computed to be $121$ and $7$.   
Notice that the theta traces were only necessary to fill the 
minus space.  
For the case $p=229$ just described 
the seven theta blocks $\TB_2(\Sigma_i)$ are:   
[ 2, 2, 3, 4, 5, 7, 7, 9, 10, 11 ],
    [ 2, 2, 3, 3, 5, 5, 7, 8, 10, 13 ], 
    [ 2, 2, 2, 3, 4, 5, 6, 8, 10, 14 ], 
    [ 1, 3, 4, 4, 5, 6, 7, 8, 11, 11 ], 
    [ 1, 3, 3, 4, 6, 6, 7, 9, 10, 11 ],
    [ 1, 3, 3, 4, 5, 7, 8, 8, 10, 11 ],
    [ 1, 3, 3, 4, 4, 5, 7, 8, 10, 13 ].  
See \cite{\refer{URL}} for full comments on these computations.  

\demo{Proof of Theorem~\refer{A2}}
This is an application of Lemmas~\refer{E1}, \refer{E3}, \refer{E7B}, \refer{E7} 
and Corollaries~\refer{E2}, and \refer{E2B}.    
Using initial Fourier expansions of a basis for  $S_2^4\left(K(p)\right)$ 
we compute initial Fourier expansions of 
linear combinations whose span contains 
$
{\Cal H}_p(2)   
$, 
$
{\Cal H}_p(2)^{\pm}    
$
and 
$
{\Cal H}_p''(2)^{\pm}  
$.  
If $\dim {\Cal H}_p''(2)^{+}=0$ then $S_2^2(K(p))^{+}$ has no nonlifts by 
Lemma~\refer{E7}, item~\therosteritem4.  
If $\dim {\Cal H}_p(2)^{-}< \dim J_{2,p}$ then $S_2^2(K(p))^{-}$ is trivial by 
Lemma~\refer{E7}, item~\therosteritem5. 
Otherwise, for various pairs of Gritsenko lifts 
$g_1,g_2 \in S_2^2\left(K(p)\right)$, we 
compute upper bounds for 
$\dim \{(H_1,H_2)\in {\Cal H}_p(2)\times {\Cal H}_p(2): H_1g_2=H_2g_1 \} $ 
by linear algebra.  
If any of these dimensions equal $\dim J_{2,p}$, then the 
injectivity of $\imath_{g_1,g_2}$ from  Lemma~\refer{E1} tells us that 
$S_2^2\left(K(p)\right)= \Grit\left( J_{2,p}\right)$.  
These computations may also be performed modulo a prime $\ell$, we 
usually take $\ell=\fav$, and we use Corollary~\refer{E2} to get the same 
conclusion: $S_2^2\left(K(p)\right)= \Grit\left( J_{2,p}\right)$.  
When this last technique was used, the 
two Gritsenko lifts $g_1$ and $g_2$ that worked were  
recorded at \cite{\refer{URL}}. 

Let us give an illustration of this last technique in the case of $p=229$. 
We have $\dim J_{2,229}=7$,
$\dim {\Cal H}_{229}(2) \le 31$ and for
$g_1=\TB_2( 2, 2, 3, 4, 5, 7, 7, 9, 10, 11 )$ and
$g_2=\TB_2(  2, 2, 3, 3, 5, 5, 7, 8, 10, 13  )$ we have
$\dim_{\F_{\fav}} \{(H_1,H_2)\in \left({\Cal H}_{229}(2)
\times {\Cal H}_{229}(2)\right)(\F_{\fav}): H_1R_{229}(g_2)=H_2R_{229}(g_1) \}\le 7$.  
Therefore we have 
$S_2^2\left(K(229)\right)= \Grit\left( J_{2,229}\right)$.   
We note that this process requires choosing $31$ linearly independent elements from 
$S_2^4(\Z)$ whose $\C$-span contains ${\Cal H}_{229}(2)$ and whose 
reductions modulo $\fav$ remain linearly independent over $\F_{\fav}$.  
 For other primes, 
we simply note which case of Lemma~\refer{E7} applies.  
We were only able to be certain that we had spanned 
 a codimension one subspace of 
$S_2^4\left(K(499)\right)$.  
Still, it can be shown that $S_2^2\left(K(499)\right)$ is spanned by 
Gritsenko lifts and this is done in Appendix~A.  
\qed
\enddemo


Although the primes $p$ for which the Gritsenko lifts span 
$S_2^2\left(K(p)\right)$ 
produce no new paramodular cusp forms--- they do give important
evidence for the Paramodular Conjecture~\refer{A1b} 
because these $p$ should also be, and so far are, primes for which there is 
no abelian surface defined over $\Q$ of conductor $p$.  
The more interesting primes $p$ are those where  
$\dim \Grit\left( J_{2,p}\right)< \dim S_2^2\left(K(p)\right)$ 
since, in order  
to test the Paramodular Conjecture, 
we must construct nonlift paramodular Hecke eigenforms.  
The first example is $p=277$.  


\proclaim{ \Equ{F1} Theorem }
We have $\dim S_2^2(K(277))=11$ whereas the dimension of 
Gritsenko lifts in $S_2^2(K(277))$ is 
$\dim J_{2,277}=10$.  
Let $G_i=\Grit\left( \TB_2(\Sigma_i) \right)$ for $1 \le i \le 10$ 
be the lifts of the $10$ theta blocks given by $\Sigma_i\in$ 
$\{${\rm [2, 4, 4, 4, 5, 6, 8, 9, 10, 14], 
    [2, 3, 4, 5, 5, 7, 7, 9, 10, 14], 
    [2, 3, 4, 4, 5, 7, 8, 9, 11, 13], 
    [2, 3, 3, 5, 6, 6, 8, 9, 11, 13], 
    [2, 3, 3, 5, 5, 8, 8, 8, 11, 13], 
    [2, 3, 3, 5, 5, 7, 8, 10, 10, 13], 
    [2, 3, 3, 4, 5, 6, 7, 9, 10, 15], 
    [2, 2, 4, 5, 6, 7, 7, 9, 11, 13], 
   [2, 2, 4, 4, 6, 7, 8, 10, 11, 12], 
    [ 2, 2, 3, 5, 6, 7, 9, 9, 11, 12] }$\}$.  
Let 
$$
\align
Q   =&
-14G_1^2-20G_8G_2+11G_9G_2+6G_2^2-30G_7G_{10}+15G_9 G_{10}+15 G_{10} G_1
-30 G_{10} G_2   \\
-&30 G_{10} G_3+5 G_4 G_5+6 G_4 G_6+17 G_4 G_7-3 G_4 G_8-5 G_4 G_9
-5 G_5 G_6 +20 G_5 G_7   \\
-&5 G_5 G_8-{10} G_5 G_9-3 G_6^2+13 G_6 G_7+3 G_6 G_8
-{10} G_6 G_9-22 G_7^2+G_7 G_8+15 G_7 G_9  \\
+&6 G_8^2-4 G_8 G_9-2 G_9^2+20 G_1 G_2
-28 G_3 G_2+23 G_4 G_2+7 G_6 G_2-31 G_7 G_2+15 G_5 G_2 \\
+&45 G_1 G_3-{10} G_1 G_5
-2 G_1 G_4-13 G_1 G_6-7 G_1 G_8+39 G_1 G_7-16 G_1 G_9-34 G_3^2  \\
+&8 G_3 G_4
+20 G_3 G_5+22 G_3 G_6+{10} G_3 G_8+21 G_3 G_9-56 G_3 G_7-3 G_4^2,  \\
L    =&      -G_4 + G_6 + 2 G_7 + G_8 - G_9 + 2 G_3 - 3 G_2 - G_1.  
\endalign
$$
Let $f ={Q}/{L}$ define a meromorphic paramodular form of weight $2$.  
The form $f$ is holomorphic and 
$S_2^2(K(277))(\Q)=\Span_{\Q}\left(f,G_1,\dots,G_{10}  \right)$.  
The form $f$ is a Hecke eigenform with spinor Euler factors
$$
\align
Q_2(f,x)  &=  1 + 2 x + 4 x^2 + 4 x^3 + 4 x^4,  
\\
Q_3(f,x)  &=  1 +   x +   x^2 + 3 x^3 + 9 x^4,  
\\
Q_5(f,x)  &=  1 +   x - 2 x^2 + 5 x^3 + 25 x^4  
.
\endalign
$$
Additionally, let
$$
\align
{\hat Q}= &
{-}55 G_5^2-13 G_8 G_2-148 G_7 G_8+17 G_4 G_9+123 G_4 G_8+160 G_4 G_7
{-}12 G_4 G_6+73 G_4 G_5  \\
-&163 G_{10} G_3-148 G_{10} G_2+211 G_{10} G_1
-62 G_9 G_{10}-94 G_7 G_{10}+145 G_9 G_2-18 G_7 G_9   \\
{+}&58 G_6 G_9
{-}138 G_6 G_8-{10} G_6 G_7+52 G_5 G_9-{10}9 G_5 G_8-73 G_5 G_6-4 G_6 G_2
+154 G_4 G_2  \\
{-}&156 G_3 G_2{+}17 G_1 G_2-5 G_8 G_9+333 G_3 G_5
+37 G_3 G_4+26 G_1 G_9+235 G_1 G_7-71 G_1 G_8  \\
{+}&49 G_1 G_6{-}34 G_1 G_4
-151 G_1 G_5+396 G_1 G_3-2 G_5 G_2-19 G_7 G_2-245 G_3 G_7
-58 G_3 G_9 \\
+&83 G_3 G_8+113 G_3 G_6-54 G_7^2-63 G_8^2+12 G_9^2
-24 G_4^2-404 G_3^2+89 G_5 G_7-196 G_1^2  \\
+&8 G_2^2-24 G_6^2
+5 G_{10} G_5+63 G_{10} G_6+9 G_8 G_{10}-3 G_4 G_{10}    ,\\
{\hat L}=  &G_{10}+15 G_1+6 G_2-5 G_3+5 G_4-6 G_5-5 G_6-13 G_7+7 G_8-8 G_9.  
\endalign
$$
We have the following identity in $S_2^8(K(277))$: 
$$
Q^2 + {\hat L} Q L + {\hat Q} L^2=0.  \tag \Equ{F2} 
$$
\endproclaim
\demo{Proof}
For the pair of Gritsenko lifts 
$$
\align
v_1 &= L= -G_4+G_6+2 G_7+G_8-G_9+2 G_3-3 G_2-G_1,  
\\
v_2 &=     G_1 - 3 G_2 + G_3 - 2 G_4 + 2 G_6 + G_7 + 2 G_8 - 2 G_9, 
\endalign
$$
we computed 
$\dim \{(H_1,H_2)\in {\Cal H}_{277}(2)\times {\Cal H}_{277}(2): H_1v_2=H_2v_1 \}\le 11 $.  
Therefore we have 
$\dim  S_2^2(K(277)) \le 11$ by Lemma~\refer{E1}.  
There is a pair $(Q,H_2)\in {\Cal H}_{277}(2)\times {\Cal H}_{277}(2)$ 
that is not an image of $\imath_{v_1,v_2}\left(\Grit\left( J_{2,277} \right) \right)$, 
that appears to satisfy $Qv_2=H_2v_1$ and that produces an initial Fourier expansion for 
$f= {Q}/{v_1}$ that is consistent with being a Hecke eigenform.  
To prove that $f$ is holomorphic, it suffices to prove the identity~\refer{F2}.  
For if we have $Q^2 + {\hat L} Q L + {\hat Q} L^2=0$, then  we have 
$$
f^2+ {\hat L}f+{\hat Q}=0,
 \text{ for  } 
f=\dfrac{Q}{L},
$$
so that $f$ is in the integral closure of 
$ M_2(K(277))$ and hence is holomorphic.  
An application of the Siegel $\phi$ map 
to $f^2+ {\hat L}f+{\hat Q}=0 $ shows that 
$(\phi f)^2=0$ so that $f$ is a cusp form.  
From $Q,L\in S_2(K(277))(\Z)$ we see that 
$f\in S_2^2(K(277))(\Q)$.  The Fourier expansion of $f$ 
is computed by the long division of $L$ into $Q$.   
The action of the Hecke operators on $ S_2^2(K(277))(\Q)$ and 
the Euler factors of $f$ are 
computed from these Fourier coefficients.  

The identity~\refer{F2} was proven by spanning $S_2^8(K(277))$.  
The space $S_2^8(K(277))$ was spanned in the same manner
as $S_2^4(K(277))$ but at greater expense. 
Products of the $56$ weight four Gritsenko lifts gave at least 1496 linearly independent
elements, at which point the Hecke operators  
$T_2$ and $T_3$   were applied, resulting in at least 1760 linearly independent elements in 
$S_2^8(K(277))^{+}$.  
  Theta traces
$\Tr\left( \vartheta_P \vartheta_Q\right)\in M_2^4\left(K(277)\right) $
for $P,Q\in {\Cal A}$ were computed 
and multiplied by $S_2^4(K(277))^{+}$
where
$$
\align
{\Cal A}=\{    
&\smfourmatNUMBERSIXTEEN{16}{18}{26}{5}{5}{7}{-5}{1}{3}, 
\smfourmatNUMBERFOURTEEN{20}{20}{22}{1}{-7}{3}{-6}{-6}{5},  
\smfourmatNUMBERTWELVE{14}{20}{30}{3}{1}{6}{-5}{-2}{1}, \\
&\smfourmatNUMBERFOURTEEN{14}{18}{30}{4}{1}{6}{4}{5}{2}, 
\smfourmatNUMBERFOURTEEN{14}{16}{32}{1}{-5}{0}{3}{-1}{6}, 
\smfourmatNUMBERSIXTEEN{16}{18}{28}{5}{5}{7}{-3}{-8}{0}  
 \}. 
\endalign
$$
Symmetrization of these with respect to $\mu$ gave an additional 57 cusp forms, 
so that $\dim S_2^8(K(277))^{+} \ge 1817$.  
Products from $S_2^4(K(277))^{+} \, S_2^4(K(277))^{-} $ gave at least 595 
linearly independent weight~8 minus forms. 
Finally, the Hecke operator $T_2$ was applied to the minus forms
computed in this manner, for an estimate  $\dim S_2^8(K(277))^{-} \ge 712$. 
These steps gave a subspace that spanned at least $2529$ dimensions over $\F_{\fav}$. 
This is the correct dimension by Ibukiyama's formula.  
Therefore, the dimension of the plus space in $S_2^8(K(277))$ is $1817$ and
the dimension of the minus space  is $712$.  
This gave a determining set of $2529$ Fourier coefficients for
$S_2^8(K(277))$.  The identity~\refer{F2} was then checked to vanish
on this determining set of $2529$ Fourier coefficients. 
\qed
\enddemo 

In connection with Theorem~\refer{F1} we mention that the hyperelliptic 
curve $C$ of genus~$2$ defined by 
$y^2+y=x^5+5x^4+8x^3+6x^2+2x$ 
has a Jacobi variety $\Jac(C)$ 
defined over $\Q$ with conductor $277$.  
We refer to the companion article \cite{\refer{BK}}  for these arithmetic results.  
The Euler factors of the Hasse-Weil $L$-function of $\Jac(C)$ 
are identical to those of the spinor $L$-function of $f$ 
for the primes $q=2,3$ and $5$.  We know further equalities of
eigenvalues but not nearly enough to prove that these $L$-functions are equal.  
The abelian surface $\Jac(C)$  has rational $15$-torsion.  
This is of a piece with the congruence in the following Theorem. 
We thank A\. Brumer for suggesting this proof.  

\proclaim{ \Equ{F3} Theorem }
Let $f$ be as in Theorem~\refer{F1}.  
We have $f\in S_2^2\left( K(277) \right)(\Z)$.  
The first Fourier Jacobi coefficient of $f$ is 
$$
\align
\phi=&  -5  \TB_2(\Sigma_5)+3  \TB_2(\Sigma_4)-3  \TB_2(\Sigma_6)+4  
\TB_2(\Sigma_7)+6  \TB_2(\Sigma_8)  \\
&+2  \TB_2(\Sigma_9)-2  
\TB_2(\Sigma_3)-2  \TB_2(\Sigma_2)- \TB_2(\Sigma_1) \in J_{2,277}. 
\endalign
$$
Let 
$
R= \Grit(\phi)= -5  G_5+3  G_4-3  G_6+4  G_7+6  G_8+2  G_9-2  G_3-2  G_2- G_1 \in S_2^2\left( K(277)
\right)(\Z).  
$
We have 
$$
\forall \, T \in {}^{277}\X2, 
a(T;f) \equiv a(T;R) \mod 15.  
$$
\endproclaim
\demo{Proof}
By a result of Shimura \cite{\refer{Shimura}}, 
we have $ Nf\in S_2^2\left( K(277) \right)(\Z)$
for some positive integer $N$.    
We choose $N$ to be minimal with this property and   
show that $N=1$.  
If not, suppose a prime $\ell $ divides $N$.    
From $f= {Q}/{L}$, we have 
$L(Nf)=NQ$.  
In $S_2^2\left( K(277) \right)(\F_{\ell})$ we have 
$R_{\ell}(L) R_{\ell}(Nf) =0$ because $Q$ has integral Fourier coefficients.    
The Fourier expansion of $L$ has unit content:
$$
L(\Omega)= 
\,e\left( \<\Omega,\tfrac12 \smtwomat{98\cdot 277}{233}{233}{2} \> \right)
+ 
\,e\left( \<\Omega,\tfrac12 \smtwomat{26\cdot 277}{120}{120}{2} \> \right)
+2 
\,e\left( \<\Omega,\tfrac12 \smtwomat{326\cdot 277}{601}{601}{4} \> \right)
+ \dots 
$$ 
where the ``$+\dots$'' indicates ${\hat \Gamma}_0(277)$-equivalent terms and 
terms of higher dyadic trace.  
Therefore $R_{\ell}(L)$ is nontrivial and cannot be a zero divisor. 
This shows that $R_{\ell}(Nf) =0$ and that 
$\frac{N}{\ell}f$ is integral, 
contradicting the minimality of $N$.  
The congruence now follows formally from the quotient 
$f= {Q}/{L}$. To see this note that
$$
\dfrac{Q}{L}= R + 15\, \dfrac{J}{L}
$$
holds identically in the $G_1,\dots,G_{10}$, considered as ten variables, if we set 
$$
\align
&J=
-G_7 G_8+G_4 G_7-2 G_{10} G_3-2 G_{10} G_2+G_{10} G_1+G_9 G_{10}-2 G_7 G_{10}
+G_9 G_2  \\
&+G_7 G_9-G_6 G_9+G_6 G_7-G_5 G_9+2 G_4 G_2-2 G_3 G_2+G_1 G_2
+2 G_3 G_5-G_1 G_9  \\
&+3 G_1 G_7-G_1 G_6-G_1 G_5+3 G_1 G_3-G_7 G_2
-4 G_3 G_7+G_3 G_9+2 G_3 G_6-2 G_7^2  \\
&-2 G_3^2+2 G_5 G_7-G_1^2.  
\endalign
$$
Therefore we have 
$R_{\ell}(L) R_{\ell}(f-R)=0 $ for $\ell=3$ and~$5$.  
The conclusion follows since $R_{\ell}(L)$ is not a zero divisor. 
\qed
\enddemo

%
%
In $\Grit\left( J_{2,p} \right) \cong S_1^2\left( \Gamma_0(p)^{*} \right)$, 
the $10$~eigenforms break into a nine dimensional piece and a rational 
eigenform: 
$2G_1+G_2-2G_3+3G_4-5G_5-3G_6+4G_7+6G_8+2G_9$.  
This rational eigenform is visibly congruent to $R$ modulo~$3$.  
Suitable multiples of the other eigenforms are each congruent to $R$ modulo a prime above~$5$.  
We have given the computations for the case $p=277$ in some detail.  
We now tabulate the results for the other exceptional primes given in Theorem~\refer{A2}.

Table~5  lists data for every  prime $p< 600$ that possibly could have a 
Hecke eigenform in $ S_2^2\left( K(p) \right) $ not in the image of the 
Gritsenko lift from $J_{2,p}$.  The existence of a  nonlift  
has been proven in some cases but remains conjectural in others.  
The $\epsilon$ and the ${\Cal O\/} $ columns indicate that there is a nonlift 
$f \in S_2^2\left( K(p) \right)^{\epsilon}( {\Cal O\/} )$.   It is  
defined by $f=Q/L$ 
for some $Q\in S_2^4\left( K(p) \right)$ and $L\in S_2^2\left( K(p) \right)$.  
There is an identity in $ S_2^8\left( K(p) \right)^{+}$, 
$Q^2 +L Q {\hat L} + L^2 {\hat Q}=0$ 
for some ${\hat Q}\in S_2^4\left( K(p) \right)$ and ${\hat L} \in S_2^2\left( K(p) \right)$, 
that would certify the holomorphicity of $f$.  
We give $Q$, ${\hat Q} $, $L$ and ${\hat L} $ at \cite{\refer{URL}} for each case.  
The verification of the weight~$8$ identity is expensive 
and a checkmark in the last column indicates that this verification is complete.  
An unchecked entry leaves open the possibility that $f$ is properly meromorphic.  
The dimension of $S_2^2\left( K(p) \right)^{\epsilon}$ would need to be reduced 
by one in each such happenstance.  
Table~5 gives some Hecke eigenvalues of these nonlifts 
and Appendix~B gives some Fourier coefficients; 
many more  Fourier coefficients may be found at \cite{\refer{URL}}. 
Euler factors of the weight two form $f$ for $q$ prime to the level $p$ may be computed as 
$$
Q_q(f,x)=
1-\lambda_{q}x+(\lambda_{q}^2-\lambda_{q^2}-1)x^2
-q\,\lambda_{q} x^3+q^{2}x^4.  
$$
For the identifications of these Euler factors 
with those of known rational abelian varieties see the 
companion article \cite{\refer{BK}}.  
Additionally, each entry displays primes $\ell$ for which we 
conjecture the congruence $f \equiv \Grit(\phi) \mod \ell$, 
where $\phi$ is the first Fourier Jacobi coefficient of $f$;    
the congruences have been proven contingent upon the existence of the nonlift $f$.  
We cannot compute enough Fourier coefficients to verify these congruences with Theorem~\refer{D10} 
but congruences can be detected by the form of $Q$ and $L$, 
as in the proof of Theorem~\refer{F3}.  
By examining minors of rational bases, 
we have also been able to prove is that these are the only possible such congruences, see \cite{\refer{URL}}.  

\smallskip
\centerline{{\bf Table 5.\/}   Hecke eigenforms in $S_2^2\left( K(p) \right)^{\epsilon}({\Cal O\/})  $
but not in $\Grit\left( J_{2,p} \right)$. }   \nopagebreak
\smallskip
\begintable
 p\| $\dim J_{2,p}$ | $\dim S_2^2\left(K(p)\right)$ | ${\Cal O\/}$ | $\epsilon$ | $\ell$ | $\lambda_2$ | $\lambda_3$ |
$\lambda_4$ | $\lambda_5$ | $\lambda_7$ | $\lambda_9$ | $\lambda_{11}$ |  proven   \cr
277
 \| 10 | 11 | $\Z$ | $+$ | $\{3,5\}$ | -2 | -1 | -1 | -1 | 1  
| -1 | -2 |  \checkmark     \cr
349
 \| 11 | 12 | $\Z$ | $+$ | $\{13\}$ | -2 | -1 | 1 | -1 | -2  
| 1 | 1 |      \cr
353
 \| 11 | 12 | $\Z$ | $+$ | $\{11\}$ | -1 | -2 | -3 | 1 | 0  
| -1 | 2 |      \cr
389
 \| 11 | 12 | $\Z$ | $+$ | $\{2,5\}$ | -1 | -2 | -2 | 1 | -3  
| 1 | 4 |      \cr
461
 \| 12 | 13 | $\Z$ | $+$ | $\{7\}$ | 0 | -3 | -3 | 1 | 0  
| 2 |  2  |     \cr
523
 \| 17 | 18 | $\Z$ | $+$ | $\{2,5\}$ | -1 | 0 | -2 | -4 | 2  
| -1 | 2  |      \cr
587
 \| 18 | 20 | $\Z$ | $+$ | $\{11\}$ | -1 | 0 | -3 | 0 | -2  
| -2 |   |      \cr
587
 \| 18 | 20 | $\Z$ | $-$ | $\{ \}$ | -3 | -4 | 3 | -2 | 0  
| 6 |  -1  |  
\endtable\smallskip

%

Finally, we remark that for $p=277$, there is a $4$ dimensional space of weight~$2$ 
Gritsenko lifts whose product with the nonlift $f$ lands in the span of the
products of the weight~$2$ Gritsenko lifts.  
For $p=353$, the analogous space is $3$ dimensional.   These curious subspaces 
of weight~$2$ Gritsenko lifts are instrinsic, as are the corresponding 
subspaces  in $J_{2,p}$ and $S_1^{3/2}\left(\Gamma_0(4p)\right)^{+}$.  
There is no corresponding space in $S_1^2\left( \Gamma_0(p)  \right)$ 
because the Shimura correspondence is noncanonical.  We do, however, have 
an intrinsic subspace of 
$S_1^2\left( \Gamma_0(p)  \right) \otimes S_1^2\left( \Gamma_0(p)  \right)$ 
constructed by composing the Witt map with the Saito-Kurokawa lift.  
It would be interesting to give an alternative characterization of these 
spaces directly in terms of elliptic modular forms.  
%
We state the following Theorem as representative of the unproven cases in Table~5.  

\proclaim{ \Equ{F4} Theorem }
We have $\dim S_2^2(K(353)) \le 12$ whereas the dimension of 
Gritsenko lifts in $S_2^2(K(353))$ is 
$\dim J_{2,353}=11$.  
Let $G_i=\Grit\left( \TB_2(\Sigma_i) \right)$ for $1 \le i \le 11$ 
be the lifts of the $11$ theta blocks given by $\Sigma_i\in$ 
$\{[ 3, 4, 4, 4, 6, 7, 8, 10, 12, 16]$, 
    $[3, 3, 4, 4, 5, 7, 7, 10, 12, 17 ]$, 
    $[ 2, 3, 5, 5, 6, 7, 9, 10, 11, 16 ]$, 
    $[ 2, 3, 5, 5, 6, 7, 8, 10, 13, 15 ]$, 
    $[ 2, 3, 5, 5, 5, 7, 10, 10, 12, 15 ]$, \newline 
    $[ 2, 3, 4, 6, 6, 7, 9, 9, 13, 15 ]$, 
    $[ 2, 3, 4, 6, 6, 7, 8, 10, 14, 14]$, 
    $[ 2, 3, 4, 5, 6, 7, 9, 11, 13, 14 ]$, \newline
    $[ 2, 3, 4, 5, 5, 7, 8, 9, 12, 17 ]$, 
    $[ 2, 3, 4, 4, 6, 8, 10, 11, 12, 14 ]$, 
    $[ 2, 3, 3, 5, 6, 9, 9, 11, 12, 14 ]\}$.  
%
Let 
$$
\align
Q   =&
-G_{10}G_{11} - 2 G_{11} G_2 + G_{10} G_3 + 4 G_{11} G_3 + 2 G_2 G_3 - 4 G_3^2 - 
5 G_{11} G_4 + 5 G_3 G_4 \\
+ &G_{11} G_5 - G_3 G_5 + 6 G_{11} G_6 - 6 G_3 G_6 - 
11 G_1 G_7 + 11 G_{10} G_7 - 9 G_{11} G_7 + 9 G_3 G_7 \\ 
+ &11 G_4 G_7 - 
11 G_6 G_7 + 11 G_7^2 + 11 G_1 G_8 - 11 G_{10} G_8 + 9 G_{11} G_8 - 
9 G_3 G_8 - 11 G_4 G_8 \\ 
+ &11 G_6 G_8 - 22 G_7 G_8 + 11 G_8^2 - 
2 G_{11} G_9 + 2 G_3 G_9, \\
L    =&      -G_{11} + G_3.  
\endalign
$$
Let $f ={Q}/{L}$ define a meromorphic form of weight $2$.  
Let {\rm 
$\{$ 
[1, 8, 12, 12], [2, 3, 4, 18], [2, 3, 12, 14], 
[2, 6, 12, 13], [3, 10, 10, 12], [4, 7, 12, 12], [5, 6, 6, 16],  
[8, 8, 9, 12]$\}$ \/} be labeled $\{\Xi_1,\dots,\Xi_8\}$.  
For  $\Xi_i=[a,b,c,d]$ and $\Xi_j=[\alpha,\beta,\gamma,\delta]$, 
define $W(i,j)=\Grit(\vartheta_{a}\vartheta_{b}\vartheta_{c}\vartheta_{d})
\Grit(\vartheta_{\alpha}\vartheta_{\beta}\vartheta_{\gamma}\vartheta_{\delta})
-\Grit\left( \TB_4(a,b,c,d,\alpha,\beta,\gamma,\delta) \right)$. Set 
$$
\align
{\hat Q}= &
22 G_1^2 - 165 G_1 G_{10} + 133 G_{10}^2 + 385 G_1 G_{11} - 646 G_{10} G_{11} + 
407 G_{11}^2 + 224 G_{10} G_2 \\ 
- &214 G_{11} G_2 + 125 G_2^2 - 286 G_1 G_3 + 
330 G_{10} G_3 - 562 G_{11} G_3 + 121 G_2 G_3 + 127 G_3^2 \\
- &165 G_1 G_4 + 
230 G_{10} G_4 - 491 G_{11} G_4 + 119 G_2 G_4 + 330 G_3 G_4 + 113 G_4^2 - 
220 G_1 G_5 \\
+ &231 G_{10} G_5 - 25 G_{11} G_5 - 22 G_2 G_5 + 113 G_3 G_5 + 
110 G_4 G_5 - 100 G_5^2 + 121 G_1 G_6 \\
- &177 G_{10} G_6 + 103 G_{11} G_6 - 
123 G_2 G_6 - 165 G_3 G_6 - 60 G_4 G_6 + 110 G_5 G_6 - 30 G_6^2 \\
- &572 G_1 G_7 + 594 G_{10} G_7 - 479 G_{11} G_7 + 99 G_2 G_7 + 457 G_3 G_7 + 
330 G_4 G_7 - 103 G_5 G_7 \\
- &66 G_6 G_7 + 205 G_7^2 + 473 G_1 G_8 - 
495 G_{10} G_8 + 644 G_{11} G_8 - 220 G_2 G_8 - 501 G_3 G_8 \\
- &352 G_4 G_8 - 
18 G_5 G_8 + 66 G_6 G_8 - 432 G_7 G_8 + 227 G_8^2 + 242 G_1 G_9 - 
214 G_{10} G_9 \\
+ &213 G_{11} G_9 - 377 G_2 G_9 - 198 G_3 G_9 + 132 G_5 G_9 + 
64 G_6 G_9 - 242 G_7 G_9 + 185 G_8 G_9 \\
+ &224 G_9^2 
+ 264 W(1,2) + 528 W(2,3) - 143 W(2,4) - 
178 W(3,4) + 264 W(1,5) \\
+ &528 W(2,5) + 528 W(3,5) - 143 W(4,5) - 
242 W(2,6) - 242 W(3,6) - 242 W(5,6) \\
- &143 W(1,7) + 35 W(3,7) - 
35 W(4,7) - 143 W(1,8) - 264 W(3,8) + 143 W(4,8) \\
 - &143 W(7,8)   ,\\
{\hat L}=  & 9 G_{10} - 3 G_{11} - 4 G_2 + G_4 + G_6 + 7 G_9.  
\endalign
$$
The following equation in $S_2^8(K(353))$ holds on all Fourier coefficients for 
$T\in {}^{277}\Xtwo$ 
satisfying $\det(2T) \le 5000$.   
If we indeed have 
$$
Q^2 + {\hat L} Q L + {\hat Q} L^2=0,  \tag \Equ{F2} 
$$
then the form $f$ is holomorphic and 
$S_2^2(K(353))(\Q)=\Span_{\Q}\left(f,G_1,\dots,G_{11}  \right)$.  
Furthermore, the form $f$ is a Hecke eigenform with spinor Euler factors
$$
\align
Q_2(f,x)  &=  1 +   x + 3 x^2 + 2 x^3 + 4 x^4,  
\\
Q_3(f,x)  &=  1 + 2 x + 4 x^2 + 6 x^3 + 9 x^4.    
\endalign
$$
\endproclaim

{\bf{\newsection{sec8}{   Examples for weights $k>2$.}} }

For weights greater than two, 
the constructions of paramodular cusp forms in this article 
become easier because the dimension formulae of Ibukiyama apply.  
Existence of a nonlift follows whenever we have 
$\dim S_2^k\left( K(p) \right) > \dim J_{k,p}^{\text{cusp}} $.  
Furthermore,  nonlifts of higher weight occur at lower prime levels 
and identities and congruences may sometimes be proven directly from 
Corollary~{\refer{D9}} and Theorem~{\refer{D10}}.   
As these nonlifts are the first examples likely to arise in related work, we give some 
examples here.  Already, there has been interest in the weight three case.  
A\. Ash, P\. Gunnells and M\. McConnell in \cite{\refer{AGM}} studied
$H^5\left( \Gamma_0(p), \C \right)$ and found cusp forms 
and computed Euler $2$ and $3$ factors at levels 
$p=61$, $73$ and $79$.  They predicted the existence of corresponding Siegel modular
cusp forms for $\Gamma_0'(p)$ and requested a construction.  
The paramodular cusp forms in $S_2^3(K(p))$ constructed here 
have macthing Euler factors at $2$ and $3$; we additionally computed the 
Euler $5$-factor.  
We thank A\. Brumer for bringing this topic to our attention.  
Recall our normalization for the $q$-Euler factor of a Hecke eigenform $ f\in S_2^k(K(p))$:    
$$
Q_q(f,x)=
1-\lambda_{q}x+(\lambda_{q}^2-\lambda_{q^2}-q^{2k-4})x^2
-\lambda_{q}\,q^{2k-3}x^3+q^{4k-6}x^4.  
$$
These Euler factors possess the symmetry 
$ x^4 q^{4k-6} Q_q\left( q^{3-2k}/x \right) = Q_q(x)$.  
The first three examples  use theta blocks of weight~$3$:
$$
\TB_3(d_1,d_2,\dots,d_{9})(\tau,z)=
\eta(\tau)^{-3}\,
\prod_{i=1}^{9}\, \vartheta(\tau, d_iz).  
$$

\subheading{Example 1} 
We have $\dim S_2^3\left( K(61) \right) =7$ and $\dim J_{3,61}^{\text{cusp}}  =6$.  
There is a nonlift Hecke eigenform $f \in S_2^3\left( K(61) \right)^{-}(\Z)$ 
with Euler factors: 
$$
\align
Q_2(f,x)  &= 1 +7x +24 x^2 +56  x^3 +64    x^4,  \\
Q_3(f,x)  &= 1 +3x +3  x^2 +81  x^3 +729   x^4,  \\
Q_5(f,x)  &= 1 -3x +85 x^2 -375 x^3 +15625 x^4.  
\endalign
$$
For $\ell=43$, $f$ is congruent to an element of 
$\Grit\left( J_{3,61}^{\text{cusp}} \right)(\Z)$ modulo $\ell$ 
and this is the only such congruence.   
We may define the nonlift $f$ via
$$
\align
f 
 &= -9 B[1] - 2 B[2] + 22 B[3] + 9 B[4] - 
      10 B[5] + 19 B[6]
- 43 B[1] B[6]/B[2],  
\endalign 
$$
where, for $1 \le i \le 6$,  the $B[i]=\Grit\left(\TB_3(\Xi_i)\right)$ are  
Gritsenko lifts of the theta blocks given by $\Xi_i=$ 
    [2, 2, 2, 3, 3, 3, 3, 5, 7],
    [2, 2, 2, 2, 3, 4, 4, 4, 7],
    [2, 2, 2, 2, 3, 3, 4, 6, 6],
    [1, 2, 3, 3, 3, 3, 4, 4, 7],
    [1, 2, 3, 3, 3, 3, 3, 6, 6],
    [1, 2, 2, 2, 4, 4, 4, 5, 6].  
The integrality of $f$ may be checked, using Theorem~{\refer{D10}}, 
by computing the $a( T; f)$ to be integers for $T\in{}^{61}\X2$ with 
$w(T) \le \frac{3}{6}\frac{61^2+1}{61+1} < 30.02$; there are $1477$ such 
${\hat \Gamma}_0(61)$-classes.  
Alternatively, we may note, as in the proof of Theorem~\refer{F3}, 
that the Fourier expansion of $B[2]$ has unit content because 
$a( \tfrac12\smtwomat{244}{22}{22}{2}; B[2])=-1$.  

\subheading{Example 2} 
We have $\dim S_2^3\left( K(73) \right) =9$ and $\dim J_{3,73}^{\text{cusp}}  =8$.  
There is a nonlift Hecke eigenform $f \in S_2^3\left( K(73) \right)^{-}(\Z)$ 
with Euler factors: 
$$
\align
Q_2(f,x)  &= 1 +6x +22  x^2 +48  x^3 +64    x^4,  \\
Q_3(f,x)  &= 1 +2x +3   x^2 +54  x^3 +729   x^4,  \\
Q_5(f,x)  &= 1     +130 x^2          +15625 x^4.  
\endalign
$$
For $\ell\in \{3,13\}$, $f$ is congruent to an element of 
$\Grit\left( J_{3,73}^{\text{cusp}} \right)(\Z)$ modulo $\ell$ 
and this is the only such congruence. 
We may define the nonlift $f= $ 
$$
\align
 & 9 B[1] + 19 B[2] + 2 B[3] - 13 B[4] + 
      34 B[5] - 15 B[6] - 12 B[7] - 10 B[8]
- 39 B[2] B[6]/B[4],  
\endalign 
$$
where the $B[i]=\Grit\left(\TB_3(\Xi_i)\right)$ are the 
Gritsenko lifts of the theta blocks given, for $1 \le i \le 8$, by $\Xi_i=$ 
[ 2, 3, 3, 3, 3, 4, 4, 5, 7 ],
    [ 2, 3, 3, 3, 3, 3, 5, 6, 6 ],
    [ 2, 2, 3, 4, 4, 4, 4, 4, 7 ],
    [ 2, 2, 3, 3, 4, 4, 4, 6, 6 ],
    [ 2, 2, 3, 3, 3, 5, 5, 5, 6 ],
    [ 2, 2, 2, 4, 4, 4, 5, 5, 6 ],
    [ 2, 2, 2, 2, 3, 4, 4, 5, 8 ],
    [ 2, 2, 2, 2, 2, 4, 5, 6, 7 ].  
From the Fourier coefficient $a( \tfrac12\smtwomat{146}{17}{17}{2}; B[4])= 1$, 
we see that the Fourier expansion of $B[4]$ has unit content, so that $f$ is integral.  

\subheading{Example 3} 
We have $\dim S_2^3\left( K(79) \right) =8$ and $\dim J_{3,79}^{\text{cusp}}  =7$.  
There is a nonlift Hecke eigenform $f \in S_2^3\left( K(79) \right)^{-}(\Z)$ 
with Euler factors: 
$$
\align
Q_2(f,x)  &= 1 +5x +14  x^2 +40   x^3 +64    x^4,  \\
Q_3(f,x)  &= 1 +5x +42  x^2 +135  x^3 +729   x^4,  \\
Q_5(f,x)  &= 1 -3x +80  x^2 -375  x^3  +15625 x^4.  
\endalign
$$
For $\ell= 2$, $f$ is congruent to an element of 
$\Grit\left( J_{3,79}^{\text{cusp}} \right)(\Z)$ modulo $\ell^5$ 
and this is the only such congruence. 
We may define the nonlift $f $ via 
$$
\align
f &=  (-32 B[1]^2 + 32 B[2]^2 + 32 B[1] B[3] - 64 B[2] B[3] + 
      32 B[3]^2 \\
&+ 26 B[1] B[4] - 38 B[2] B[4] + 19 B[3] B[4] + 3 B[4]^2 + 
      32 B[1] B[5] \\
&- 17 B[4] B[5] - 32 B[2] B[6] + 32 B[3] B[6] + 
      27 B[4] B[6] + 64 B[1] B[7] \\
&+ 64 B[2] B[7] - 96 B[3] B[7] - 
      68 B[4] B[7] - 32 B[5] B[7] - 32 B[6] B[7]  )/B[4], 
\endalign 
$$
where the $B[i]=\Grit\left(\TB_3(\Xi_i)\right)$ are the 
Gritsenko lifts of the theta blocks given, for $1 \le i \le 7$, by $\Xi_i=$ 
[ 2, 2, 3, 3, 3, 3, 5, 5, 8 ],
    [ 2, 2, 2, 3, 4, 4, 4, 5, 8 ],
    [ 2, 2, 2, 2, 4, 4, 5, 6, 7 ],
    [ 2, 2, 2, 2, 2, 4, 4, 5, 9 ],
    [ 1, 3, 3, 3, 3, 4, 4, 5, 8 ],
    [ 1, 2, 3, 4, 4, 4, 4, 4, 8 ],
    [ 1, 2, 3, 3, 3, 4, 5, 6, 7 ].  
From the Fourier coefficient $a( \tfrac12\smtwomat{1106}{47}{47}{2}; B[4])=-1$, 
we see that the Fourier expansion of $B[4]$ has unit content and so $f$ is integral.

\subheading{Example 4} 
We have $\dim S_2^4\left( K(83) \right)^{-} =1$.  
This is the lowest prime level for which a weight four minus form occurs.
There is a nonlift Hecke eigenform $f \in S_2^4\left( K(83) \right)^{-}(\Z)$ 
with Euler factors: 
$$
\align
Q_2(f,x)  &= 1 +17x +132  x^2 +544   x^3 +1024    x^4,  \\
Q_3(f,x)  &= 1 +23x +270  x^2 +5589  x^3 +59049   x^4.  
\endalign
$$
We may define the nonlift $f $ via 
$$
f=
\frac1{48}\left( 
\Tr\left( \vartheta_Q^2 \right) -
\Tr\left( \vartheta_Q^2 \right) \vert \mu
\right) 
\text{ for } 
Q=\smfourmatNUMBERFOUR{6}{8}{44}{1}{1}{2}{0}{-2}{3}.  
$$
To prove the integrality of $f$, 
we use forms of determinant $83^2$ and level $83$:   
$$
B=\smfourmatNUMBERFOUR{6}{16}{22}{1}{0}{3}{1}{1}{4}; \ \  
C=\smfourmatNUMBERTWO{2}{42}{42}{0}{-1}{0}{0}{-1}{0}; \ \ 
D=\smfourmatNUMBERTWO{12}{16}{22}{1}{-1}{3}{0}{-2}{3}.  
$$
By Theorem~\refer{B12}, the form 
$\frac1{4}\left( 
\Tr\left( \vartheta_Q \vartheta_B\right) -
\Tr\left( \vartheta_Q \vartheta_B\right) \vert \mu
\right) $ is integral; it has leading coefficient $51$.  
Also by Theorem~\refer{B12}, the form 
$\frac1{4}\left( 
\Tr\left( \vartheta_C \vartheta_D\right) -
\Tr\left( \vartheta_C \vartheta_D\right) \vert \mu
\right) $ is integral; it has leading coefficient $16 \cdot 23$. 
Some integer linear combination of the above two integral forms has leading coefficient~$1$ 
because $51$ and $16 \cdot 23$ are relatively prime.  Since the leading coefficient of 
$f$ is $1$ and the space $ S_2^4\left( K(83) \right)^{-}$ is one dimensional,  
$f$ is itself integral.  

{\bf{\newsection{sec9}{  Appendix A:  $p=499$.}} }

The proof of Theorem~\refer{A2} in \refer{sec6} 
was incomplete in the case of $p=499$ because 
we were unable to span $S_2^4\left( K(499) \right)$.   
The proof is completed in this Appendix.  
Ibukiyama's dimension formula yields $\dim S_2^4\left( K(499) \right)=503$ 
but we were only able to generate a plus subspace ${\Cal S}^{+}$ of dimension 
at least~$422$ and a minus subspace ${\Cal S}^{-}$ of dimension 
at least~$80$.  Still, a subspace ${\Cal S}={\Cal S}^{+}\oplus{\Cal S}^{-}$ 
of codimension at most one in $ S_2^4\left( K(499) \right) $ is sufficient 
information to prove that $ S_2^2\left( K(499) \right) $ is spanned by 
Gritsenko lifts. 
Recall the notation ${}^N\Xtwo$, ${}^N\Xtwo'$,   ${\Cal H}_N(2)^{\pm}$ and ${\Cal H}_N'(2)^{\pm}$ 
from section~6.  
For a prime $\ell$ and for any subspace $ Z \subseteq M_2^k\left( K(N) \right)$ 
defined over $\Q$, we set 
$Z(\F_{\ell})= R_{\ell}\left( Z \cap M_2^k\left( K(N) \right)(\Z) \right)$.  
We use the following Lemma.  

\proclaim{ \Equ{GG1} Lemma}
Let $ N \in \N$.  
Let $\ell \in \N $ be prime.  
Let $g_1$, $g_2 \in S_2^2\left( K(N) \right)^{+}(\Z)$ be linearly indepenent 
over $\F_{\ell}$.  
Let $g_3$, $g_4 \in S_2^2\left( K(N) \right)^{+}(\Z)$ be linearly indepenent 
over $\F_{\ell}$.    
Let ${\Cal S } \subseteq S_2^4 \left( K(N) \right)$ be a subspace of codimension 
at most one defined over $\Q$.  
Set ${\Cal S}^{\pm}(2) = {\Cal H}_N(2)^{\pm} \cap {\Cal S}$.  
Set ${\Cal S}^{+}(2)' = {\Cal H}_N'(2)^{+} \cap {\Cal S}$.  
If $ S_2^2\left( K(N) \right)^{+} \ne \Grit\left( J_{2,N}^{\text{cusp}} \right)   $ 
then there exist $(\alpha, \beta)$, $(\gamma, \delta)\in \Pj^1(\F_{\ell})$ 
such that the equation
$$
H_2\left( \alpha R_{\ell}(g_1) + \beta R_{\ell}(g_2) \right) = 
H_1\left( \gamma R_{\ell}(g_3) + \delta R_{\ell}(g_4) \right)
 \tag \Equ{GG2} 
$$
has a nontrivial solution 
$(H_1,H_2) \in {\Cal S}^{+}(2)'(\F_{\ell}) \times {\Cal S}^{+}(2)'(\F_{\ell})$.  
\endproclaim
\demo{Proof}
We know 
${\Cal S}^{+}(2)'$ is a rational subspace of ${\Cal H}_N'(2)^{+}$ 
with codimension at most one.  In the more difficult case, 
when the codimension equals one, there exists $ b \in {\Cal H}_N'(2)^{+}(\Q)$ 
such that ${\Cal S}^{+}(2)'$ and $b$ span ${\Cal H}_N'(2)^{+}$.  
If $ S_2^2\left( K(N) \right)^{+} \ne \Grit\left( J_{2,N}^{\text{cusp}} \right)   $ 
then there exists an $ f \in S_2^2\left( K(N) \right)^{+}(\Z) \setminus \Grit\left( J_{2,N}^{\text{cusp}} \right)   $ 
with zero as its first Fourier-Jacobi coefficient  and with $R_{\ell}(f)\ne 0$.  
For $i=1,2,3,4$, we have $g_i f \in {\Cal S}^{+}(2)'+r_i b$ 
for some $r_i \in \Q$.  Therefore, there are  pairs of relatively prime integers 
$(\alpha, \beta)$ and  $(\gamma, \delta)$ with 
$\left( \alpha g_1 + \beta g_2 \right)f$, 
$\left( \gamma g_3 + \delta g_4 \right)f \in {\Cal S}^{+}(2)'(\Z)$.  
These elements are nontrivial in ${\Cal S}^{+}(2)'(\F_{\ell})$ because the ring 
$M_2(K(N))(\F_{\ell})$ does not have zero divisors.  Calling the first image $H_1$ and the 
second $H_2$, we have a nontrivial solution to equation~\refer{GG2}.  
In the easier case, when ${\Cal S}^{+}(2)'={\Cal H}_N'(2)^{+}$, 
equation~\refer{GG2} has nontrivial solutions for every pair 
$(\alpha, \beta)$, $(\gamma, \delta)\in \Pj^1(\F_{\ell})$.  
\qed
\enddemo

\proclaim{ \Equ{GG3} Proposition}
We have 
$S_2^2\left( K(499) \right)^{-} = \{0\}$ and 
$ S_2^2\left( K(499) \right)^{+} = \Grit\left( J_{2,499}^{\text{cusp}} \right)   $.  
\endproclaim
\demo{Proof}
Since $ \dim J_{2,499}^{\text{cusp}}=18$, there are $18$ 
linearly independent Gritsenko lifts $g_1, \dots, g_{18}$.  
We will use Lemma~\refer{GG1}.  
Let $g_1$, $g_2$, $g_3$, $g_4 \in S_2^2\left( K(N) \right)^{+}(\Z)$ be choosen as the 
Gritsenko lifts of the theta blocks 
$[1,1,1,1,    1,4,4,31]$, 
$[1,1,1,1,    1,11,14,26]$,  
$[1,1,1,1,    4,5,13,28]$,  
$[1,1,1,1,    4,7,20,23]$, respectively.  
The values  for  
$\left( a\left( \smtwomat{499}{9/2}{9/2}1; g_i \right), a\left( \smtwomat{499}{0}{0}1; g_i \right) \right)$ 
are $(-1,0)$, $(1,-20)$, $(-3,0)$, $(-1,12)$, for $i=1,2,3,4$ respectively, 
 so that 
$g_1$ and $g_2 $ are linearly indepenent over $\F_{19}$,  
as are $g_3$ and $g_4$.   

Define ${\Cal S}^{+}= \Span\left( T_j(g_{\alpha}g_{\beta}): 
1 \le \alpha \le \beta \le 18, j=1,2,\dots,7\right) \subseteq  S_2^4\left( K(499) \right)^{+}$.  
Computing the action of these Hecke operators required $565,010$ Fourier coefficients.  
Define ${\Cal S}^{-}= \Span\left( \Tr( \vartheta_P\vartheta_Q)-\Tr( \vartheta_P\vartheta_Q)| \mu: 
P,Q \in {\Cal A }\right) \subseteq  S_2^4\left( K(499) \right)^{-}$ 
for the $120$ pairs $(P,Q)$ given by
$$
\align
{\Cal A}=\{    
&\smfourmatNUMBEREIGHTEEN{20}{28}{32}{1}{0}{2}{2}{5}{13}, 
\smfourmatNUMBERSIXTEEN{16}{28}{40}{1}{1}{3}{-3}{2}{10},  
\smfourmatNUMBERTWENTY{22}{24}{34}{7}{6}{8}{-4}{-2}{5}, 
\smfourmatNUMBEREIGHTEEN{24}{24}{34}{4}{5}{10}{-6}{3}{0}, 
\smfourmatNUMBERSIXTEEN{22}{26}{36}{4}{5}{10}{-2}{1}{0}, \\
&\smfourmatNUMBERSIXTEEN{18}{30}{36}{3}{6}{8}{0}{-2}{1}, 
\smfourmatNUMBERSIXTEEN{26}{28}{30}{5}{2}{8}{5}{10}{7},   
\smfourmatNUMBERFOURTEEN{24}{24}{38}{0}{-1}{2}{-2}{11}{7},
\smfourmatNUMBERSIXTEEN{26}{30}{32}{5}{5}{10}{3}{10}{13},
\smfourmatNUMBERFOURTEEN{18}{20}{52}{1}{1}{1}{-1}{2}{6},  \\
&\smfourmatNUMBEREIGHTEEN{20}{24}{34}{1}{5}{3}{-2}{-3}{6}, 
\smfourmatNUMBERTWENTY{22}{28}{32}{3}{6}{9}{-2}{9}{11}, 
\smfourmatNUMBERTWENTY{24}{26}{26}{4}{1}{2}{9}{5}{5}, 
\smfourmatNUMBERTWENTY{22}{22}{32}{0}{-3}{5}{2}{1}{9}, 
\smfourmatNUMBEREIGHTEEN{20}{24}{34}{3}{-4}{2}{-7}{-1}{3}
 \}. 
\endalign
$$
Let $E$ be the $1515$ ${\hat \Gamma}_0(499)$-classes of 
indices $T\in {}^{499}\X2$ with $\det(2T) \le 716$ and let $\pi_E$ 
be the projection onto these Fourier coefficients.  
From $\dim \pi_E \circ R_{19}( {\Cal S}^{+} ) =422$,   
we see that $\dim{\Cal S}^{+} \ge 422$. 
Let $F$ be the $610$ ${\hat \Gamma}_0(499)$-classes of 
indices $T\in {}^{499}\X2$ with $\det(2T) \le 375$ and let $\pi_F$ 
be the projection onto these Fourier coefficients.  
From $\dim \pi_F \circ R_{19}( {\Cal S}^{-} ) =80$,   
we see that $\dim{\Cal S}^{-} \ge 80$.  
Since $\dim S_2^4\left( K(499) \right)=503$ , 
${\Cal S}= {\Cal S}^{+} \oplus {\Cal S}^{-}$ is a subspace of codimension 
at most one;  
we have $\dim{\Cal S}^{+}=422 + \epsilon$ and 
$\dim{\Cal S}^{-}=80 + 1- \epsilon$ where $\epsilon = 0 $ or $1$ 
is the (unknown to us) codimension of ${\Cal S}^{+}$ in $S_2^4\left( K(499) \right)^{+}$.  
For later use, note that 
$ \dim \ker \pi_F \circ R_{19} \le 1 - \epsilon \le 1$.  

Let ${  {{\Cal S}^{+}(2)'_E}}= \{ f \in {\Cal S}^{+}: 
\forall T \in {}^{499}\X2 \setminus ( {}^{499}\X2+{}^{499}\X2'):  
[T] \in E,   a(T;f)=0 \}$.  We have 
$ {\Cal S}^{+}(2)' \subseteq {\Cal S}^{+}(2)'_E \subseteq {\Cal S}^{+} $.  
We now use Lemma~\refer{GG1} to show that 
$ S_2^2\left( K(499) \right)^{+} = \Grit\left( J_{2,499}^{\text{cusp}} \right)   $. 
We proceed by contradiction:  
if $ S_2^2\left( K(499) \right)^{+} \ne \Grit\left( J_{2,499}^{\text{cusp}} \right) $
then there exist 
$(\alpha, \beta)$, $(\gamma, \delta)\in \Pj^1(\F_{19})$ 
such that the equation~\refer{GG2} 
has a nontrivial solution 
$(H_1,H_2) \in {\Cal S}^{+}(2)'(\F_{19}) \times {\Cal S}^{+}(2)'(\F_{19})$.  
For each of the $(19+1)^2=400$ lists $( (\alpha, \beta),(\gamma, \delta)) 
\in \Pj^1(\F_{19})\times \Pj^1(\F_{19})$ we checked that, 
for 
$(H_1,H_2) \in {\Cal S}^{+}(2)'_E(\F_{19}) \times {\Cal S}^{+}(2)'_E(\F_{19})$,  
equation~\refer{GG2}  
has only the trivial solution. 
This contradiction proves the second assertion.  

Next we show that $S_2^2\left( K(499) \right)^{-} = \{0\}$.  
We calculated $\dim \pi_F \circ R_{19} \left({\Cal S}^{-}(2) \right) =9$; 
therefore $ \dim  {\Cal S}^{-}(2) \le 9 + 1-\epsilon \le 10$ 
since the kernel of $ \pi_F \circ R_{19}$ has dimension at most $1-\epsilon$.  
Since ${\Cal S}^{-}(2)\subseteq {\Cal H}^{-}_{499}(2)$ has codimension at most one, 
we have $ \dim {\Cal H}^{-}_{499}(2) \le  11$.  
We show $S_2^2\left( K(499) \right)^{-} = \{0\}$ by contradiction:  
If there were a nontrivial $ f\in S_2^2\left( K(499) \right)^{-} $, 
then there would be an $18$-dimensional subspace 
$\Span (f g_i) \subseteq {\Cal H}^{-}_{499}(2)$; however, 
${\Cal H}^{-}_{499}(2)$ is a space of dimension at most~$11$.  
Therefore $S_2^2\left( K(499) \right)^{-} = \{0\}$.  \qed
\enddemo

{\bf{\newsection{sec10}{  Appendix B:  Tables of Fourier coefficients.}} }

In this Appendix, we explain how to read the Tables of Fourier coeffcients.  
Each form $ f\in S_2^k\left( K(N) \right)$ has a Fourier expansion 
$$
f(\Omega) = \sum_{T \in {}^N\Xtwo} a(T)\, e\left(\<\Omega, T\> \right) .
$$
For $T = \tfrac12 \smtwomat{a}bbc $, the column ``Unreduced form" gives $a,b,c$.  
The column ``Coeff" gives the value of $ a\left( \tfrac12 \smtwomat{a}bbc \right)$.  
The ``Det" column gives $\det(2T)$.  
For example,  the ninth entry row of the Table for $p=277$ gives
$$
a\left( \tfrac12 \smtwomat{277554}{1825}{1825}{12} \right)=-5  
\text{ and } 
\det    \smtwomat{277554}{1825}{1825}{12}  = 23.   
$$
Note that $277554= 1002 \cdot 277$ and so 
$$
\tfrac12 \smtwomat{277554}{1825}{1825}{12} \in {}^{277}\Xtwo.  
$$

%
%
For $ U \in {\hat \Gamma}_0(N)= \{   \smtwomat{a}bcd \in \GL_2(\Z) : 
N|c \}$, the Fourier coefficients satisfy 
$$
a( U' T U) = \det(U)^k\, a(T).  
$$
Since $k=2$ in these tables, the Fourier coefficients are ${\hat \Gamma}_0(N)$-invariant.  
The ``Unreduced form" gives one representative from each 
${\hat \Gamma}_0(N)$-class in ${}^{N}\Xtwo$. 
The entries are in order of increasing determinant so that one can see when the 
${\hat \Gamma}_0(N)$-classes from ${}^{N}\Xtwo$ for a particular determinant
are exhausted.  Since $D=\det(2T)=ac-b^2 \equiv -b^2 \mod N$, 
only determinants $D$ with $-D$ a quadratic residue modulo $N$ appear.

\subheading{ Plus and Minus Spaces}  
Each eigenform in the Tables satisfies $ f | \mu = \epsilon f$ with $\epsilon = \pm 1$.  
The Fourier coefficients then satisfy $ a(\Twin(T))= \epsilon a(T)$, where 
 $\Twin(T)= \frac12 \smtwomat{Nc}{-b}{-b}{a/N}$. 
For example, for $D=23$, we have seen 
$$
a\left( \tfrac12 \smtwomat{277\cdot 1002}{1825}{1825}{12} \right)=-5 .     
$$
Since the nonlift  for $p=277$ is in the $\mu$-plus space we also have 
$$
a\left( \tfrac12 \smtwomat{277\cdot 12}{-1825}{-1825}{1002} \right)=-5 .     
$$
However, since the entries in the ``Unreduced form" column are noncanonical, 
it is not immediately clear which entry this corresponds to.  

To remedy this we let $\GL_2(\Z)$ act on $\Xtwo \times \Pj(\Z/N\Z)$ as 
$(T,v) \mapsto (U'TU, U\inv v)$.  
Recall that 
$
\Pj(\Z/N\Z)= \{ \pmatrix a \\ b \endpmatrix \in \Z^2: \operatorname{gcd}(a,b)=1 \}/\sim 
$,  
where 
$\pmatrix a \\ b \endpmatrix \sim \pmatrix c \\ d \endpmatrix$ 
if and only if 
$\exists \mu \in \Z:$ $\mu a \equiv c \mod N$ and $\mu b \equiv d \mod N$.  
If we identify ${}^N\Xtwo$ with 
${}^N\Xtwo \times \{ \pmatrix 1 \\ 0 \endpmatrix \} \subseteq \Xtwo \times \Pj(\Z/N\Z)$, 
then for $T_1,T_2 \in {}^N\Xtwo$, we have that 
$(T_1, \pmatrix 1 \\ 0 \endpmatrix)$ is $\GL_2(\Z)$-equivalent to 
$(T_2, \pmatrix 1 \\ 0 \endpmatrix)$ if and only if $T_1$ is 
$ {\hat \Gamma}_0(N) $-equivalent to $T_2$.  
For each  $T_1 \in {}^N\Xtwo$, we $\GL_2(\Z)$ reduce $(T_1, \pmatrix 1 \\ 0 \endpmatrix)$ 
to a $(T,v)$ with $T$ Legendre reduced, that is $ 0 \le 2b \le a \le c$.  
The column ``Reduced form" lists $a,b,c \quad \alpha,\beta$ for 
$(T,v)= \left( \frac12 \smtwomat{a}{b}{b}{c}, \pmatrix \alpha \\ \beta \endpmatrix \right)$.  
This is canonical except for the small number of choices for $v$: 
We have the  $\GL_2(\Z)$-equivalence 
$(T,v_1) \sim (T,v_2) \iff \exists \,U \in \Aut_{\Z}(T):\, v_1=Uv_2$.  

Going back to the example, we reduce 
$$
\left(  \tfrac12 \smtwomat{277\cdot 12}{-1825}{-1825}{1002} , \pmatrix 1 \\ 0 \endpmatrix \right) \sim 
\left(  \tfrac12 \smtwomat{4}{1}{1}{6} , \pmatrix 1 \\ 27 \endpmatrix \right).  
$$
To see this, use $U=\smtwomat{11}{-6}{20}{-11}\in \GL_2(\Z)$ so that 
$$
U' \smtwomat{277\cdot 12}{-1825}{-1825}{1002} U=
\smtwomat{11}{20}{-6}{-11}
 \smtwomat{277\cdot 12}{-1825}{-1825}{1002}
\smtwomat{11}{-6}{20}{-11}=
\smtwomat{4}{1}{1}{6} \text{ and   }
$$
$$
U\inv \pmatrix 1 \\ 0 \endpmatrix= 
\smtwomat{11}{-6}{20}{-11} \pmatrix 1 \\ 0 \endpmatrix= 
\pmatrix 11 \\ 20 \endpmatrix \equiv \pmatrix 1 \\ 27 \endpmatrix 
\text{ in } \Pj(\F_{277}) 
$$
since $\vmatrix 11 & 1 \\ 20 & 27 \endvmatrix = 277$.  
This explains why the Fourier coefficient for the tenth entry 
$D=23$; and Reduced form ``$4,1,6 \quad 1,27$" is also $-5$.  

If one wishes to move from a reduced form $(T,v)$ to an 
unreduced form, one simply takes $U'TU$ for a 
$U \in \GL_2(\Z)$ whose first column is $v \mod N$.  
For example, for $D=23$ again, take the reduced form from row eleven:  
$$
(T,v)=\left( \tfrac12 \smtwomat4116 , \pmatrix 1 \\ 65 \endpmatrix\right).  
$$
For $U=\smtwomat10{65}1$ we have 
$$
U' T U= 
\smtwomat1{65}{0}1 \smtwomat4116 \smtwomat10{65}1= 
\smtwomat{25484}{391}{391}{6},
$$
which is exactly the unreduced entry because  
the ``Unreduced form" column was constructed in just this way.  
In each of the following tables, 
for each determinant~$D$, 
either all or none of the ${\hat \Gamma}_0(N)$-classes with determinant~$D$ 
are represented. 
For $(p,\epsilon)=(277,+)$, $(349,+)$, $(353,+)$, $(389,+)$, $(461,+)$, $(523,+)$, 
$(587,-)$ and $(587,+)$, the following short tables give Fourier coefficients for 
a nonlift $f \in S_2^2\left( K(p) \right)^{\epsilon}(\Z)$.  
Larger tables of Fourier coefficients are at \cite{\refer{URL}}.

\newpage
\centerline{LEVEL   277    }
\smallskip
\begintable
\   Coeff \ &\  Det   \ &\    Reduced form   \ &\   Unreduced form\crthick
$     -3$ &  $   3$  & $ 2, 1, 2; 1, 116$ &    $27146, 233, 2$\nr
$     -2$ &  $   4$  & $ 2, 0, 2;  1, 60$ &    $7202, 120, 2$\nr
$     -1$ &  $   7$  & $ 2, 1, 4; 1, 150$ &    $90302, 601, 4$\nr
$     -3$ &  $  12$  & $ 2, 0, 6; 1, 107$ &    $68696, 642, 6$\nr
$      6$ &  $  12$  & $ 4, 2, 4; 1, 116$ &    $54292, 466, 4$\cr
$     -5$ &  $  16$  & $ 2, 0, 8;  1, 30$ &    $7202, 240, 8$\nr
$      6$ &  $  16$  & $ 4, 0, 4;  1, 60$ &    $14404, 240, 4$\nr
$     -2$ &  $  19$  & $2, 1, 10;  1, 52$ &    $27146, 521, 10$\nr
$     -5$ &  $  23$  & $2, 1, 12; 1, 152$ &    $277554, 1825, 12$\nr
$     -5$ &  $  23$  & $ 4, 1, 6;  1, 27$ &    $4432, 163, 6$\cr
$     10$ &  $  23$  & $ 4, 1, 6;  1, 65$ &    $25484, 391, 6$\nr
$      0$ &  $  27$  & $2, 1, 14; 1, 129$ &    $233234, 1807, 14$\nr
$      6$ &  $  27$  & $ 6, 3, 6; 1, 116$ &    $81438, 699, 6$\nr
$     -9$ &  $  28$  & $2, 0, 14; 1, 112$ &    $175618, 1568, 14$\nr
$      4$ &  $  28$  & $ 4, 2, 8; 1, 150$ &    $180604, 1202, 8$\cr
$      8$ &  $  36$  & $2, 0, 18;  1, 20$ &    $7202, 360, 18$\nr
$     -7$ &  $  36$  & $4, 2, 10; 1, 130$ &    $169524, 1302, 10$\nr
$      2$ &  $  36$  & $ 6, 0, 6;  1, 60$ &    $21606, 360, 6$\nr
$      4$ &  $  39$  & $2, 1, 20;  1, 21$ &    $8864, 421, 20$\nr
$      4$ &  $  39$  & $4, 1, 10;  1, 42$ &    $17728, 421, 10$\cr
$      4$ &  $  39$  & $4, 1, 10; 1, 124$ &    $154012, 1241, 10$\nr
$    -11$ &  $  39$  & $ 6, 3, 8; 1, 224$ &    $402758, 1795, 8$\nr
$     -3$ &  $  40$  & $2, 0, 20; 1, 103$ &    $212182, 2060, 20$\nr
$     -3$ &  $  40$  & $4, 0, 10;  1, 71$ &    $50414, 710, 10$
\endtable
\begintable
\   Coeff \ &\  Det   \ &\    Reduced form   \ &\   Unreduced form\crthick
$      3$ &  $  47$  & $2, 1, 24; 1, 105$ &    $264812, 2521, 24$\cr
$      3$ &  $  47$  & $4, 1, 12; 1, 113$ &    $153458, 1357, 12$\nr
$    -12$ &  $  47$  & $4, 1, 12; 1, 210$ &    $529624, 2521, 12$\nr
$      3$ &  $  47$  & $ 6, 1, 8;  1, 31$ &    $7756, 249, 8$\nr
$      3$ &  $  47$  & $ 6, 1, 8;  1, 38$ &    $11634, 305, 8$\nr
$      6$ &  $  48$  & $2, 0, 24;  1, 85$ &    $173402, 2040, 24$\cr
$      0$ &  $  48$  & $4, 0, 12; 1, 107$ &    $137392, 1284, 12$\nr
$     -9$ &  $  48$  & $ 6, 0, 8;  1, 22$ &    $3878, 176, 8$\nr
$      3$ &  $  48$  & $ 8, 4, 8; 1, 116$ &    $108584, 932, 8$\nr
$     -5$ &  $  52$  & $2, 0, 26;  1, 74$ &    $142378, 1924, 26$\nr
$     10$ &  $  52$  & $4, 2, 14; 1, 179$ &    $449294, 2508, 14$\cr
$      7$ &  $  55$  & $2, 1, 28;  1, 16$ &    $7202, 449, 28$\nr
$     -8$ &  $  55$  & $4, 1, 14;  1, 32$ &    $14404, 449, 14$\nr
$      7$ &  $  55$  & $4, 1, 14;  1, 47$ &    $31024, 659, 14$\nr
$     -8$ &  $  55$  & $ 8, 3, 8;  1, 82$ &    $54292, 659, 8$\nr
$      1$ &  $  59$  & $2, 1, 30;  1, 44$ &    $58170, 1321, 30$\cr
$      1$ &  $  59$  & $6, 1, 10;  1, 34$ &    $11634, 341, 10$\nr
$      1$ &  $  59$  & $6, 1, 10; 1, 132$ &    $174510, 1321, 10$\nr
$      4$ &  $  63$  & $2, 1, 32;  1, 39$ &    $48752, 1249, 32$\nr
$      4$ &  $  63$  & $4, 1, 16;  1, 78$ &    $97504, 1249, 16$\nr
$      4$ &  $  63$  & $4, 1, 16;  1, 95$ &    $144594, 1521, 16$\cr
$      1$ &  $  63$  & $6, 3, 12; 1, 150$ &    $270906, 1803, 12$\nr
$    -11$ &  $  63$  & $ 8, 1, 8; 1, 156$ &    $195008, 1249, 8$\nr
$     -1$ &  $  64$  & $2, 0, 32;  1, 15$ &    $7202, 480, 32$\nr
$      4$ &  $  64$  & $4, 0, 16;  1, 30$ &    $14404, 480, 16$\nr
$     -4$ &  $  64$  & $ 8, 0, 8;  1, 60$ &    $28808, 480, 8$\cr
$     -1$ &  $  64$  & $8, 4, 10;   1, 7$ &    $554, 74, 10$
\endtable

\newpage
\centerline{LEVEL   349    }
\smallskip
\begintable
\   Coeff \ &\  Det   \ &\    Reduced form   \ &\   Unreduced form\crthick
$      3$ &  $   3$  & $ 2, 1, 2; 1, 122$ &    $30014, 245, 2$\nr
$      2$ &  $   4$  & $ 2, 0, 2; 1, 136$ &    $36994, 272, 2$\nr
$      1$ &  $  12$  & $ 2, 0, 6; 1, 151$ &    $136808, 906, 6$\nr
$     -6$ &  $  12$  & $ 4, 2, 4; 1, 122$ &    $60028, 490, 4$\nr
$      6$ &  $  15$  & $ 2, 1, 8; 1, 110$ &    $97022, 881, 8$\cr
$     -7$ &  $  15$  & $ 4, 1, 4; 1, 220$ &    $194044, 881, 4$\nr
$      3$ &  $  16$  & $ 2, 0, 8;  1, 68$ &    $36994, 544, 8$\nr
$     -6$ &  $  16$  & $ 4, 0, 4; 1, 136$ &    $73988, 544, 4$\nr
$      4$ &  $  19$  & $2, 1, 10;  1, 22$ &    $4886, 221, 10$\nr
$      7$ &  $  20$  & $2, 0, 10;  1, 81$ &    $65612, 810, 10$\cr
$     -6$ &  $  20$  & $ 4, 2, 6; 1, 251$ &    $379014, 1508, 6$\nr
$      0$ &  $  23$  & $2, 1, 12; 1, 148$ &    $263146, 1777, 12$\nr
$      0$ &  $  23$  & $ 4, 1, 6; 1, 169$ &    $171708, 1015, 6$\nr
$      0$ &  $  23$  & $ 4, 1, 6; 1, 296$ &    $526292, 1777, 6$\nr
$     -2$ &  $  27$  & $2, 1, 14;  1, 97$ &    $131922, 1359, 14$\cr
$     -6$ &  $  27$  & $ 6, 3, 6; 1, 122$ &    $90042, 735, 6$\nr
$     -9$ &  $  31$  & $2, 1, 16;  1, 54$ &    $46766, 865, 16$\nr
$      4$ &  $  31$  & $ 4, 1, 8; 1, 108$ &    $93532, 865, 8$\nr
$      4$ &  $  31$  & $ 4, 1, 8; 1, 328$ &    $861332, 2625, 8$\nr
$      5$ &  $  36$  & $2, 0, 18;  1, 71$ &    $90740, 1278, 18$\cr
$     -8$ &  $  36$  & $4, 2, 10; 1, 221$ &    $489298, 2212, 10$\nr
$     -2$ &  $  36$  & $ 6, 0, 6; 1, 136$ &    $110982, 816, 6$\nr
$      2$ &  $  48$  & $2, 0, 24;  1, 99$ &    $235226, 2376, 24$\nr
$      4$ &  $  48$  & $4, 0, 12; 1, 151$ &    $273616, 1812, 12$\nr
$    -11$ &  $  48$  & $ 6, 0, 8;  1, 52$ &    $21638, 416, 8$\cr
$      3$ &  $  48$  & $ 8, 4, 8; 1, 122$ &    $120056, 980, 8$
\endtable

\begintable
\   Coeff \ &\  Det   \ &\    Reduced form   \ &\   Unreduced form\crthick
$     -4$ &  $  51$  & $2, 1, 26;  1, 43$ &    $48162, 1119, 26$\nr
$      9$ &  $  51$  & $6, 3, 10; 1, 167$ &    $279898, 1673, 10$\nr
$     -6$ &  $  56$  & $2, 0, 28;  1, 18$ &    $9074, 504, 28$\nr
$      7$ &  $  56$  & $4, 0, 14;  1, 36$ &    $18148, 504, 14$\cr
$     -6$ &  $  56$  & $6, 2, 10;  1, 89$ &    $79572, 892, 10$\nr
$      7$ &  $  56$  & $6, 2, 10; 1, 120$ &    $144486, 1202, 10$\nr
$    -10$ &  $  60$  & $2, 0, 30;  1, 82$ &    $201722, 2460, 30$\nr
$      2$ &  $  60$  & $4, 2, 16; 1, 110$ &    $194044, 1762, 16$\nr
$      3$ &  $  60$  & $6, 0, 10; 1, 103$ &    $106096, 1030, 10$\cr
$      2$ &  $  60$  & $ 8, 2, 8; 1, 220$ &    $388088, 1762, 8$\nr
$     -6$ &  $  64$  & $2, 0, 32;  1, 34$ &    $36994, 1088, 32$\nr
$      0$ &  $  64$  & $4, 0, 16;  1, 68$ &    $73988, 1088, 16$\nr
$      8$ &  $  64$  & $ 8, 0, 8; 1, 136$ &    $147976, 1088, 8$\nr
$     -6$ &  $  64$  & $8, 4, 10; 1, 170$ &    $290368, 1704, 10$\cr
$     13$ &  $  67$  & $2, 1, 34; 1, 190$ &    $1227782, 6461, 34$\nr
$      4$ &  $  68$  & $2, 0, 34;  1, 75$ &    $191252, 2550, 34$\nr
$     -9$ &  $  68$  & $4, 2, 18;  1, 64$ &    $73988, 1154, 18$\nr
$      4$ &  $  68$  & $6, 2, 12;  1, 20$ &    $4886, 242, 12$\nr
$      4$ &  $  68$  & $6, 2, 12;  1, 96$ &    $110982, 1154, 12$\cr
$     -5$ &  $  75$  & $2, 1, 38; 1, 271$ &    $2791302, 10299, 38$\nr
$      8$ &  $  75$  & $6, 3, 14;  1, 12$ &    $2094, 171, 14$\nr
$     -3$ &  $  75$  & $10, 5, 10; 1, 122$ &    $150070, 1225, 10$\nr
$     -3$ &  $  76$  & $2, 0, 38;  1, 30$ &    $34202, 1140, 38$\nr
$     -8$ &  $  76$  & $4, 2, 20;  1, 22$ &    $9772, 442, 20$\cr
$     -3$ &  $  76$  & $8, 2, 10;  1, 44$ &    $19544, 442, 10$\nr
$     10$ &  $  76$  & $8, 2, 10; 1, 165$ &    $272918, 1652, 10$
\endtable

\newpage
\centerline{LEVEL   353    }
\smallskip
\begintable
\   Coeff \ &\  Det   \ &\    Reduced form   \ &\   Unreduced form\crthick
$      2$ &  $   4$  & $ 2, 0, 2;  1, 42$ &    $3530, 84, 2$\nr
$      1$ &  $   8$  & $ 2, 0, 4;  1, 23$ &    $2118, 92, 4$\nr
$      1$ &  $  11$  & $ 2, 1, 6; 1, 246$ &    $363590, 1477, 6$\nr
$      0$ &  $  15$  & $ 2, 1, 8; 1, 183$ &    $268280, 1465, 8$\nr
$      0$ &  $  15$  & $ 4, 1, 4;  1, 13$ &    $706, 53, 4$\cr
$      3$ &  $  16$  & $ 2, 0, 8;  1, 21$ &    $3530, 168, 8$\nr
$     -4$ &  $  16$  & $ 4, 0, 4;  1, 42$ &    $7060, 168, 4$\nr
$     -3$ &  $  19$  & $2, 1, 10; 1, 238$ &    $566918, 2381, 10$\nr
$      4$ &  $  23$  & $2, 1, 12;  1, 65$ &    $50832, 781, 12$\nr
$      4$ &  $  23$  & $ 4, 1, 6; 1, 105$ &    $66364, 631, 6$\cr
$     -7$ &  $  23$  & $ 4, 1, 6; 1, 130$ &    $101664, 781, 6$\nr
$     -4$ &  $  32$  & $2, 0, 16; 1, 165$ &    $435602, 2640, 16$\nr
$     -2$ &  $  32$  & $ 4, 0, 8;  1, 23$ &    $4236, 184, 8$\nr
$      7$ &  $  32$  & $ 6, 2, 6; 1, 148$ &    $132022, 890, 6$\nr
$      1$ &  $  35$  & $2, 1, 18;  1, 57$ &    $58598, 1027, 18$\cr
$      1$ &  $  35$  & $ 6, 1, 6;  1, 64$ &    $24710, 385, 6$\nr
$      1$ &  $  36$  & $2, 0, 18;  1, 14$ &    $3530, 252, 18$\nr
$      1$ &  $  36$  & $4, 2, 10;  1, 25$ &    $6354, 252, 10$\nr
$     -4$ &  $  36$  & $ 6, 0, 6;  1, 42$ &    $10590, 252, 6$\nr
$     -5$ &  $  39$  & $2, 1, 20;  1, 56$ &    $62834, 1121, 20$\cr
$     -5$ &  $  39$  & $4, 1, 10;  1, 29$ &    $8472, 291, 10$\nr
$      6$ &  $  39$  & $4, 1, 10; 1, 112$ &    $125668, 1121, 10$\nr
$      6$ &  $  39$  & $ 6, 3, 8;  1, 36$ &    $10590, 291, 8$\nr
$     -5$ &  $  43$  & $2, 1, 22; 1, 130$ &    $372062, 2861, 22$
\endtable

\begintable
\   Coeff \ &\  Det   \ &\    Reduced form   \ &\   Unreduced form\crthick
$     -3$ &  $  44$  & $2, 0, 22;  1, 38$ &    $31770, 836, 22$\cr
$     -1$ &  $  44$  & $4, 2, 12; 1, 246$ &    $727180, 2954, 12$\nr
$      8$ &  $  44$  & $ 6, 2, 8;  1, 16$ &    $2118, 130, 8$\nr
$     -3$ &  $  44$  & $ 6, 2, 8; 1, 160$ &    $205446, 1282, 8$\nr
$     -2$ &  $  47$  & $2, 1, 24; 1, 215$ &    $1109832, 5161, 24$\nr
$     -2$ &  $  47$  & $4, 1, 12;  1, 77$ &    $71306, 925, 12$\cr
$     -2$ &  $  47$  & $4, 1, 12; 1, 217$ &    $565506, 2605, 12$\nr
$      9$ &  $  47$  & $ 6, 1, 8; 1, 149$ &    $177912, 1193, 8$\nr
$     -2$ &  $  47$  & $ 6, 1, 8; 1, 292$ &    $682702, 2337, 8$\nr
$      0$ &  $  60$  & $2, 0, 30;  1, 20$ &    $12002, 600, 30$\nr
$      0$ &  $  60$  & $4, 2, 16; 1, 183$ &    $536560, 2930, 16$\cr
$      0$ &  $  60$  & $6, 0, 10;  1, 60$ &    $36006, 600, 10$\nr
$      0$ &  $  60$  & $ 8, 2, 8;  1, 13$ &    $1412, 106, 8$\nr
$     -5$ &  $  64$  & $2, 0, 32; 1, 166$ &    $881794, 5312, 32$\nr
$      1$ &  $  64$  & $4, 0, 16;  1, 21$ &    $7060, 336, 16$\nr
$     -2$ &  $  64$  & $ 8, 0, 8;  1, 42$ &    $14120, 336, 8$\cr
$      6$ &  $  64$  & $8, 4, 10; 1, 245$ &    $602218, 2454, 10$\nr
$      8$ &  $  68$  & $2, 0, 34;  1, 78$ &    $206858, 2652, 34$\nr
$     -3$ &  $  68$  & $4, 2, 18; 1, 108$ &    $210388, 1946, 18$\nr
$     -3$ &  $  68$  & $6, 2, 12;  1, 73$ &    $64246, 878, 12$\nr
$     -3$ &  $  68$  & $6, 2, 12; 1, 162$ &    $315582, 1946, 12$\cr
$      4$ &  $  72$  & $2, 0, 36; 1, 110$ &    $435602, 3960, 36$\nr
$      4$ &  $  72$  & $4, 0, 18; 1, 133$ &    $318406, 2394, 18$\nr
$     -4$ &  $  72$  & $6, 0, 12;  1, 23$ &    $6354, 276, 12$\nr
$     -2$ &  $  76$  & $2, 0, 38;  1, 51$ &    $98840, 1938, 38$\nr
$      3$ &  $  76$  & $4, 2, 20; 1, 238$ &    $1133836, 4762, 20$\cr
$     -2$ &  $  76$  & $8, 2, 10; 1, 123$ &    $151790, 1232, 10$\nr
$     -2$ &  $  76$  & $8, 2, 10; 1, 159$ &    $253454, 1592, 10$
\endtable

\centerline{LEVEL   389    }
\smallskip
\begintable
\   Coeff \ &\  Det   \ &\    Reduced form   \ &\   Unreduced form\crthick
$     -2$ &  $   4$  & $ 2, 0, 2; 1, 115$ &    $26452, 230, 2$\nr
$      1$ &  $   7$  & $ 2, 1, 4;  1, 46$ &    $8558, 185, 4$\nr
$     -1$ &  $  11$  & $ 2, 1, 6;  1, 59$ &    $21006, 355, 6$\nr
$     -2$ &  $  16$  & $ 2, 0, 8; 1, 137$ &    $150154, 1096, 8$\nr
$      4$ &  $  16$  & $ 4, 0, 4; 1, 115$ &    $52904, 460, 4$\cr
$     -3$ &  $  19$  & $2, 1, 10;  1, 44$ &    $19450, 441, 10$\nr
$     -5$ &  $  20$  & $2, 0, 10;  1, 33$ &    $10892, 330, 10$\nr
$      5$ &  $  20$  & $ 4, 2, 6; 1, 204$ &    $250516, 1226, 6$\nr
$     -1$ &  $  24$  & $2, 0, 12;  1, 18$ &    $3890, 216, 12$\nr
$     -1$ &  $  24$  & $ 4, 0, 6;  1, 36$ &    $7780, 216, 6$\cr
$      5$ &  $  28$  & $2, 0, 14;  1, 82$ &    $94138, 1148, 14$\nr
$     -3$ &  $  28$  & $ 4, 2, 8;  1, 46$ &    $17116, 370, 8$\nr
$      4$ &  $  35$  & $2, 1, 18;  1, 41$ &    $30342, 739, 18$\nr
$     -6$ &  $  35$  & $ 6, 1, 6; 1, 123$ &    $91026, 739, 6$\nr
$      0$ &  $  36$  & $2, 0, 18; 1, 168$ &    $508034, 3024, 18$\cr
$      0$ &  $  36$  & $4, 2, 10; 1, 242$ &    $586612, 2422, 10$\nr
$      4$ &  $  36$  & $ 6, 0, 6; 1, 115$ &    $79356, 690, 6$\nr
$      3$ &  $  44$  & $2, 0, 22; 1, 103$ &    $233400, 2266, 22$\nr
$      1$ &  $  44$  & $4, 2, 12;  1, 59$ &    $42012, 710, 12$\nr
$      3$ &  $  44$  & $ 6, 2, 8; 1, 283$ &    $641850, 2266, 8$\cr
$     -7$ &  $  44$  & $ 6, 2, 8; 1, 300$ &    $721206, 2402, 8$\nr
$      3$ &  $  52$  & $2, 0, 26;  1, 48$ &    $59906, 1248, 26$\nr
$     -7$ &  $  52$  & $4, 2, 14;  1, 89$ &    $111254, 1248, 14$ 
\endtable

\begintable
\   Coeff \ &\  Det   \ &\    Reduced form   \ &\   Unreduced form\crthick
$      3$ &  $  55$  & $2, 1, 28; 1, 285$ &    $2274872, 7981, 28$\nr
$      3$ &  $  55$  & $4, 1, 14; 1, 181$ &    $459020, 2535, 14$\cr
$      3$ &  $  55$  & $4, 1, 14; 1, 319$ &    $1425296, 4467, 14$\nr
$     -7$ &  $  55$  & $ 8, 3, 8; 1, 122$ &    $119812, 979, 8$ \nr
$     -1$ &  $  59$  & $2, 1, 30; 1, 166$ &    $827014, 4981, 30$\nr
$     -1$ &  $  59$  & $6, 1, 10; 1, 109$ &    $119034, 1091, 10$\nr
$     -1$ &  $  59$  & $6, 1, 10; 1, 202$ &    $408450, 2021, 10$\cr
$     -5$ &  $  63$  & $2, 1, 32; 1, 323$ &    $3339176, 10337, 32$\nr
$     -5$ &  $  63$  & $4, 1, 16; 1, 257$ &    $1057302, 4113, 16$\nr
$      5$ &  $  63$  & $4, 1, 16; 1, 375$ &    $2250754, 6001, 16$\nr
$     -2$ &  $  63$  & $6, 3, 12;  1, 46$ &    $25674, 555, 12$\nr
$      5$ &  $  63$  & $ 8, 1, 8; 1, 125$ &    $125258, 1001, 8$\cr
$      2$ &  $  64$  & $2, 0, 32; 1, 126$ &    $508034, 4032, 32$\nr
$     -2$ &  $  64$  & $4, 0, 16; 1, 137$ &    $300308, 2192, 16$\nr
$      0$ &  $  64$  & $ 8, 0, 8; 1, 115$ &    $105808, 920, 8$\nr
$      2$ &  $  64$  & $8, 4, 10; 1, 141$ &    $199946, 1414, 10$\nr
$     -2$ &  $  67$  & $2, 1, 34; 1, 269$ &    $2460814, 9147, 34$\cr
$      3$ &  $  68$  & $2, 0, 34;  1, 31$ &    $32676, 1054, 34$\nr
$      3$ &  $  68$  & $4, 2, 18;  1, 71$ &    $91026, 1280, 18$\nr
$      3$ &  $  68$  & $6, 2, 12; 1, 301$ &    $1088422, 3614, 12$\nr
$     -7$ &  $  68$  & $6, 2, 12; 1, 347$ &    $1446302, 4166, 12$
\endtable

\newpage
\centerline{LEVEL   461    }
\smallskip
\begintable
\   Coeff \ &\  Det   \ &\    Reduced form   \ &\   Unreduced form\crthick
$      2$ &  $   4$  & $ 2, 0, 2;  1, 48$ &    $4610, 96, 2$\nr
$      1$ &  $  16$  & $ 2, 0, 8;  1, 24$ &    $4610, 192, 8$\nr
$     -2$ &  $  16$  & $ 4, 0, 4;  1, 48$ &    $9220, 192, 4$\nr
$     -1$ &  $  19$  & $2, 1, 10; 1, 202$ &    $408446, 2021, 10$\nr
$     -3$ &  $  20$  & $2, 0, 10;  1, 44$ &    $19362, 440, 10$\cr
$      4$ &  $  20$  & $ 4, 2, 6;  1, 73$ &    $32270, 440, 6$\nr
$     -2$ &  $  23$  & $2, 1, 12; 1, 104$ &    $130002, 1249, 12$\nr
$     -2$ &  $  23$  & $ 4, 1, 6;  1, 99$ &    $59008, 595, 6$\nr
$      5$ &  $  23$  & $ 4, 1, 6; 1, 208$ &    $260004, 1249, 6$\nr
$     -1$ &  $  24$  & $2, 0, 12;  1, 98$ &    $115250, 1176, 12$\cr
$     -1$ &  $  24$  & $ 4, 0, 6; 1, 196$ &    $230500, 1176, 6$\nr
$      2$ &  $  36$  & $2, 0, 18;  1, 16$ &    $4610, 288, 18$\nr
$      2$ &  $  36$  & $4, 2, 10; 1, 213$ &    $454546, 2132, 10$\nr
$     -6$ &  $  36$  & $ 6, 0, 6;  1, 48$ &    $13830, 288, 6$\nr
$      4$ &  $  39$  & $2, 1, 20; 1, 133$ &    $354048, 2661, 20$\cr
$     -3$ &  $  39$  & $4, 1, 10; 1, 266$ &    $708096, 2661, 10$\nr
$      4$ &  $  39$  & $4, 1, 10; 1, 287$ &    $824268, 2871, 10$\nr
$     -3$ &  $  39$  & $ 6, 3, 8; 1, 128$ &    $131846, 1027, 8$\nr
$      4$ &  $  43$  & $2, 1, 22;  1, 31$ &    $21206, 683, 22$\nr
$      0$ &  $  56$  & $2, 0, 28; 1, 162$ &    $734834, 4536, 28$\cr
$      0$ &  $  56$  & $4, 0, 14; 1, 137$ &    $262770, 1918, 14$\nr
$      0$ &  $  56$  & $6, 2, 10; 1, 269$ &    $724692, 2692, 10$\nr
$      0$ &  $  56$  & $6, 2, 10; 1, 376$ &    $1415270, 3762, 10$\nr
$      4$ &  $  59$  & $2, 1, 30;  1, 54$ &    $87590, 1621, 30$\nr
$     -3$ &  $  59$  & $6, 1, 10; 1, 162$ &    $262770, 1621, 10$\cr
$     -3$ &  $  59$  & $6, 1, 10; 1, 391$ &    $1529598, 3911, 10$
\endtable

\begintable
\   Coeff \ &\  Det   \ &\    Reduced form   \ &\   Unreduced form\crthick
$      0$ &  $  64$  & $2, 0, 32;  1, 12$ &    $4610, 384, 32$\nr
$      2$ &  $  64$  & $4, 0, 16;  1, 24$ &    $9220, 384, 16$\nr
$     -4$ &  $  64$  & $ 8, 0, 8;  1, 48$ &    $18440, 384, 8$\nr
$      0$ &  $  64$  & $8, 4, 10;  1, 38$ &    $14752, 384, 10$\cr
$      2$ &  $  67$  & $2, 1, 34; 1, 347$ &    $4094602, 11799, 34$\nr
$     -2$ &  $  68$  & $2, 0, 34;  1, 76$ &    $196386, 2584, 34$\nr
$      5$ &  $  68$  & $4, 2, 18;  1, 10$ &    $1844, 182, 18$\nr
$     -2$ &  $  68$  & $6, 2, 12;  1, 15$ &    $2766, 182, 12$\nr
$     -2$ &  $  68$  & $6, 2, 12; 1, 292$ &    $1024342, 3506, 12$\cr
$      2$ &  $  76$  & $2, 0, 38; 1, 112$ &    $476674, 4256, 38$\nr
$      0$ &  $  76$  & $4, 2, 20; 1, 202$ &    $816892, 4042, 20$\nr
$     -5$ &  $  76$  & $8, 2, 10; 1, 241$ &    $581782, 2412, 10$\nr
$      2$ &  $  76$  & $8, 2, 10; 1, 404$ &    $1633784, 4042, 10$\nr
$      2$ &  $  80$  & $2, 0, 40;  1, 22$ &    $19362, 880, 40$\cr
$     -4$ &  $  80$  & $4, 0, 20;  1, 44$ &    $38724, 880, 20$\nr
$      2$ &  $  80$  & $6, 2, 14; 1, 392$ &    $2152870, 5490, 14$\nr
$      2$ &  $  80$  & $6, 2, 14; 1, 398$ &    $2219254, 5574, 14$\nr
$     -5$ &  $  80$  & $8, 0, 10;  1, 88$ &    $77448, 880, 10$\nr
$      3$ &  $  80$  & $8, 4, 12;  1, 73$ &    $64540, 880, 12$\cr
$      5$ &  $  84$  & $2, 0, 42;  1, 30$ &    $37802, 1260, 42$\nr
$     -2$ &  $  84$  & $4, 2, 22; 1, 141$ &    $437950, 3104, 22$\nr
$     -2$ &  $  84$  & $6, 0, 14;  1, 90$ &    $113406, 1260, 14$\nr
$      5$ &  $  84$  & $10, 4, 10;  1, 58$ &    $34114, 584, 10$\nr
$      2$ &  $  87$  & $2, 1, 44; 1, 137$ &    $826112, 6029, 44$\cr
$     -5$ &  $  87$  & $4, 1, 22; 1, 145$ &    $462844, 3191, 22$\nr
$      2$ &  $  87$  & $4, 1, 22; 1, 274$ &    $1652224, 6029, 22$\nr
$      2$ &  $  87$  & $6, 3, 16;  1, 84$ &    $113406, 1347, 16$\nr
$      2$ &  $  87$  & $8, 3, 12; 1, 112$ &    $151208, 1347, 12$\nr
$     -5$ &  $  87$  & $8, 3, 12; 1, 118$ &    $167804, 1419, 12$
\endtable

\centerline{LEVEL   523    }
\smallskip
\begintable
\   Coeff \ &\  Det   \ &\    Reduced form   \ &\   Unreduced form\crthick
$     -3$ &  $   3$  & $ 2, 1, 2;  1, 60$ &    $7322, 121, 2$\nr
$     -2$ &  $   8$  & $ 2, 0, 4;  1, 28$ &    $3138, 112, 4$\nr
$     -1$ &  $  12$  & $ 2, 0, 6; 1, 134$ &    $107738, 804, 6$\nr
$      3$ &  $  12$  & $ 4, 2, 4;  1, 60$ &    $14644, 242, 4$\nr
$      1$ &  $  20$  & $2, 0, 10; 1, 124$ &    $153762, 1240, 10$\cr
$      1$ &  $  20$  & $ 4, 2, 6;  1, 32$ &    $6276, 194, 6$\nr
$      2$ &  $  27$  & $2, 1, 14; 1, 250$ &    $875502, 3501, 14$\nr
$      3$ &  $  27$  & $ 6, 3, 6;  1, 60$ &    $21966, 363, 6$\nr
$     -2$ &  $  32$  & $2, 0, 16;  1, 14$ &    $3138, 224, 16$\nr
$      4$ &  $  32$  & $ 4, 0, 8;  1, 28$ &    $6276, 224, 8$\cr
$     -2$ &  $  32$  & $ 6, 2, 6;  1, 37$ &    $8368, 224, 6$\nr
$     -2$ &  $  35$  & $2, 1, 18; 1, 228$ &    $936170, 4105, 18$\nr
$     -2$ &  $  35$  & $ 6, 1, 6;  1, 13$ &    $1046, 79, 6$\nr
$      3$ &  $  39$  & $2, 1, 20;  1, 25$ &    $12552, 501, 20$\nr
$      3$ &  $  39$  & $4, 1, 10;  1, 50$ &    $25104, 501, 10$\cr
$      3$ &  $  39$  & $4, 1, 10; 1, 159$ &    $253132, 1591, 10$\nr
$     -7$ &  $  39$  & $ 6, 3, 8; 1, 193$ &    $299156, 1547, 8$\nr
$      7$ &  $  47$  & $2, 1, 24;  1, 80$ &    $153762, 1921, 24$\nr
$     -3$ &  $  47$  & $4, 1, 12; 1, 160$ &    $307524, 1921, 12$\nr
$     -3$ &  $  47$  & $4, 1, 12; 1, 450$ &    $2430904, 5401, 12$\cr
$     -3$ &  $  47$  & $ 6, 1, 8; 1, 152$ &    $185142, 1217, 8$\nr
$     -3$ &  $  47$  & $ 6, 1, 8; 1, 240$ &    $461286, 1921, 8$\nr
$      0$ &  $  48$  & $2, 0, 24;  1, 67$ &    $107738, 1608, 24$\nr
$     -2$ &  $  48$  & $4, 0, 12; 1, 134$ &    $215476, 1608, 12$\nr
$      0$ &  $  48$  & $ 6, 0, 8; 1, 201$ &    $323214, 1608, 8$\cr
$      6$ &  $  48$  & $ 8, 4, 8;  1, 60$ &    $29288, 484, 8$
\endtable

\begintable
\   Coeff \ &\  Det   \ &\    Reduced form   \ &\   Unreduced form\crthick
$      5$ &  $  51$  & $2, 1, 26;  1, 47$ &    $57530, 1223, 26$\nr
$     -5$ &  $  51$  & $6, 3, 10; 1, 122$ &    $149578, 1223, 10$\nr
$      4$ &  $  55$  & $2, 1, 28;  1, 22$ &    $13598, 617, 28$\nr
$      4$ &  $  55$  & $4, 1, 14;  1, 44$ &    $27196, 617, 14$\cr
$     -6$ &  $  55$  & $4, 1, 14; 1, 180$ &    $453964, 2521, 14$\nr
$     -6$ &  $  55$  & $ 8, 3, 8; 1, 184$ &    $271960, 1475, 8$\nr
$      0$ &  $  56$  & $2, 0, 28;  1, 41$ &    $47070, 1148, 28$\nr
$      0$ &  $  56$  & $4, 0, 14;  1, 82$ &    $94140, 1148, 14$\nr
$      0$ &  $  56$  & $6, 2, 10;  1, 10$ &    $1046, 102, 10$\cr
$      0$ &  $  56$  & $6, 2, 10; 1, 408$ &    $1666278, 4082, 10$\nr
$     -4$ &  $  59$  & $2, 1, 30; 1, 100$ &    $300202, 3001, 30$\nr
$     -4$ &  $  59$  & $6, 1, 10; 1, 300$ &    $900606, 3001, 10$\nr
$      6$ &  $  59$  & $6, 1, 10; 1, 432$ &    $1867110, 4321, 10$\nr
$      5$ &  $  67$  & $2, 1, 34; 1, 214$ &    $1557494, 7277, 34$\cr
$      0$ &  $  72$  & $2, 0, 36; 1, 165$ &    $980102, 5940, 36$\nr
$      0$ &  $  72$  & $4, 0, 18; 1, 193$ &    $670486, 3474, 18$\nr
$      4$ &  $  72$  & $6, 0, 12;  1, 28$ &    $9414, 336, 12$\nr
$     -3$ &  $  75$  & $2, 1, 38; 1, 126$ &    $603542, 4789, 38$\nr
$     -3$ &  $  75$  & $6, 3, 14;  1, 43$ &    $26150, 605, 14$\cr
$     12$ &  $  75$  & $10, 5, 10;  1, 60$ &    $36610, 605, 10$\nr
$     -8$ &  $  79$  & $2, 1, 40; 1, 324$ &    $4199690, 12961, 40$\nr
$      2$ &  $  79$  & $4, 1, 20; 1, 125$ &    $312754, 2501, 20$\nr
$     -8$ &  $  79$  & $4, 1, 20; 1, 241$ &    $1162106, 4821, 20$\nr
$     12$ &  $  79$  & $8, 1, 10; 1, 250$ &    $625508, 2501, 10$\cr
$      2$ &  $  79$  & $8, 1, 10; 1, 482$ &    $2324212, 4821, 10$
\endtable

\newpage
\centerline{LEVEL   587    minus space}
\smallskip
\begintable
\   Coeff \ &\  Det   \ &\    Reduced form   \ &\   Unreduced form\crthick
$      0$ &  $   8$  & $ 2, 0, 4; 1, 190$ &    $144402, 760, 4$\nr
$      0$ &  $  11$  & $ 2, 1, 6; 1, 485$ &    $1412322, 2911, 6$\nr
$      1$ &  $  15$  & $ 2, 1, 8;  1, 17$ &    $2348, 137, 8$\nr
$     -1$ &  $  15$  & $ 4, 1, 4;  1, 34$ &    $4696, 137, 4$\nr
$      0$ &  $  19$  & $2, 1, 10;  1, 64$ &    $41090, 641, 10$\cr
$      1$ &  $  20$  & $2, 0, 10;  1, 86$ &    $73962, 860, 10$\nr
$     -1$ &  $  20$  & $ 4, 2, 6;  1, 52$ &    $16436, 314, 6$\nr
$     -1$ &  $  23$  & $2, 1, 12;  1, 67$ &    $54004, 805, 12$\nr
$      0$ &  $  23$  & $ 4, 1, 6; 1, 134$ &    $108008, 805, 6$\nr
$      1$ &  $  23$  & $ 4, 1, 6; 1, 257$ &    $396812, 1543, 6$\cr
$     -1$ &  $  24$  & $2, 0, 12; 1, 275$ &    $907502, 3300, 12$\nr
$      1$ &  $  24$  & $ 4, 0, 6;  1, 37$ &    $8218, 222, 6$\nr
$     -1$ &  $  32$  & $2, 0, 16;  1, 95$ &    $144402, 1520, 16$\nr
$      0$ &  $  32$  & $ 4, 0, 8; 1, 190$ &    $288804, 1520, 8$\nr
$      1$ &  $  32$  & $ 6, 2, 6; 1, 253$ &    $385072, 1520, 6$\cr
$     -1$ &  $  35$  & $2, 1, 18;  1, 18$ &    $5870, 325, 18$\nr
$      1$ &  $  35$  & $ 6, 1, 6;  1, 54$ &    $17610, 325, 6$\nr
$      0$ &  $  39$  & $2, 1, 20;  1, 71$ &    $100964, 1421, 20$\nr
$      1$ &  $  39$  & $4, 1, 10; 1, 142$ &    $201928, 1421, 10$\nr
$      0$ &  $  39$  & $4, 1, 10; 1, 210$ &    $441424, 2101, 10$\cr
$     -1$ &  $  39$  & $ 6, 3, 8; 1, 324$ &    $841758, 2595, 8$\nr
$      0$ &  $  44$  & $2, 0, 22; 1, 269$ &    $1591944, 5918, 22$\nr
$      0$ &  $  44$  & $4, 2, 12; 1, 485$ &    $2824644, 5822, 12$\nr
$      1$ &  $  44$  & $ 6, 2, 8; 1, 434$ &    $1508590, 3474, 8$\nr
$     -1$ &  $  44$  & $ 6, 2, 8; 1, 446$ &    $1593118, 3570, 8$

\endtable

\begintable
\   Coeff \ &\  Det   \ &\    Reduced form   \ &\   Unreduced form\crthick
$     -1$ &  $  52$  & $2, 0, 26; 1, 266$ &    $1839658, 6916, 26$\nr
$      1$ &  $  52$  & $4, 2, 14;   1, 9$ &    $1174, 128, 14$ \nr
$      1$ &  $  56$  & $2, 0, 28; 1, 155$ &    $672702, 4340, 28$\nr
$      0$ &  $  56$  & $4, 0, 14; 1, 277$ &    $1074210, 3878, 14$\nr
$     -1$ &  $  56$  & $6, 2, 10; 1, 199$ &    $396812, 1992, 10$\cr
$      0$ &  $  56$  & $6, 2, 10; 1, 505$ &    $2552276, 5052, 10$\nr
$      0$ &  $  60$  & $2, 0, 30;  1, 30$ &    $27002, 900, 30$\nr
$     -1$ &  $  60$  & $4, 2, 16;  1, 17$ &    $4696, 274, 16$\nr
$      0$ &  $  60$  & $6, 0, 10;  1, 90$ &    $81006, 900, 10$\nr
$      1$ &  $  60$  & $ 8, 2, 8;  1, 34$ &    $9392, 274, 8$\cr
$      1$ &  $  71$  & $2, 1, 36; 1, 224$ &    $1806786, 8065, 36$\nr
$     -1$ &  $  71$  & $4, 1, 18; 1, 448$ &    $3613572, 8065, 18$\nr
$     -1$ &  $  71$  & $4, 1, 18; 1, 465$ &    $3892984, 8371, 18$\nr
$     -1$ &  $  71$  & $6, 1, 12;  1, 85$ &    $86876, 1021, 12$\nr
$      1$ &  $  71$  & $6, 1, 12; 1, 404$ &    $1959406, 4849, 12$\cr
$      1$ &  $  71$  & $8, 3, 10;  1, 15$ &    $2348, 153, 10$\nr
$      0$ &  $  71$  & $8, 3, 10; 1, 454$ &    $2063892, 4543, 10$\nr
$      1$ &  $  72$  & $2, 0, 36; 1, 259$ &    $2414918, 9324, 36$\nr
$     -1$ &  $  72$  & $4, 0, 18;  1, 69$ &    $85702, 1242, 18$\nr
$      0$ &  $  72$  & $6, 0, 12; 1, 190$ &    $433206, 2280, 12$\cr
$      0$ &  $  76$  & $2, 0, 38; 1, 250$ &    $2375002, 9500, 38$\nr
$      0$ &  $  76$  & $4, 2, 20;  1, 64$ &    $82180, 1282, 20$\nr
$      0$ &  $  76$  & $8, 2, 10; 1, 128$ &    $164360, 1282, 10$\nr
$      0$ &  $  76$  & $8, 2, 10; 1, 576$ &    $3320072, 5762, 10$
\endtable

\newpage
\centerline{LEVEL   587    plus space}
\smallskip
\begintable
\   Coeff \ &\  Det   \ &\    Reduced form   \ &\   Unreduced form\crthick
$      2$ &  $   8$  & $ 2, 0, 4; 1, 190$ &    $144402, 760, 4$\nr
$      2$ &  $  11$  & $ 2, 1, 6; 1, 485$ &    $1412322, 2911, 6$\nr
$     -1$ &  $  15$  & $ 2, 1, 8;  1, 17$ &    $2348, 137, 8$\nr
$     -1$ &  $  15$  & $ 4, 1, 4;  1, 34$ &    $4696, 137, 4$\nr
$      2$ &  $  19$  & $2, 1, 10;  1, 64$ &    $41090, 641, 10$\cr
$      1$ &  $  20$  & $2, 0, 10;  1, 86$ &    $73962, 860, 10$\nr
$      1$ &  $  20$  & $ 4, 2, 6;  1, 52$ &    $16436, 314, 6$\nr
$     -3$ &  $  23$  & $2, 1, 12;  1, 67$ &    $54004, 805, 12$\nr
$      8$ &  $  23$  & $ 4, 1, 6; 1, 134$ &    $108008, 805, 6$\nr
$     -3$ &  $  23$  & $ 4, 1, 6; 1, 257$ &    $396812, 1543, 6$\cr
$      1$ &  $  24$  & $2, 0, 12; 1, 275$ &    $907502, 3300, 12$\nr
$      1$ &  $  24$  & $ 4, 0, 6;  1, 37$ &    $8218, 222, 6$\nr
$      3$ &  $  32$  & $2, 0, 16;  1, 95$ &    $144402, 1520, 16$\nr
$     -4$ &  $  32$  & $ 4, 0, 8; 1, 190$ &    $288804, 1520, 8$\nr
$      3$ &  $  32$  & $ 6, 2, 6; 1, 253$ &    $385072, 1520, 6$\cr
$      1$ &  $  35$  & $2, 1, 18;  1, 18$ &    $5870, 325, 18$\nr
$      1$ &  $  35$  & $ 6, 1, 6;  1, 54$ &    $17610, 325, 6$\nr
$     -6$ &  $  39$  & $2, 1, 20;  1, 71$ &    $100964, 1421, 20$\nr
$      5$ &  $  39$  & $4, 1, 10; 1, 142$ &    $201928, 1421, 10$\nr
$     -6$ &  $  39$  & $4, 1, 10; 1, 210$ &    $441424, 2101, 10$\cr
$      5$ &  $  39$  & $ 6, 3, 8; 1, 324$ &    $841758, 2595, 8$\nr
$     -6$ &  $  44$  & $2, 0, 22; 1, 269$ &    $1591944, 5918, 22$\nr
$     -2$ &  $  44$  & $4, 2, 12; 1, 485$ &    $2824644, 5822, 12$\nr
$      5$ &  $  44$  & $ 6, 2, 8; 1, 434$ &    $1508590, 3474, 8$\nr
$      5$ &  $  44$  & $ 6, 2, 8; 1, 446$ &    $1593118, 3570, 8$
\endtable

\begintable
\   Coeff \ &\  Det   \ &\    Reduced form   \ &\   Unreduced form\crthick
$      1$ &  $  52$  & $2, 0, 26; 1, 266$ &    $1839658, 6916, 26$\nr
$      1$ &  $  52$  & $4, 2, 14;   1, 9$ &    $1174, 128, 14$\nr
$      7$ &  $  56$  & $2, 0, 28; 1, 155$ &    $672702, 4340, 28$\nr
$     -4$ &  $  56$  & $4, 0, 14; 1, 277$ &    $1074210, 3878, 14$\nr
$      7$ &  $  56$  & $6, 2, 10; 1, 199$ &    $396812, 1992, 10$\cr
$     -4$ &  $  56$  & $6, 2, 10; 1, 505$ &    $2552276, 5052, 10$\nr
$     -6$ &  $  60$  & $2, 0, 30;  1, 30$ &    $27002, 900, 30$\nr
$      3$ &  $  60$  & $4, 2, 16;  1, 17$ &    $4696, 274, 16$\nr
$     -6$ &  $  60$  & $6, 0, 10;  1, 90$ &    $81006, 900, 10$\nr
$      3$ &  $  60$  & $ 8, 2, 8;  1, 34$ &    $9392, 274, 8$\cr
$     -1$ &  $  71$  & $2, 1, 36; 1, 224$ &    $1806786, 8065, 36$\nr
$     -1$ &  $  71$  & $4, 1, 18; 1, 448$ &    $3613572, 8065, 18$\nr
$     -1$ &  $  71$  & $4, 1, 18; 1, 465$ &    $3892984, 8371, 18$\nr
$     -1$ &  $  71$  & $6, 1, 12;  1, 85$ &    $86876, 1021, 12$\nr
$     -1$ &  $  71$  & $6, 1, 12; 1, 404$ &    $1959406, 4849, 12$\cr
$     -1$ &  $  71$  & $8, 3, 10;  1, 15$ &    $2348, 153, 10$\nr
$     10$ &  $  71$  & $8, 3, 10; 1, 454$ &    $2063892, 4543, 10$\nr
$      1$ &  $  72$  & $2, 0, 36; 1, 259$ &    $2414918, 9324, 36$\nr
$      1$ &  $  72$  & $4, 0, 18;  1, 69$ &    $85702, 1242, 18$\nr
$     -4$ &  $  72$  & $6, 0, 12; 1, 190$ &    $433206, 2280, 12$\cr
$     -6$ &  $  76$  & $2, 0, 38; 1, 250$ &    $2375002, 9500, 38$\nr
$     -2$ &  $  76$  & $4, 2, 20;  1, 64$ &    $82180, 1282, 20$\nr
$     16$ &  $  76$  & $8, 2, 10; 1, 128$ &    $164360, 1282, 10$\nr
$     -6$ &  $  76$  & $8, 2, 10; 1, 576$ &    $3320072, 5762, 10$
\endtable

\enddocument

\enddocument

\enddocument

\enddocument

\enddocument

\enddocument

\enddocument

\enddocument

\enddocument

\enddocument

\enddocument

\enddocument

\enddocument